%% file: stat-svd0.tex
\documentclass[11pt]{article}
\usepackage{amsfonts}
\usepackage{amsgen,amsmath,amstext,amsbsy,amsopn,amssymb}
\usepackage{subfigure}

\usepackage[dvips]{graphicx}
\usepackage{comment}
\usepackage[shortlabels]{enumitem}
\usepackage{natbib}
\usepackage{booktabs}
\usepackage{blkarray}
\usepackage{algorithm}
\usepackage{algpseudocode}
\usepackage{url}
\usepackage{longtable}
\usepackage{stmaryrd}
\usepackage{mwe}

\usepackage[usenames]{color}

 \allowdisplaybreaks

\newtheorem{Theorem}{Theorem}

\newtheorem{Lemma}{Lemma}
\newtheorem{Remark}{Remark}

\newtheorem{Proposition}{Proposition}

\newcommand{\subscript}[2]{$#1 _ #2$}

\newcommand{\cU}{\mathcal{U}}
\newcommand{\cV}{\mathcal{V}}

\newcommand{\cI}{{\mathcal{I}}}
\newcommand{\cJ}{{\mathcal{J}}}
\renewcommand{\cU}{{\mathcal{U}}}
\renewcommand{\cV}{{\mathcal{V}}}
\newcommand{\cA}{{\mathcal{A}}}
\newcommand{\cB}{{\mathcal{B}}}

\newcommand{\br}{{\boldsymbol{r}}}
\newcommand{\bp}{{\boldsymbol{p}}}
\newcommand{\bs}{{\boldsymbol{s}}}

\newcommand{\X}{{\mathbf{X}}}

\renewcommand{\S}{{\mathbf{S}}}

\newcommand{\W}{{\mathbf{W}}}

\newcommand{\Y}{{\mathbf{Y}}}
\newcommand{\Z}{{\mathbf{Z}}}

\newcommand{\Proj}{{\mathbb{P}}}

\newcommand{\rank}{{\rm rank}}

\newcommand{\diag}{{\rm diag}}

\newcommand{\SVD}{{\rm SVD}}
\newcommand{\SSVD}{{\rm SSVD}}

\newcommand{\argmin}{\mathop{\rm arg\min}}
\newcommand{\argmax}{\mathop{\rm arg\max}}

\newcommand{\supp}{{\rm supp}}

\allowdisplaybreaks

\begin{document}

\title{Optimal Sparse Singular Value Decomposition for High-dimensional High-order Data}
\author{Anru Zhang ~ and ~ Rungang Han\footnote{Anru Zhang is Assistant Professor, Department of Statistics, University of Wisconsin-Madison, Madison, WI 53706, E-mail: anruzhang@stat.wisc.edu; Rungang Han is PhD student, Department of Statistics, University of Wisconsin-Madison, Madison, WI 53706, E-mail: rhan32@wisc.edu. The research of Anru Zhang is supported in part by NSF Grant DMS-1811868.}\\
University of Wisconsin-Madison}
\maketitle

\bigskip

\begin{abstract}
	In this article, we consider the sparse tensor singular value decomposition, which aims for dimension reduction on high-dimensional high-order data with certain sparsity structure. A method named \underline{s}parse \underline{t}ensor \underline{a}lternating \underline{t}hresholding for \underline{s}ingular \underline{v}alue \underline{d}ecomposition (STAT-SVD) is proposed. The proposed procedure features a novel double projection \& thresholding scheme, which provides a sharp criterion for thresholding in each iteration. Compared with regular tensor SVD model, STAT-SVD permits more robust estimation under weaker assumptions. Both the upper and lower bounds for estimation accuracy are developed. The proposed procedure is shown to be minimax rate-optimal in a general class of situations. Simulation studies show that STAT-SVD performs well under a variety of configurations. We also illustrate the merits of the proposed procedure on a longitudinal tensor dataset on European country mortality rates.
\end{abstract}


\input{stat-svd1-intro}

\input{stat-svd2-procedure}

\input{stat-svd3-theory}

\input{stat-svd4-simu}

\input{stat-svd5-real}

\input{stat-svd6-proof}

\section*{Acknowledgment}

The authors would like to thank the Editor, Associate Editor, and anonymous referees for their constructive comments, which have greatly helped to improve the presentation of the paper.

\bibliographystyle{apa}
\bibliography{reference}

\newpage

\appendix

\setcounter{page}{1}
\setcounter{section}{0}

\begin{center}
	{\LARGE Supplement to ``Optimal Sparse Singular Value Decomposition for}
	
	\bigskip	
	{\LARGE  High-dimensional High-order Data"	
		\footnote{Anru Zhang is Assistant Professor, Department of Statistics, University of Wisconsin-Madison, Madison, WI 53706, E-mail: anruzhang@stat.wisc.edu; Rungang Han is PhD student, Department of Statistics, University of Wisconsin-Madison, Madison, WI 53706, E-mail: rhan32@wisc.edu. The research of Anru Zhang is supported in part by NSF Grant DMS-1811868.}
	}

	\bigskip\medskip
	{Anru Zhang ~ and ~ Rungang Han}
\end{center}

\input{stat-svd7-appendix}

\end{document}

%% file: stat-svd1-intro.tex
\section{Introduction}\label{sec:intro}

High-dimensional high-order data, i.e., values arranged in large-scale tensors along three or more directions, commonly occur in a broad range of applications due to revolutionary developments in science and technology. These data possess distinct characteristics compared with the traditional low-dimensional or low-order data and pose unprecedented challenges to various communities, including statistics, machine learning, applied mathematics, and electrical engineering. To better summarize, visualize, and analyze high-dimensional high-order data, a sufficient dimension reduction often becomes the crucial first step. Therefore, how to effectively exploit the low-rank structure from high-dimensional high-order observations is often an important task. 

To this end, the framework of tensor SVD (or tensor PCA) has been introduced and extensively studied recently \citep{allen2012regularized,richard2014statistical,anandkumar2016homotopy,zhang2017tensor,liu2017characterizing,wang2017tensor}. Suppose one is interested in an order-$d$ low-rank tensor of dimension $p_1\times \cdots \times p_d$, which is observed with entry-wise additive noise $\Z$ as $\Y = \X + \Z$. Assume the fixed tensor $\X$ is low-rank in the sense that the fibers of $\X$ along different directions (i.e., the counterpart of rows and columns for matrix) all lie in low-dimensional subspaces, say $U_1,\ldots, U_d$. The goal of tensor SVD is to estimate the loadings $U_1,\ldots, U_d$ and underlying low-rank tensor $\X$. Under the regular tensor SVD setting, several practical methods have been introduced and studied, including High-order SVD (HOSVD) \citep{de2000multilinear}, High-order Orthogonal Iteration (HOOI) \citep{de2000best}, sum-of-square scheme \citep{hopkins2015tensor}, homotopy or continuation method \citep{anandkumar2016homotopy}. Using lower bound arguments, \cite{richard2014statistical} and \cite{zhang2017tensor} showed that the signal-noise-ratio (SNR) $\geq C\max\{p_1, p_2, p_3\}^{1/2}$ is required to ensure that consistent estimation is statistically possible; and SNR $\geq C\max\{p_1,p_2,p_3\}^{3/4}$ may be further necessary for computationally efficient methods. However in many applications, such conditions, i.e., SNR is no less than a polynomial of dimension $p$, are often too restrictive to satisfy.

Moreover, in many applications, the leading singular/eigenvectors of the high-dimensional high-order data may satisfy intrinsic structural assumptions along certain ways. We have seen the need for singular value decomposition in a number of modern tensor data applications, where sparsity plays an essential role. 
For example, in high-order longitudinal study, since observations often come as multivariate functions of time, (e.g. country-wise fertility and death rates by the calendar year and age \citep{wilmoth2016human}), the leading singular vectors along the mode of calendar year or age are expected to be smooth, and therefore becomes sparse after differential transformation; 
in Electroencephalogram (EEG) data analysis, the brain electrical activities are measured and stored as multi-way data with three or more modes representing channels, time, and patients, etc. It is often believed that some parts of brain region are more active and vary through time smoothly, then the leading singular vectors may be sparse along the channel mode and smooth along time mode \citep{miwakeichi2004decomposing}; in imaging ensemble analysis, facial images are often stored as high-dimensional high-order tensors. To sufficiently reduce the dimension, one looks for low-dimensional subspaces that can best explain the possibly sparse facial features and suppress illumination effects such as shadows and highlights \citep{vasilescu2003multilinear}. How to incorporate these structural assumptions wisely to improve the performance of subsequent statistical analyses is crucial for singular value decomposition in tensor data analysis. Such a problem, however, has not been well studied or understood in previous literature. 

In this article, we aim to fill this gap by developing methodology and theory for \emph{sparse tensor SVD}.
In addition to the regular tensor SVD model, we assume that the underlying low-rank structure $\X$ satisfies some sparsity constraints. As mentioned above, the data are not necessarily sparse along all modes in practice (for example, it is not reasonable to assume sparsity for the patient mode in EEG data or subject mode in high-order longitudinal data). To allow more flexibility, we suppose there exists a subset $J_s \subseteq \{1,\ldots, d\}$ such that part of the loadings $\{U_k: k\in J_s\}$ contains certain row-wise sparsity structures. The detailed formulation of the sparse tensor SVD model is introduced in Section \ref{sec:formulation}.

To better illustrate the nature and difficulty of sparse tensor SVD problem, it is also helpful to discuss its order-2 counterpart, matrix sparse singular value decomposition, for comparison. The framework of matrix sparse singular value decomposition, which focuses on extracting simultaneously sparse and low-rank matrix structure from high-dimensional matrix data, has been introduced and extensively studied during the past decade (see \cite{lee2010biclustering,yang2014sparse,yang2016rate} and the references therein). In addition, sparse principal component analysis, a closely connected topic, has also been considered in \cite{zou2006sparse,shen2008sparse,johnstone2009consistency,cai2013sparse}. In contrast, sparse tensor SVD is much more involved and difficult than sparse matrix SVD and regular tensor SVD in many aspects. First, classical methods for matrix data are often not directly applicable to high-order data. Many previous works approach the tensor problem by vectorizing or matricizing high-order data (or intuitively speaking, “stretching” the data cubes into matrices or vectors) so that high-order problems are transformed into vector or matrix ones. However, since high-order structures can get lost in the process of simple vectorizing or matricizing, one may only obtain sub-optimal results in the subsequent analyses. Second, some straightforward extensions from sparse matrix SVD methods, such as sparse HOSVD, sparse HOOI, or a single projection \& thresholding scheme, does not perform optimally in general.
Third, as pointed out by the seminal work of \cite{hillar2013most}, many basic concepts or methods for matrix data cannot be directly generalized to the high-order ones. Naive extensions of concepts such as operator norm, singular values, and eigenvalues are mathematically possible but computationally NP-hard.

To overcome these difficulties, we propose a procedure named \emph{\underline{S}parse \underline{T}ensor \underline{A}lternating \underline{T}runcation for \underline{S}ingular \underline{V}alue \underline{D}ecomposition} (STAT-SVD) for sparse tensor SVD in this paper. The method consists of two steps: (i) a thresholded spectral initialization and (ii) an iterative alternating updating scheme. One crucial part of the procedure is a novel \emph{double projection \& thresholding} scheme, which provides a sharp criterion for thresholding in each iteration. Since each step of STAT-SVD only involves basic matrix and tensor operations, such as matricization, multiplication, matrix SVD, and thresholding, the proposed procedure can be implemented efficiently.

We study both the theoretical and numerical properties of the proposed procedure. We prove by an upper bound argument that the STAT-SVD estimator can recover the low-rank structures accurately. A lower bound is further developed to show that the proposed estimator is rate optimal for a general class of simultaneously sparse and low-rank tensors. To the best of our knowledge, we are among the first ones to study the method and theory for sparse tensor SVD with matching upper and lower bound results. 
The numerical results show that the STAT-SVD outperforms other more naive methods, such as regular HOOI, HOSVD, sparse HOOI, and sparse HOSVD, by achieving significantly smaller estimation errors within much shorter running time. We also illustrate the merit of STAT-SVD in the analyses of high-order mortality rate data.

The rest of this article is organized as follows. After a brief introduction of the notations and preliminaries, we formally introduce the sparse tensor SVD model in Section \ref{sec:formulation}. Then we propose the methodology for sparse tensor SVD in Section \ref{sec:method}. The theoretical properties of the proposed procedure are developed in Section \ref{sec:theory}. The data-driven hyperparameter selection is discussed in Section \ref{sec:parameter}. Numerical performance of the proposed methods is studied through both simulation studies and the real data analyses on high-order mortality rate dataset in Section \ref{sec:simu}. 
Finally, further discussions, proofs of the technical results, and supporting theoretical tools are postponed to the supplementary materials.

\section{Problem Formulation}\label{sec:formulation}

\subsection{Notations and Preliminaries}

We start this section with notations and preliminaries that will be used throughout the paper. The lowercase letters, e.g. $x, y, u, v$, are used to denote scalars or vectors. For any $a, b\in \mathbb{R}$, let $a\wedge b$ and $a\vee b$ be the minimum and maximum of $a$ and $b$, respectively. For convenience, we denote $\bp = (p_1,\ldots, p_d)$, $\br = (r_1,\ldots, r_d)$, $\bs = (s_1,\ldots, s_d)$, and $p = p_1\cdots p_d$, $s = s_1\cdots s_d, r = r_1\cdots r_d$. For any $1\leq k \leq p$, we also note $p_{-k} = \prod_{l \neq k} p_l, s_{-k} = \prod_{l \neq k}s_l, r_{-k} = \prod_{l\neq k} r_l$. We use $C, c$ and the variations to denote generic constants, whose actual values may change from line to line.

The uppercase letters are used to note matrices. For $A\in \mathbb{R}^{m\times n}$, we assume the singular value decomposition is $A = \sum_{i} \sigma_i u_i v_i^\top$, where $\sigma_1 \geq \sigma_2 \geq \cdots \geq 0$ are the singular values in descending order. Denote $\sigma_i(A) = \sigma_i$ as the $i$-th largest singular value of $A$. Particularly, the largest and smallest non-trivial singular values: $\sigma_{\max}(A) = \sigma_1(A)$ and $\sigma_{\min}(A) = \sigma_{m \wedge n}(A)$ play important roles in our analysis. We also define $\SVD_r(A) = [u_1,\ldots, u_r]$ as the matrix comprised of the top $r$ left singular vectors of $A$. Let the collections of regular and sparse orthogonal matrices be $\mathbb{O}_{p, r} = \{U\in \mathbb{R}^{p\times r}: U^\top U = I_r\}$, $\mathbb{O}_{p, r}(s) = \{U\in \mathbb{R}^{p\times r}: U^\top U = I_r, \|U\|_0 = \sum_{i=1}^p 1_{\{U_{[i,:]}\neq 0\}}\leq s\}$. These orthogonal matrices are extensively used in later narratives.

In addition, the boldface capital letters, e.g. $\X$, $\Y$, $\Z$, are used to represent tensors of order-3 or higher. For any $\S \in \mathbb{R}^{r_1\times \cdots \times r_d}$ and $U_k\in \mathbb{R}^{p_k\times r_k}$, the mode-$k$ tensor-matrix product is defined as
$$\S \times_k U_k\in \mathbb{R}^{r_1\times \cdots \times r_{k-1}\times p_k\times r_{k+1}\times \cdots \times r_d}, \quad \left(\S \times_k U_k\right)_{i_1,\ldots, i_d} = \sum_{j_k = 1}^{r_k} \S_{i_1,\ldots, j_{k}, \ldots, i_d} U_{i_k, j_k}.$$
Multiplication along different directions is commutative invariant, i.e., $\left(\S \times_{k_1} U_{k_1}\right) \times_{k_2} U_{k_2} = \left(\S \times_{k_2} U_{k_2}\right) \times_{k_1} U_{k_1}$ for $k_1\neq k_2$. For convenience, the tensor product along all $d$ modes is applied
\begin{equation*}
\S \times_1 U_1 \times \cdots \times_d U_d := \llbracket \S; U_1,\ldots, U_d\rrbracket.
\end{equation*}
Matricization is a basic tensor operation that transforms tensors to matrices. Particularly the mode-$1$ matricization for any $\X \in \mathbb{R}^{p_1\times \cdots \times p_d}$ is defined as
\begin{equation*}
\mathcal{M}_1(\X)\in \mathbb{R}^{p_1\times p_{-1}}, \quad \text{where} \quad [\mathcal{M}_1(\X)]_{i_1, i_2+p_2(i_3-1)+\cdots + p_2\cdots p_{d-1}(i_d-1)} = \X_{i_1,\ldots, i_p}.
\end{equation*}
The general mode-$k$ matricization can be defined similarly. The readers are referred to \cite{kolda2009tensor} for a more comprehensive tutorial for tensor algebra.

We also use the R syntax to denote sub-vectors, -matrices, and -tensors. For example, $A_{[I_1, I_2]}$ is used to note the sub-matrix with row indices $I_1$ and column indices $I_2$ of $A$. In addition, for any integers $a \leq b$, we use ``$a:b$" to represent consecutive sequence $\{a,\cdots, b\}$; and ``:" alone represents the entire index set. Thus, $A_{[:,1:r]}$ represents the first $r$ columns of $A$ while $A_{[(s_1+1):p_1,:]}$ represents the $\{(s_1 + 1), \ldots, p_1\}$-th rows with indices $\{(s_1+1),\ldots, p_1\}$. Similar notations are also applied to sub-vectors and sub-tensors.

\subsection{Sparse Tensor SVD Model}

\begin{figure}
	\centering
	\includegraphics[width=0.8\linewidth]{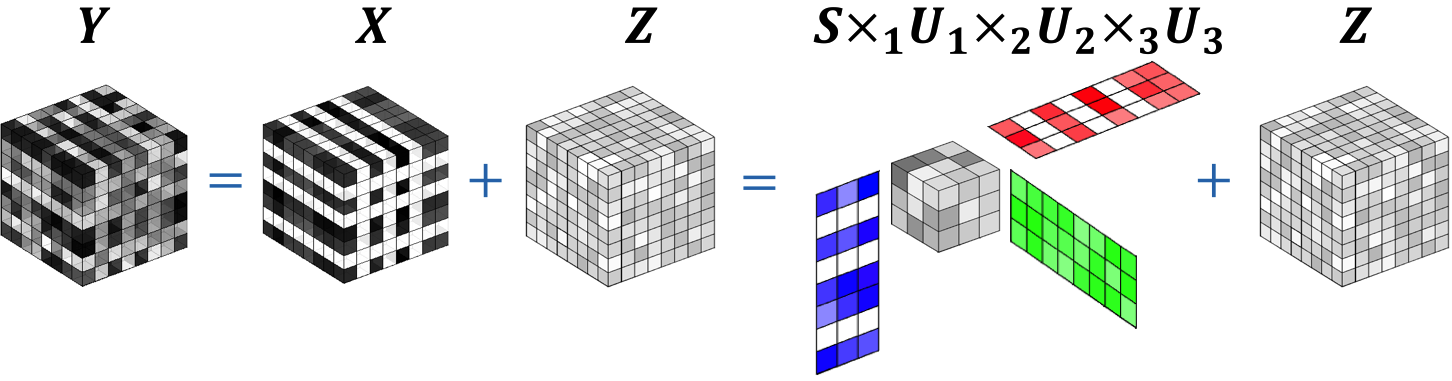}
	\caption{Illustration of sparse tensor SVD model. Here, $\X$ is sparse along Modes-1 and -3.}
	\label{fig:sparse-tensor-svd}
\end{figure}

Now we formally introduce the sparse tensor SVD model. Suppose one observes a $(p_1\times \cdots \times p_d)$-dimensional tensor $\Y$ with additive noise, say $\Y = \X + \Z$. Here, $\X$ is a low-rank tensor in the sense that all fibers along each direction lie in some low-dimensional subspace, $\Z$ consists of i.i.d. Gaussian noises with mean zero and variance $\sigma^2$. Given the connection between Tucker low-rank and Tucker decomposition \citep{tucker1966some,kolda2009tensor}, the model can be further written as
\begin{equation}\label{eq:model}
\Y = \X + \Z = \S\times_1 U_1\times\cdots \times_d U_d + \Z = \llbracket \S; U_1,\ldots, U_d \rrbracket + \Z,
\end{equation}
with $\S\in \mathbb{R}^{r_1\times \cdots \times r_d}$ and $\{U_k \in \mathbb{O}_{p_k, r_k}\}_{k=1}^d$ being the unknown core tensor and Mode-$1, \ldots, k$ loadings, respectively. Especially, $U_k$ is also the singular subspace of $\X$ along the $k$-th direction.
Recall $\mathbb{O}_{p, r}$ and $\mathbb{O}_{p, r}(s)$ are the classes of regular and sparse orthogonal columns. Given previous discussions, suppose a subset $J_s\subseteq \{1,\ldots, d\}$ is known \emph{a priori}, such that the mode-$k$ loading $U_k$ for any $k\in J_s$ is $s_k$-sparse,
\begin{equation*}
U_k \in \left\{
\begin{array}{ll}
\mathbb{O}_{p_k, r_k}(s_k), & k\in J_s;\\
\mathbb{O}_{p_k, r_k}, & k\notin J_s.
\end{array}
\right.
\end{equation*}
To be more flexible, we set $s_k = p_k$ if there is no sparsity constraint in the $k$-th Mode, i.e., $k\notin J_s$. The central goal is to estimate singular subspaces $U_1,\ldots, U_d$, and the original low-rank tensor $\X$ based on observations $\Y$. See Figure \ref{fig:sparse-tensor-svd} for a pictorial illustration of sparse tensor SVD model.
\begin{Remark}\rm
It is worth mentioning that our framework is related but distinct from the CP low-rank-based decomposition framework in literature (see \cite{de2008tensor,kolda2009tensor,anandkumar2014tensor,anandkumar2014guaranteed,sun2015provable,sun2017dynamic,wang2017three,hao2018sparse} and the references therein). Since the main focus of tensor SVD analysis is to sufficiently reduce the high-order data onto low-dimensional subspaces, the Tucker low-rank structure provides a more natural fit based on the interpretation of Tucker decomposition that we have mentioned above. Additionally, it is known that any CP rank-$r$ tensor must be of Tucker rank at most $r$, but the opposite of this statement is not generally true  \citep{kolda2009tensor}. Therefore, we mainly pursue the more general Tucker low-rank model \eqref{eq:model} rather than the CP one, with distinct methodologies and theories developed for the rest of this paper.
\end{Remark}

%% file: stat-svd2-procedure.tex
\section{The STAT-SVD Procedure}\label{sec:method}

Next, we propose the procedure for sparse tensor singular value decomposition. Note that $U_k$ in the sparse tensor SVD model \eqref{eq:model} is not only the mode-$k$ singular subspace of $\X$ but also the left singular subspace of $\mathcal{M}_k(\X)$, a straightforward idea for initialization is
\begin{equation*}
\tilde{U}_k^{(0)} = \SVD_{r_k}\left(\mathcal{M}_k(\Y)\right).
\end{equation*}
This method, originally introduced by \cite{de2000multilinear}, has been widely referred to as high-order SVD (HOSVD) and used in numerous scenarios. However, $\tilde{U}_k^{(0)}$ may fail to provide any consistent estimations or even warm starts, unless one has strong SNR $\lambda/\sigma \geq Cp^{3/2}$ \citep{zhang2017tensor}, where $\lambda = \min_{k=1,2,3}\sigma_{r_k}(\mathcal{M}_k(\X))$. An alternative idea is to apply $\ell_1$ regularized estimation to encourage sparsity \citep{allen2012sparse,allen2012regularized},
$$\hat{U}_1, \hat{U}_2, \hat{U}_3 = \argmin_{U_1,U_2,U_3} \|\Y - \S \times_1 U_1 \times_2 U_2 \times_3 U_3\|_F^2 + \lambda\sum_{k\in J_s} \|U_k\|_1.$$
Both the computation and theoretical analysis for the regularized MLE may be difficult. Instead, we propose a procedure with a novel double thresholding \& projection scheme as follows.
\begin{enumerate}[leftmargin=*]
	\item[Step 1] (Initialization: support) Let $Y_k = \mathcal{M}_k(\Y)$ be the matricization of $\Y$ for $k=1,\ldots, d$. Recall $p_{-k} = p_1\cdots p_d/p_k$, $p=p_1\cdots p_d$. For $k=1,\ldots, d$, select the index sets for all sparse modes,
	\begin{equation}\label{eq:I_k^{(0)}}
	\begin{split}
	\hat{I}_k^{(0)} = & \Bigg\{i: \left\|(Y_k)_{[i, :]}\right\|_2^2 \geq \sigma^2 \left(p_{-k} + 2\sqrt{p_{-k}\log p} + 2 \log p\right)\\
	& \quad \text{or } \max_{j} |(Y_k)_{[i,j]}| \geq 2\sigma\sqrt{\log p}\Bigg\},\quad k \in J_s.
	\end{split}
	\end{equation}
	Ideally speaking, $\hat{I}_k^{(0)}$ provides an initial estimate for the support of $\{U_k: k\in J_s\}$, which captures the significant signals of $\X$ but may miss the weak ones. For $k\notin J_s$, we select $\hat{I}_k^{(0)} = \{1,\ldots, p_k\}$.
	\item[Step 2] (Initialization: loadings) Construct
	\begin{equation*}
	\tilde{\Y} \in \mathbb{R}^{p_1\times \cdots \times p_d}, \quad \tilde{\Y}_{[i_1,\ldots, i_d]} = \left\{\begin{array}{ll}
	\Y_{[i_1,\ldots, i_d]}, & (i_1,\ldots, i_d) \in \hat{I}_1^{(0)}\otimes \cdots \otimes \hat{I}_d^{(0)};\\
	0, & \text{otherwise},
	\end{array}\right.
	\end{equation*}
	and initialize
	\begin{equation*}
	\hat{U}_k^{(0)} = \SVD_{r_k}\left(\mathcal{M}_k(\tilde{\Y})\right),\quad k=1,\ldots, d.
	\end{equation*}
	Then $\hat{U}_k^{(0)}$ provides a rough estimate for $U_k$. The initialization step is illustrated in Figure \ref{fig:procedure} (a).
\begin{figure}\label{fig:procedure}
	\begin{center}
		\subfigure[Initilaization]{\includegraphics[width=0.7\linewidth]{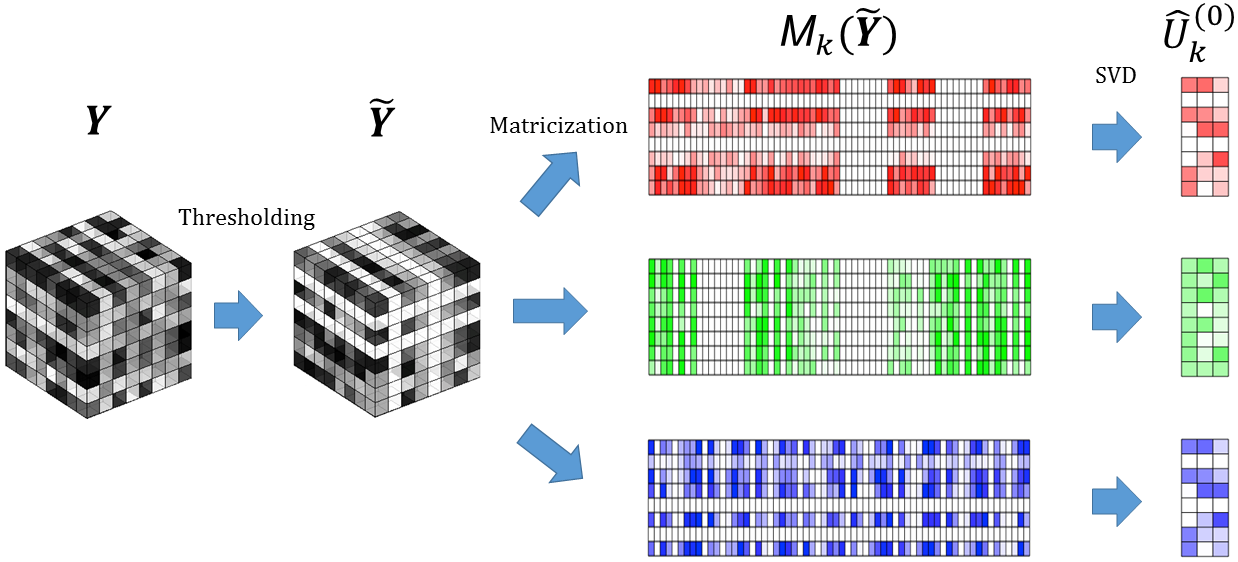}}\\\vskip.5cm
		\subfigure[Dense mode update]{\includegraphics[width=0.7\linewidth]{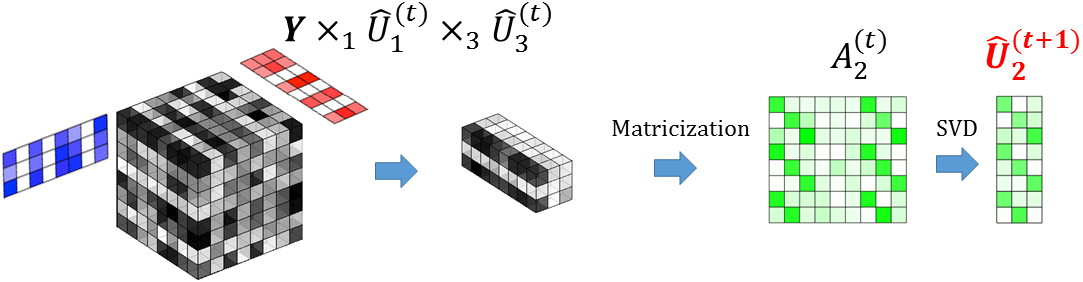}}\\\vskip.5cm
		\subfigure[Sparse mode update]{\includegraphics[width=0.7\linewidth]{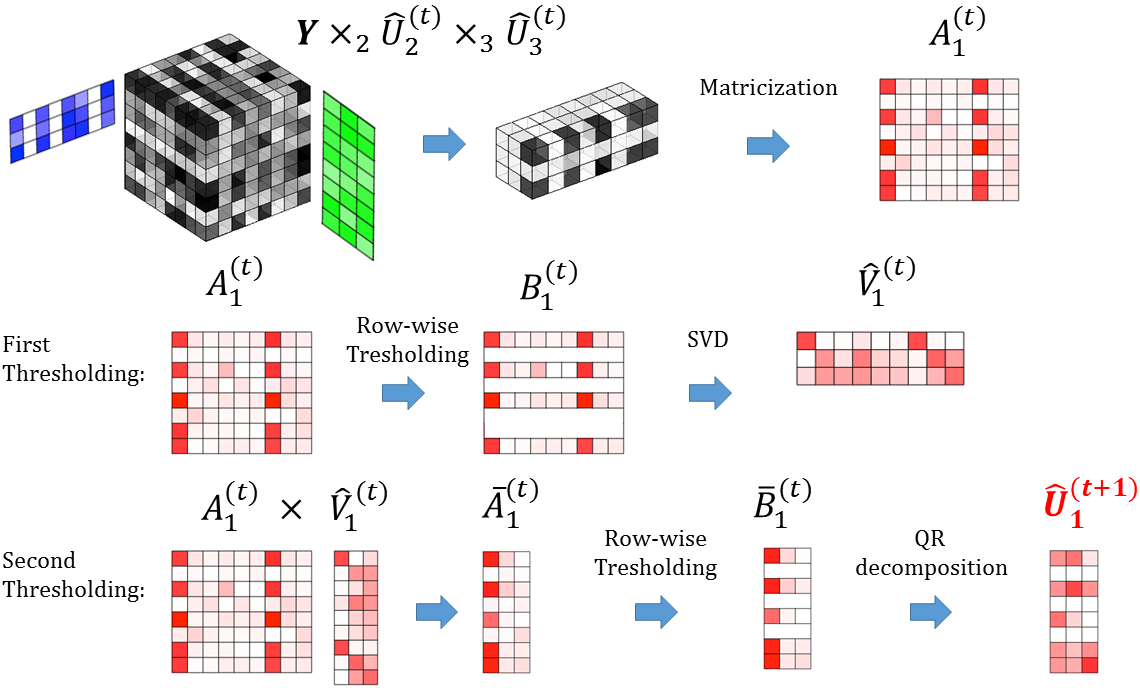}}
	\end{center}
	\caption{Illustration of STAT-SVD procedure.}
\end{figure}

	\item[Step 3] (Power Iteration and Alternating Thresholding) Provided a reasonable initialization by thresholded spectral method, we consider an iterative alternating scheme to refine the estimation. Suppose we aim to update $\hat{U}_k^{(t)} \to \hat{U}_k^{(t+1)}$ for some specific $k=1,\ldots, d$ and $t=0,1,\ldots$. The updating scheme is discussed under two scenarios.	
	\begin{itemize}[leftmargin=*]
		\item If the targeting mode, say Mode-2, is non-sparse, we can update by orthogonal projection
		\begin{equation}\label{eq:update-non-sparse-projection}
		A_2^{(t)} = \mathcal{M}_2\left(\left\llbracket \Y; (\hat{U}_1^{(t+1)})^\top, I_{p_2}, (\hat{U}_{3}^{(t)})^\top, \ldots, (\hat{U}_{d}^{(t)})^\top \right\rrbracket\right)\in \mathbb{R}^{p_2\times r_{-2}}.
		\end{equation}
		Ideally speaking, the dimension of $\Y$ is  significantly reduced from $p_1\cdots p_d$ to $p_2\times r_{-2}$, while the majority of signal $\X$ can be preserved in $A_2^{(t)}$ after such a projection. Then we update
		$$\hat{U}_2^{(t+1)} = \SVD_{r_2}\big(A_2^{(t)}\big).$$
		See Figure \ref{fig:procedure} (b) for an illustration.

		\item If the targeting mode, say Mode-1, is sparse, we apply a \emph{double projection \& thresholding scheme} for refinement (see Figure \ref{fig:procedure} (c)). We still perform projection first,
		\begin{align*}
		\text{(1st Projection)} \quad & A_1^{(t)} = \mathcal{M}_1\left(\left\llbracket \Y; I_{p_1}, (\hat{U}_2^{(t)})^\top, \ldots, (\hat{U}_{d}^{(t)})^\top \right\rrbracket\right)\in \mathbb{R}^{p_1\times r_{-1}}.
		\intertext{Then $A_1^{(t)}$ is denoised by row-wise hard thresholding,}
		\text{(1st Thresholding)}\quad & B_1^{(t)}\in \mathbb{R}^{p_1\times r_{-1}},\quad  B_{1, [i,:]}^{(t)} = A_{1, [i,:]}^{(t)}1_{\{\|A_{1,[i,:]}^{(t)}\|_2^2 \geq \eta_1\}}, \quad 1\leq i \leq p_1,
		\end{align*}
		where the thresholding level $\eta_1 = \sigma^2(r_{-1} + 2\sqrt{r_{-1}\log p}+2\log p)$ is determined by the quantile of $\|A_{1,[i,:]}^{(t)}\|_2^2$ if the $i$-th row of $A_1^{(t)}$ does not contain any signal. Since $\|A_{1, [i, :]}\|_2^2 \sim \sigma^2\chi^2_{r^2} \approx \sigma^2 r^2$ when $A_{1, [i, :]}$ are pure noise, an $\eta_1$ induced by $\|A_{1, [i, :]}\|_2^2$ can be as large as $\approx \sigma^2r_{-1}$, which may falsely kill many rows with weak signals. In order to lower the large thresholding level, we introduce the second projection and thresholding: let the leading $r_1$ right singular vectors of $B_1^{(t)}$ be $\hat{V}_1^{(t)}$, i.e., $\hat{V}_1^{(t)} = \SVD_{r_1}(B_1^{(t)\top})$, then we further reduce the dimension of $p_1$-by-$r_{-1}$ matrix $A_1^{(t)}$ to $p_1$-by-$r_1$ matrix $\bar{A}_1^{(t)}$ by performing right projection,
		\begin{align*}
		\text{(2nd Projection)}\quad & \bar{A}_1^{(t)} = A_1^{(t)}\hat{V}_k^{(t)} \in \mathbb{R}^{p_1\times r_1},
		\intertext{and apply the second thresholding}
		\text{(2nd Thresholding)}\quad & \bar{B}^{(t)}\in \mathbb{R}^{p_1\times r_1}, \quad \bar{B}^{(t)}_{1, [i, :]} = \bar{A}_{1, [i, :]} 1_{\{\|A_{1,[i,:]}^{(t)}\|_2^2 \geq \bar{\eta}_1\}},
		\end{align*}
		with much smaller thresholding value $\bar{\eta}_1 = \sigma^2\left(r_1 + 2\sqrt{r_1\log p} + 2\log p\right)$, given the reduced dimension of $\bar{A}_1^{(t)}$. As we will illustrate in theoretical analysis, the double projection \& thresholding scheme provides more accurate denoising performance. It is also noteworthy that a similar version of double thresholding appears in recent high-dimensional clustering literature \citep{jin2016influential}, although their problem, method, and theory were all different from ours. Finally, we update $\hat U_1$ by QR decomposition
		$$\hat{U}_1^{(t+1)} = {\rm QR}(\bar{B}_{1}^{(t)}).$$
	\end{itemize}

	\item[Step 4] The iteration is stopped until the maximum number of iteration is reached (i.e. $t> t_{\max}$), or convergence, i.e. the following criterion holds,
	\begin{equation*}
	\left\|\bar{B}_k^{(t)}\right\|_F^2 - \left\|\bar{B}_k^{(t-5)}\right\|_F^2 < \varepsilon_{tol}.
	\end{equation*}
	Here $\varepsilon_{tol}$ is the maximum tolerance, which can be chosen empirically. 	Finally, we propose the denoising estimator for $\X$ as
	\begin{equation*}
	\hat{\X} = \Y \times_1 (\hat{U}_1\hat{U}_1^\top) \times\cdots \times_d (\hat{U}_d\hat{U}_d^\top).
	\end{equation*}
\end{enumerate}
The pseudo-code for the proposed STAT-SVD method is provided in Algorithm \ref{al:procedure}. We particularly summarize the double projection \& thresholding scheme for the sparse mode update in Algorithm \ref{al:refinement-sparse}.
\begin{algorithm}
	\caption{Sparse tensor alternating thresholding - SVD (STAT-SVD)}
	\begin{algorithmic}[1]
		\State Input: order-$d$ tensor data $\Y\in \mathbb{R}^{\bp}$, rank $\br$, noise level $\sigma^2$, set of sparse modes $J_s$.
		\State (Initialization) Let $Y_k = \mathcal{M}_k(\Y) $. Select subset $\hat{I}_k^{(0)}$ by
		\begin{equation*}
		\hat{I}_k^{(0)} = \left\{\begin{array}{ll}\Big\{1\leq i\leq p_k: \left\|(Y_k)_{[i, :]}\right\|_2^2 \geq \sigma^2 \left(p_{-k} + 2\sqrt{p_{-k}\log p} + 2 \log p\right)\\
		\quad \quad \text{or } \max_{j} |(Y_k)_{[i,j]}| \geq 2\sigma\sqrt{\log p}\Big\},  & k\in J_s;\\
		\{1,\ldots, p_k\}, & k\notin J_s.
		\end{array}\right.
		\end{equation*}
		\State Set $t=0$. Calculate the initializations
		$$\hat{U}^{(0)}_k = \SVD_{r_k}\left(\mathcal{M}_k(\tilde{\Y})\right), \quad \text{where} \quad \tilde{\Y} =\left\{\begin{array}{ll}
		\Y_{i_1,\ldots,i_d}, & (i_1,\ldots, i_d)\in \hat{I}_1^{(0)}\otimes \cdots \otimes \hat{I}_d^{(0)};\\
		0, & \text{otherwise}.
		\end{array}\right.$$
		\While{$t < t_{\max}$ or convergence criterion not satisfied}
			\For{$k = 1,\ldots, d$}
				\If{$k\in J_s$} 				
				\State (Sparse mode update) Update $\hat{U}_k^{(t+1)}$ via Algorithm \ref{al:refinement-sparse};
				\Else
				\State (Dense mode update) Calculate
				$$A_k^{(t+1)} = \mathcal{M}_k\left(\llbracket \Y;  (\hat{U}_1^{(t+1)})^\top,  \ldots, (\hat{U}_{k-1}^{(t+1)})^\top, I_{p_k}, (\hat{U}_{k+1}^{(t)})^\top, \ldots, (\hat{U}_{d}^{(t)})^\top \rrbracket \right).$$
				\State Update as $\hat{U}_k^{(t+1)} = \SVD_{r_k}\left(A_k^{(t+1)}\right).$
				\EndIf
			\EndFor
			\State $t = t+1$.
		\EndWhile
	\end{algorithmic}\label{al:procedure}
\end{algorithm}

\begin{algorithm}
	\caption{Update Scheme in STAT-SVD -- Sparse Mode}
	\begin{algorithmic}[1]
	\State Input: $\Y$, $\hat{U}_{k+1}^{(t)}, \ldots, \hat{U}_d^{(t)}, \hat{U}_1^{(t+1)},\ldots, \hat{U}_{k-1}^{(t+1)}$, rank $r_k$.
	\State (First Projection) Calculate
	$$A_k^{(t)} = \mathcal{M}_k\left(\llbracket \Y;  (\hat{U}_1^{(t+1)})^\top,  \ldots, (\hat{U}_{k-1}^{(t+1)})^\top, I_{p_k}, (\hat{U}_{k+1}^{(t)})^\top, \ldots, (\hat{U}_{d}^{(t)})^\top \rrbracket \right).$$
	\State (First Thresholding) Perform row-wise thresholding for $A_k^{(t+1)}$:
	$$B_k^{(t+1)} \in \mathbb{R}^{p_k\times r_{-k}}, \quad B^{(t)}_{k, [i, :]} = A^{(t+1)}_{k, [i, :]} 1_{\left\{\|A^{(t+1)}_{k, [i, :]}\|_2^2 \geq \eta_k\right\}}, \quad 1\leq i \leq p_k,$$
	where $\eta_k = \sigma^2\left(r_{-k}+ 2\sqrt{r_{-k} \log p} + 2\log p\right)$.
	\State (Second Projection) Extract the leading $r_k$ right singular vectors of $B_k^{(t)}$
	$$\hat{V}_k^{(t+1)} = \SVD_{r_k}\left(B_k^{(t+1)\top}\right) \in \mathbb{O}_{p_k, r_k},$$
	then project as $\bar{A}_k^{(t+1)} = A_k^{(t+1)} \hat{V}_k^{(t+1)} \in\mathbb{R}^{p_k\times r_k}$.
	\State (Second thresholding) Perform thresholding on $\bar{B}_k^{(t+1)}$,
	$$\bar{B}_k^{(t+1)} \in \mathbb{R}^{p_k\times r_{k}}, \quad \bar{B}^{(t+1)}_{k, [i, :]} = \bar{A}^{(t+1)}_{k, [i, :]} 1_{\left\{\|A^{(t+1)}_{k, [i, :]}\|_2^2 \geq \bar{\eta}_k\right\}}, \quad 1\leq i \leq p_k,$$
	where $\bar{\eta}_k = \sigma^2\left(r_k + 2\sqrt{r_k\log p} + 2\log p\right)$.
	\State Apply QR decomposition to $\bar{B}_k^{(t+1)}$, and assign the $Q$ part to $\hat{U}_k^{(t+1)}\in \mathbb{O}_{p_k, r_k}$.
	\end{algorithmic}
	\label{al:refinement-sparse}
\end{algorithm}

%% file: stat-svd3-theory.tex
\section{Theoretical Analysis}\label{sec:theory}

In this section, we analyze the theoretical properties of the proposed procedure in the previous section. The $\sin\Theta$ distances are adopted to quantify the singular subspace estimation errors. Particularly for any $U, V \in \mathbb{O}_{p, r}$, the principal angles between $U$ and $V$ is defined as an $r$-by-$r$ diagonal matrix: $\Theta(U, V) = \diag(\cos^{-1}(\sigma_1(U^\top V)),\ldots, \cos^{-1}(\sigma_r(U^\top V)))$. Then the $\sin\Theta$ Frobenius norm $\|\sin\Theta(U, V)\|_F = \sqrt{r - \|U^\top V\|_F^2}$ can be used to characterize the distance between $U$ and $V$. The readers are referred to \cite[Lemma 1]{cai2016rate} for more discussions on the properties of $\sin\Theta$ distance.

\begin{Theorem}[Upper bound]\label{th:upper_bound}
	Suppose $\br, \sigma^2$ are known, $\log{p_1} \asymp \cdots \asymp \log{p_d}$, and for $1\leq k \leq d$, $\lambda_k = \sigma_{r_k}\left(\mathcal{M}_k(\X)\right) \geq C_0\sigma\left((s \log p)^{1/2} + \sum_{k} s_k r_k + \max_{1\leq k\leq d} r_{-k}\right).$ Then after at most $O(ds\log p)$ iterations, the proposed Algorithm \ref{al:procedure} yields the following estimation error upper bound,
	\begin{equation}\label{ineq:hat_X-X}
	\left\|\hat{\X} - \X\right\|_F^2\leq C\sigma^2 \left(\prod_{l=1}^d r_l + \sum_k s_kr_k + \sum_{k\in J_s} s_k\log{p_k} \right).
	\end{equation}
	\begin{equation}
	\begin{split}
	\left\|\sin\Theta\left(\hat{U}_k, U_k\right)\right\|^2_F \leq \left\{\begin{array}{ll}
	\left\{C\sigma^2\left(s_kr_k + s_k \log p_k\right)/\lambda_k^2\right\} \wedge r_k, &  k \in J_s,\\
	\left\{C\sigma^2 p_kr_k/\lambda_k^2\right\}\wedge r_k, & k\notin J_s,
	\end{array}\right.
	\end{split}
	\end{equation}
	with probability at least 
	$1 - \frac{C(p_1+\cdots+p_d)\log(s\log p)}{p_1\cdots p_d}$. Here $C>0$ is some uniform constant, which does not depend on $\sigma, \bp, \br, \bs$.
\end{Theorem}

\begin{Remark}\rm
	The estimation error upper bound \eqref{ineq:hat_X-X} is comprised of three terms: $\sigma^2\prod_{l=1}^d r_l$, $\sigma^2 s_kr_k$, and $\sigma^2s_k\log(p_k)$, which correspond to the estimation complexity for the core tensor $\S$, the values of loading $U_k$, and the support of $U_k$ (only for sparse modes), respectively. On the other hand, the signal strength assumption $\sigma_{r_k}(X_k)\geq C\sigma\left((s\log p)^{1/2} + \sum_k s_kr_k + \max_{1\leq k \leq d} r_{-k}\right)$ involves $p$ only in the logarithmic term. Compared with the assumption $\sigma_{r_k}(X_k) \geq \sigma\max\{p_1, p_2, p_3\}^{3/4}$ that is required in regular tensor SVD, our proposed algorithm is able to handle high-dimensional settings under much weaker conditions.
\end{Remark}

\begin{Remark}[Proof sketch for Theorem \ref{th:upper_bound}]\rm
	Since the proof of Theorem \ref{th:upper_bound} is lengthy and highly non-trivial, we briefly discuss the sketch here. First, a number of conditions are introduced as the baseline assumptions (Step 1). Under these conditions, we try to establish the upper bound for the initial estimate:
	$$\|\sin\Theta(\hat{U}_k^{(0)}, U_k)\|\leq c\quad \text{and}\quad \|(\hat{U}_{k\perp}^{(0)})^\top \mathcal{M}_k(\X)\|_F \leq C\sqrt{s\log p}$$
	for constants $0<c<1/2$ and $C>0$ (Step 2), then develop the upper bound for estimates after each iteration:
	$$ \left\|\sin\Theta(\hat{U}_k^{(t+1)}, U_k)\right\|_F \leq \frac{\sqrt{s_kr_k} + \sqrt{s_k\log p}\cdot 1_{\{k\in J_s\}}}{\lambda_k/\sigma} + \frac{C\sqrt{s\log p}}{\lambda_k/\sigma}\cdot \tau^t,$$ 
	$$\text{and}\quad \left\|(\hat{U}_{k\perp}^{(t)})^\top X_k\right\|_F \leq C\sigma\sqrt{s_kr_k} + C\sigma\sqrt{s_k \log p}\cdot 1_{\{k\in J_s\}} + C\sigma\sqrt{s\log p}\cdot \tau^t.$$
	for some constant $0<\tau\leq1/2$ (Step 3). Then after a number of iterations, the error rate sufficiently decays. One obtain final estimators $\hat{U}_k^{(t_{\max})}$ and $\hat{\X}^{(t_{\max})}$ that achieve the targeting upper bounds of estimation error (Steps 4 and 5). Finally, we use a coupling scheme to show that the introduced conditions hold with high probability (Step 6).
\end{Remark}

\begin{Remark}[Performance of Single Projection \& Truncation Scheme]	\rm
	As discussed earlier, the proposed double projection \& thresholding scheme is crucial to the performance of STAT-SVD. Such a scheme is essential in the sense that the simple single projection \& thresholding, which may be a more straightforward extension from matrix sparse SVD method, may yield sub-optimal results. To be specific, suppose $\tilde{\X}$ is the estimator with single thresholding \& projection (see Algorithm \ref{al:refinement-sparse-single-thresholding} in Appendix for detailed explanations), Theorem \ref{th:single-thresholding-projection} in the Appendix shows that there exists a low-rank tensor $\X$ that satisfies all assumptions in Theorem \ref{th:upper_bound}, but $\tilde{\X}$ yields a higher rate of convergence in the sense that
	\begin{equation*}
	\mathbb{E}\|\tilde{\X} - \X\|_F^2 \geq C\sigma^2\left(\prod_l r_l + \sum_{k} s_k r_k + \sum_{k \in J_s} s_k\log p_k + \sum_{k \in J_s} s_k\left(r_{-k}\log p_k\right)^{1/2}\right).
	\end{equation*}
\end{Remark}

Next, we study the statistical lower bound for sparse tensor SVD. Consider the following class of sparse and low-rank tensors,
\begin{equation*}
\begin{split}
\mathcal{F}_{\bp, \br, \bs}(J_s) =
\left\{\X = \llbracket \S; U_1,\ldots, U_d\rrbracket \in \mathbb{R}^{p_1\times \cdots \times p_d}: \begin{array}{l}
\rank(\X) \leq (r_1,\ldots, r_d),\\ \left\|\supp(U_k)\right\|_0 \leq s_k \text{ for } k\in J_s
\end{array}\right\}.
\end{split}
\end{equation*}
The following lower bound results hold for sparse tensor SVD.
\begin{Theorem}[Lower Bound: Subspace Estimation]\label{th:lower_bound_subspace}
Suppose $s_k \geq 3r_k$, consider the following $k$ classes of sparse and low-rank tensors with least singular value constraint on matricization,
	\begin{equation*}
	\mathcal{F}_{\bp, \br, \bs}^{(k)}(J_s, \lambda_k) = \left\{\X\in \mathcal{F}_{\bp, \br, \bs}(J_s): \sigma_{r_k}\left(\mathcal{M}_k(\X)\right)\geq \lambda_k\right\}, \quad k=1,\ldots, d.
	\end{equation*}
	Then, for any fixed $k$, we have
	\begin{equation*}
	\inf_{\hat{U}_k} \sup_{\X \in \mathcal{F}_{\bp, \br, \bs}^{(k)}(J_s, \lambda_k)} \left\|\sin\Theta\left(\hat{U}_k, U_k\right) \right\|_F^2 \geq c\left(\frac{\sigma^2s_kr_k}{\lambda_k^2} + \frac{\sigma^2s_k \log(p_k/s_k)}{\lambda_k^2} \cdot I_{\{k \in J_s\}} \right)\wedge r_k.
	\end{equation*}
\end{Theorem}

\begin{Theorem}[Lower bound: Tensor Recovery]\label{th:lower_bound_recovery}
	Suppose $s_k\geq 3r_k$ for any $k$. Consider the tensor recovery over the class of sparse and low-rank tensors $\mathcal F_{\bp, \br, \bs}(J_s)$, there exists uniform constant $c>0$ such that
	\begin{equation*}
	\inf_{\hat{\X}} \sup_{\X \in \mathcal{F}_{\bp, \br, \bs}(J_s)} \left\|\hat{\X} - \X \right\|_F^2 \geq c\sigma^2\left(\prod_{k=1}^d r_k + \sum_{k=1}^d s_kr_k + \sum_{k\in J_s} s_k\log {(p_k/s_k)} \right).
	\end{equation*}
\end{Theorem}
	Theorems \ref{th:upper_bound}, \ref{th:lower_bound_subspace}, and \ref{th:lower_bound_recovery} together show the proposed STAT-SVD method achieves the optimal rate of convergence in the general class of sparse low-rank tensors when $\log {(p_k/s_k)}\asymp \log{p_k}$, for $k \in J_s$.\\

\section{Data-driven Hyperparameter Selection}\label{sec:parameter}

In practice, the proposed procedure requires the input of hyperparameters $\sigma^2$ and $(r_1,r_2,\ldots,r_d)$. Since $\Y = \X + \Z$ and a significant portion of entries of $\X$ are zeros, the data-driven median estimator \citep{yang2016rate} can be used to estimate $\sigma$: $\hat{\sigma} = {\rm Median}(|\Y|)/0.6744$. Here $0.6744$ is the 75\% quantile of standard normal distribution. We have the following theoretical guarantee for $\hat{\sigma}$.
\begin{Proposition}[Concentration Inequality of $\hat{\sigma}^2$]\label{prop:sigma_estimate}
	Let $\hat{\sigma} = {\rm Median}(|\Y|)/z_{0.75}$. If $s = o(\sqrt{p\log p})$, then there exists universal constants $C$ and $c$, such that
	\begin{equation}\label{ineq:hat_sigma}
		\mathbb{P}\left(\frac{|\hat{\sigma} - \sigma|}{\sigma} \leq c \sqrt{\frac{\log p}{p}} \right) \geq 1 - C/p.
	\end{equation}

\end{Proposition}
Now we consider the data-driven selections of $(r_1,\ldots,r_d)$. Recall that we directly threshold on each mode and obtain $\tilde{\Y}$ in the initialization step of Algorithm 1. We propose to select $(r_1,\ldots, r_d)$ based on the singular values of $\tilde{\Y}$,
\begin{equation}\label{eq-rank-estimate}
	\hat r_k = \max \left\{r:\sigma_r(\mathcal{M}_k(\tilde{\Y}_k)) \geq \hat \sigma \delta_{|\hat I^{(0)}_k|,|\hat I^{(0)}_{-k}|}\right\}.
\end{equation}
Here, $\delta_{ij} = 1.02\left(\sqrt{i} + \sqrt{j} + \sqrt{2i\log \frac{ep_k}{i} + 2j\log \frac{ep_{-k}}{j} + 4\log p_k}\right)$ and $\hat{\sigma}$ is the estimated standard deviation. Under regularity conditions, one can show that $(\hat r_1, \ldots,\hat r_d)$ match the true ranks with high probability.
\begin{Proposition}\label{prop:rank_estimate}
Under the conditions of Theorem \ref{th:upper_bound} and Proposition \ref{prop:sigma_estimate}, we have $\hat r_k = r_k$ for each $k=1,\ldots,d$ with probability at least $1 - O((p_1+\cdots+p_d)/p^{1-\delta})$ for any $\delta>0$. 
\end{Proposition}

If neither $\sigma$ nor $\{r_k\}_{k=1}^d$ are known, we can first estimate $\sigma$ by MAD estimator $\hat{\sigma}$, then estimate $\{r_k\}_{k=1}^d$ by \eqref{eq-rank-estimate}, and feed all the estimations to Algorithm \ref{al:procedure}. We have the following theoretical guarantee for this fully data-driven method.
\begin{Theorem}\label{th:emperical-upper-bound}
Suppose we set $\sigma = \hat\sigma$ and $\br_k = \hat \br_k$ in Algorithm 1. If $s = o(\sqrt{p\log p})$, the conclusion of Theorem \ref{th:upper_bound} holds with probability at least $1-O\left((p_1+\cdots+p_d)/p^{1-\delta}\right)$ for any $\delta > 0$.
\end{Theorem}

\begin{Remark}\rm 
Similarly as some previous works on distribution-based methods for principal component number selection \citep{choi2017selecting}, the performance of the proposed $\hat\sigma$ and $\hat{r}_k$ may rely on specific Gaussian noise assumptions. In scenarios with general noise, some more empirical schemes can be applied. For example, to estimate $\sigma$, one can trim a portion of entries from $\Y$ with largest absolute values, and evaluate the \emph{trimmed variance} $\hat{\sigma}^2$ as the sample variance for remaining entries  \citep{serfling1984generalized}; for $r_k$, we can apply the \emph{cumulative percentage of total variation} criterion \citep[Chapter 6.1.1]{jolliffe2002principal} -- a commonly used in principal component analysis literature:
\begin{equation*}
\hat{r}_{k} = \argmin\left\{r: \frac{\sum_{i=1}^{r}\sigma_{i}^{2}(\mathcal{M}_{k}(\Y))}{\sum_{i=1}^{p_k} \sigma_{i}^{2}(\mathcal{M}_{k}(\Y))}>\rho\right\}.
\end{equation*}
Here, $\rho \in (0,1)$ is an empirical thresholding level. 
\end{Remark}

%% file: stat-svd4-simu.tex
\section{Numerical Analysis}\label{sec:simu}

\subsection{Simulation Studies}

We evaluate the numerical performance of the proposed STAT-SVD method by simulation studies on various synthetic datasets. In each setting, we first generate an $r_1$-by-$r_2$-by-$r_3$ tensor $\bar{\S}$ with i.i.d. standard normal entries, then rescale $\bar{\S}$ as $\S = \bar{\S} \cdot \lambda/\min_{1\leq k \leq 3}{\sigma_{r_k}(\mathcal{M}_k(\bar{\S}))}$
to ensure that $\sigma_{r_k}(\mathcal{M}_k(\S))\geq \lambda$. For $k = 1,2,3$, we generate singular subspaces $\bar{U}_k$ and indices subsets $\Omega_t$ with cardinality $s_k$ uniformly at random from $\mathbb{O}_{s_k, r_k}$ and $\{1,\ldots, p\}$, respectively. The combination of $\bar{U}_k$ and $\Omega_k$
$$(U_k)_{[i, :]} =
\left\{\begin{array}{ll}
(\bar{U}_k)_{[j, :]}, & i \in \Omega_k, \text{ and $i$ is the $j$-th element of $\Omega_k$};\\
0, & i\notin \Omega_k.\\
\end{array}\right.$$
yields a uniformly random sparse singular subspace in $\mathbb{O}_{p_k, r_k}(s_k)$.
Let $\X = \S\times_1 U_1\times_2 U_2\times_3 U_3$, $\Z \overset{iid}{\sim}N(0, \sigma^2)$, and $\Y = \X + \Z$ be the underlying parameter, noise, and observation tensors. To examine the performance of STAT-SVD, we apply the proposed Algorithm \ref{al:procedure} along with four baseline methods, HOSVD, HOOI, sparse HOSVD (S-HOSVD), and sparse HOOI (S-HOOI), on the same synthetic data for comparison. Here, HOSVD and HOOI are classical methods \citep{de2000multilinear,de2000best} that have been widely used in literature. S-HOSVD and S-HOOI are the sparse modifications of HOSVD and HOOI -- intuitively, S-HOSVD and S-HOOI are performed by replacing all regular matrix SVD steps in HOSVD and HOOI by sparse matrix SVD \citep{yang2014sparse}. The detailed implementation of S-HOSVD and S-HOOI are summarized in Section \ref{baseline-algorithm} of the supplementary materials. The experiments are repeated for 100 times in each setting.

\begin{figure} \centering
	\subfigure{
		\includegraphics[width =0.48\linewidth,height=2.0in]{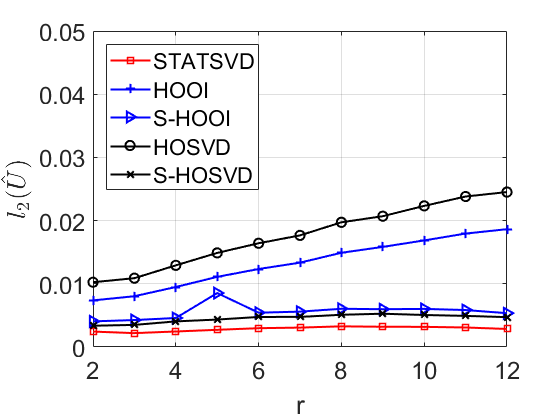}
	}
	\subfigure{
		\includegraphics[width =0.48\linewidth,height=2.0in]{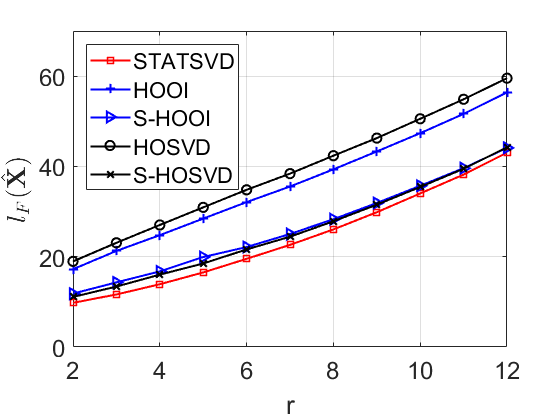}
	}
	\caption{Estimation error of $U_k$ (left panel) and $\X$ (right panel) for different methods. Here, $p_{1}=p_{2}=p_{3}=50$, $s_{1}=s_{2}=s_{3}=15$, $\sigma=1$, $\lambda=70$, and $r_{1}=r_{2}=r_{3}=r$ vary.}
	\label{fig:rank}
	\subfigure{
		\includegraphics[width =0.48\linewidth,height=2.0in]{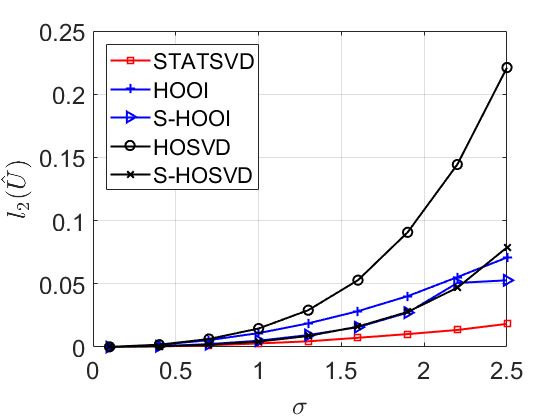}
	}
	\subfigure{
		\includegraphics[width =0.48\linewidth,height=2.0in]{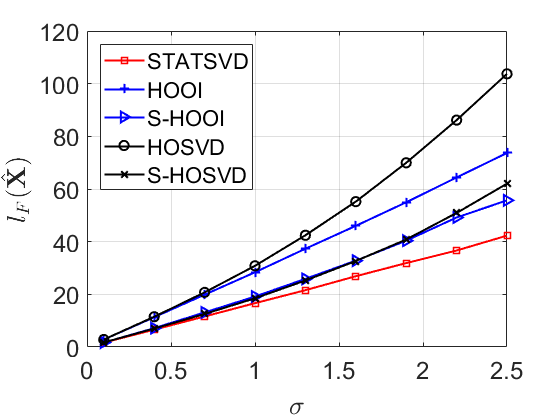}
	}
	\caption{Estimation error of $U_k$ (left panel) and $\X$ (right panel) for different methods. Here, $p_{1}=p_{2}=p_{3}=50$, $s_{1}=s_{2}=s_{3}=15$, $r_1=r_2=r_3=5$, $\lambda=70$, $\sigma$ varies.}
	\label{fig:noise}
\end{figure}

In the first simulation study, we compare the estimation errors of STAT-SVD and baseline methods (HOOI, S-HOOI, HOSVD, S-HOSVD) in average Frobenius norm,
$$l_2(\hat{U}) = \frac{1}{3}\sum_{k=1}^3 \left\|\sin\Theta(\hat{U}_k, U_k)\right\|_F,\quad l(\hat{\X}) = \|\hat{\X} - \X\|_F.$$
For hyperparameters, we use the median estimator $\hat{\sigma}$ in STAT-SVD and all baseline methods; since $\sin\Theta$ distance is one of the most important error quantification in tensor SVD analysis and one requires the correct $r_k$ to evaluate the $\sin \Theta$ distance between singular subspaces, we use the true ranks $r_k$ for all implementations. Fix $p_{1}=p_{2}=p_{3}=50$, $s_{1}=s_{2}=s_{3}=15$, $\lambda=70$, $r=r_1=r_2=r_3$, we specifically consider two scenarios: (1) $\sigma=1$, $r$ varies from 2 to 12; (2) $r=5$, $\sigma$ ranges from 0.1 to 2.5. 
Although the signal-to-noise ratio (SNR, $\lambda/\sigma$) here seems large, $\lambda$ represents the singular value of the each matricization that measures the signal strength of the \emph{whole tensor}; while $\sigma$ is the standard deviation of each $Z_{ijk}$ that quantifies the noise level of each \emph{single entry}. By random matrix theory \citep{vershynin2010introduction}, the singular values of $\mathcal{M}_k(\Z)$ are around $\sigma\sqrt{p_{-k}} \approx 5$ to $125$, which is comparable to the signal $\mathcal{M}_k(\X)$. As one can see from the numerical results in Figures \ref{fig:rank} and \ref{fig:noise}, although all methods yield smaller estimation error with smaller noise level, STAT-SVD significantly outperforms all other schemes in estimations of both subspaces and original tensors. 

We also consider the accuracy of hyperparameters' estimation. Under the previous simulation setting, we examine the performance of MAD estimator $\hat{\sigma}$ and rank estimator $\hat{r}_k$ in Section \ref{sec:parameter}. The results are provided in Table \ref{tab:rank}. One can see that the proposed $\hat{\sigma}$ and $\hat{r}_k$ provide reasonable estimations. The rank estimation is particularly accurate when noise level is moderate.
\begin{table}
	\centering
	\begin{tabular}{lllllllll}
		\hline
		$\sigma$        & 0.1 & 0.3  & 0.5 & 0.7  & 0.9 & 1.1 & 1.3 & 1.5		\\ \hline
$|\hat \sigma - \sigma|$ & 0.003 & 0.008 & 0.013 & 0.018  & 0.021 & 0.025 & 0.026 & 0.028
\\ \hline
$|\hat r - r|$ & 0 & 0 & 0 & 0 & 0 & 0 & 0.54 & 0.98 \\\hline\hline 
$r_k$        & 3 & 4  & 5 & 6  & 7 & 8  & 9 & 10 \\\hline
$|\hat \sigma - \sigma|$ & 0.019 & 0.021 & 0.022 & 0.023  & 0.024 & 0.025 & 0.025 & 0.026 \\\hline 
$|\hat r - r|$ & 0 & 0 & 0 & 0 & 0 & 0 & 0 & 0 \\\hline 
	\end{tabular}
	\caption{Rank and noise level estimation accuracy. Here, $p_{1}=p_{2}=p_{3}=50$, $s_{1}=s_{2}=s_{3}=15$,  $\lambda=70$. The first three rows explore the setting where $r_1=r_2=r_3=5$ with $\sigma$ varies; the last three rows consider the case with fixed $\sigma=1$ and varying rank. The errors are averaged over 100 repetitions.}
	\label{tab:rank}
\end{table}

\begin{figure}
	\subfigure{
		\includegraphics[width =0.48\linewidth,height=2.0in]{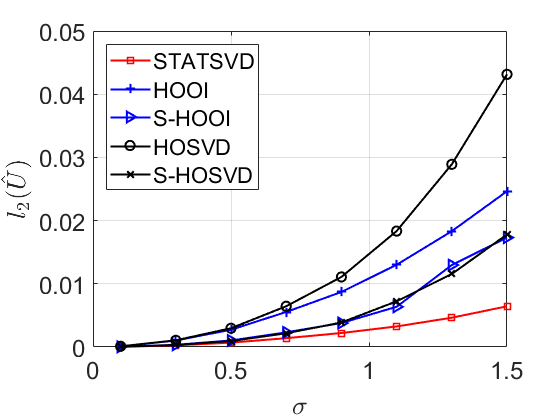}
	}
	\subfigure{
		\includegraphics[width =0.48\linewidth,height=2.0in]{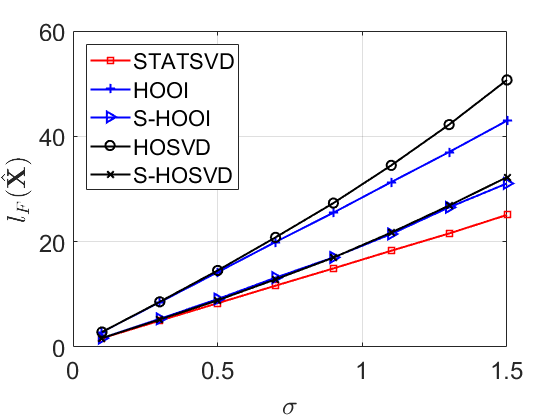}
	}
	\caption{Estimation error of $U_k$ (left panel) and $\X$ (right panel) with uniform distributed noise.}
	\label{fig:noise_unif}
\end{figure}

In previous sections, the presentation and analysis were mostly focused on Gaussian noise case. Next, we consider the setting that the noises are uniformly distributed on $[-\sigma\sqrt{3}, \sigma\sqrt{3}]$. Let $p_{1}=p_{2}=p_{3}=50$, $s_{1}=s_{2}=s_{3}=15$, $r_1=r_2=r_3=5$, $\lambda=70$, $\sigma$ varies from 0.1 to 1.5. As we can see from the estimation error results in Figure \ref{fig:noise_unif}, STAT-SVD still achieves significantly better performance than other methods.

\begin{figure} \centering
	\subfigure {
		\includegraphics[width =0.48\linewidth,height=2.0in]{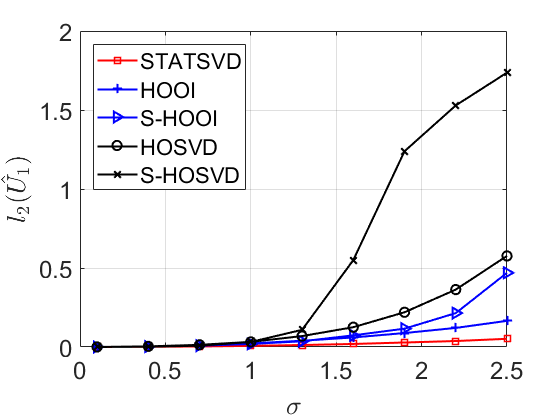}
	}
	\subfigure {
		\includegraphics[width =0.48\linewidth,height=2.0in]{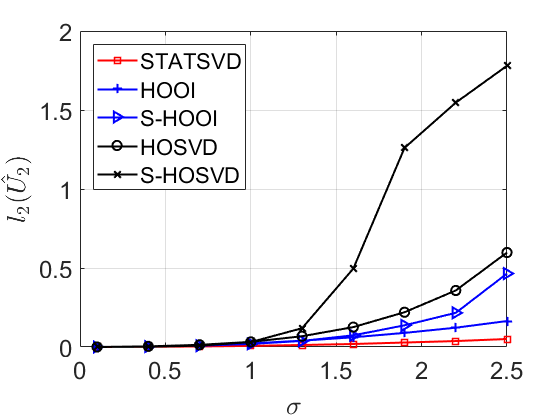}
	}\\
	\subfigure {
		\includegraphics[width =0.48\linewidth,height=2.0in]{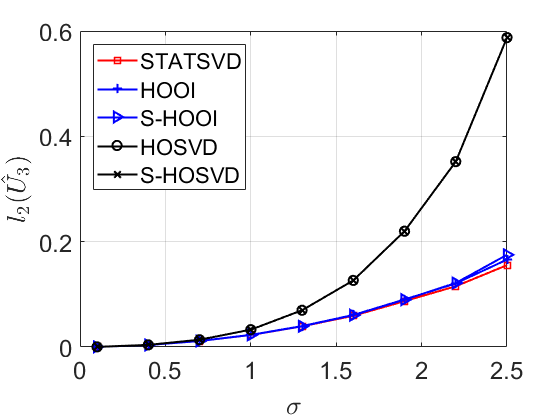}
	}
	\subfigure {
		\includegraphics[width =0.48\linewidth,height=2.0in]{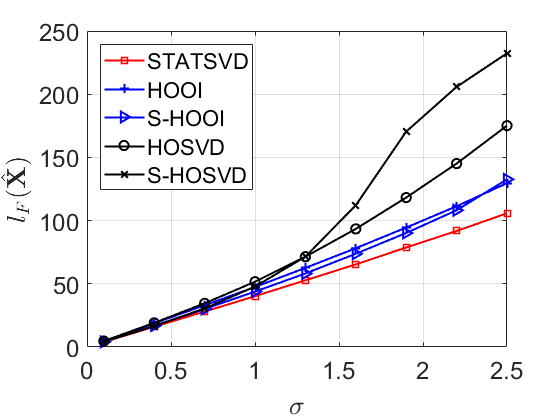}
	}
	\caption{Estimation error of $U_1$ (left top), $U_2$ (right top), $U_3$ (left bottom) and $\X$ (right bottom) in partial sparse setting. Here, $\lambda=60$,  $p_{1}=p_{2}=p_{3}=50$, $s_{3}=p_{3}$, $s_{1}=s_{2}=20$, $r_{1}=r_{2}=r_{3}=10$.
	}
	\label{fig:partial}
\end{figure}

As we have discussed before, in many applications, the tensor dataset possesses sparsity structure in only part of directions. Thus, we turn to a setting that $\X$ contains both sparse and dense modes. Specifically, we set $p_{1}=p_{2}=p_{3}=50$, $r_{1}=r_{2}=r_{3}=10$, $s_{1}=s_{2}=20$, $\lambda=60$, and $s_3 = p_3$, so that $\X$ is sparse along mode-1, -2 and dense along mode-3. The singular subspace estimation error with varying noise level is evaluated and shown in Figure \ref{fig:partial}. We can see STAT-SVD significantly outperforms the other methods on the estimation of sparse modes ($\hat{U}_{1}$ and $\hat{U}_{2}$). More interestingly, STAT-SVD also estimates the non-sparse singular subspace $U_{3}$ slightly more accurately, especially when $\sigma^2$ is large. In fact, different modes and singular subspaces of any specific tensor dataset are a unity rather than separate objects, so more accurate estimation of dense mode $U_3$ is possible when one can fully utilize the sparsity in $U_1$ and $U_2$. Especially for STAT-SVD, with the proposed double projection \& thresholding scheme (Algorithm \ref{al:refinement-sparse}), one gets significant better estimations on $U_{1}$ and $U_{2}$ than the baseline methods in each iteration so more precise projection \eqref{eq:update-non-sparse-projection} can be achieved when updating dense mode singular subspace $\hat{U}_3^{(t)}$, and a slightly better final estimation of $U_3$ can be achieved.

As the time cost is another critical issue in high-dimensional data analysis, we summarize the computational complexity for both initialization and each iteration of STAT-SVD and baseline methods into Table \ref{tab:complexity}. We also compare the running time of STAT-SVD and other algorithms by simulations. As one can see from Tables \ref{tab:complexity} and \ref{tab:speed}, HOOI, HOSVD, S-HOOI, and S-HOSVD are all slower than STAT-SVD -- this is because the embedded sparse matrix SVD or regular matrix SVD in baselines are computationally expensive, especially for the high-dimensional cases. In contrast, STAT-SVD is much faster, as the efficient truncation makes sure that only the significant parts of the data are extensively applied in computation, i.e., we only need to perform SVD on the $s_k$-by-$s_{-k}$ submatrix instead of the original $p_k$-by-$p_{-k}$ one.

\begin{table}
\centering
\begin{tabular}{lllllll}
\hline
       & STAT-SVD  & HOSVD  & HOOI  & S-HOSVD  & S-HOOI   \\ \hline
initialization & $O((p_1^d+s_1^{d+1}))$ & $O(p_1^{d+1})$ & $O(p_1^{d+1})$  & $O(Tp_1^{d+1})$  & $O(Tp_1^{d+1})$     \\ \hline
per-iteration   & $O(p_1^{d-1}r_1 + s_1^{d+1})$ & 0 & $O(p_1^{d+1})$ & 0 &  $O(Tp_1^{d+1})$  \\ \hline
\end{tabular}
\caption{Computational complexity of STAT-SVD and baseline methods for order-$d$ sparse tensor SVD. Here, each mode share the same $p_i$, $s_i$, and $r_i$. $T$ denotes the number of iteration times. }
\label{tab:complexity}
\end{table}

\begin{table}
\centering
\begin{tabular}{lllllllllll}
\hline
$p_i$       & 50  & 80  & 110  & 140  & 170  & 200  & 230  & 260   & 290   & 320   \\ \hline
STAT-SVD & 0.1 & 0.3 & 0.7  & 1.3  & 2.3  & 3.5  & 5.9  & 8.3   & 12.8  & 18.9  \\ \hline
HOSVD    & 0.1 & 0.2 & 0.9  & 2.0  & 4.7  & 8.4 & 13.0 & 21.2  & 33.2  & 64.5  \\ \hline
HOOI     & 0.1 & 0.3 & 1.0  & 2.4  & 5.2  & 9.4 & 14.6 & 22.9 & 36.3  & 66.7  \\ \hline
S-HOSVD  & 1.7 & 3.2 & 6.8 & 11.3 & 17.2 & 54.5 & 83.0 & 590.3 & 1340.7 & 1918.2 \\ \hline
S-HOOI   & 1.9 & 3.8 & 7.0 & 11.9 & 16.8 & 56.9 & 91.1 & 482.7 & 1255.5 & 1916.7 \\ \hline
\end{tabular}
\caption{Running time (unit: seconds) of STAT-SVD and baseline methods. Here, $\lambda = 70$, $\sigma = 1$, $s_1=s_2=s_3=10$, $r_1=r_2=r_3=5$, $p_1=p_2=p_3$ varies.} 
\label{tab:speed}
\end{table}

%% file: stat-svd5-real.tex
\subsection{Mortality Rate Data Analysis}
We illustrate the power of STAT-SVD through a demographic example. The mortality rate, i.e. the number of deaths divided by the total number of population, provides interesting insights to demographic information of the certain area, period, and age span. The Berkeley Human Mortality Database \citep{wilmoth2016human} contains a good source of morality rate data aligned by countries, ages, and years. 
We aim to analyze the mortality rate among 26 European countries for ages 0 to 95 from 1959 to 2010. The data tensor $\Y$ is of dimension $96 \times 52 \times 26$ with three modes representing age, year, and country, respectively. Since the mortality rate is relatively steady from teenagers to adults (see, e.g., \cite{minino2013death}), it is reasonable to assume that the underlying loadings of Mode age are sparse after taking differential transformation. Specifically, let $D \in \mathbb{R}^{96\times96}$ be a secondary difference matrix: $D_{ij} = -2$ if $i=j$; $D_{ij} = 1$ if $|i-j|=1$; $D_{ij}=0$ if $|i-j|\geq 2$. We pre-process the data by multiplying the secondary difference matrix along Modes age: $\tilde{\Y} = \Y \times_1 D$. 

We aim to apply sparse tensor SVD on $\tilde{\Y}$ with $J_s = \{1\}$. To achieve more robust performance, the hyperparameters are selected via empirical methods instead of the Gaussian-noise-based procedure discussed in Section \ref{sec:parameter}. $r_k$ is selected via cumulative percentage of total variation criterion \citep[Chapter 6.1.1]{jolliffe2002principal},
\begin{equation*}
\hat{r}_{k} = \argmin\left\{r: \frac{\sum_{i=1}^{r}\sigma_{i}^{2}(\mathcal{M}_{k}(\X))}{\sum_{i=1}^{p_k} \sigma_{i}^{2}(\mathcal{M}_{k}(\X))}>0.5\right\}.
\end{equation*}
For the noise level, we trim $15\%$ largest entries of $\Y$ in absolute value, then set $\hat{\sigma}$ as the sample standard deviation of the remaining entries of $\Y$. The estimated rank and noise level of $\tilde{\Y}$ are $(\hat r_1 ,\hat{r}_2, \hat{r}_3) = (2,5,4)$ and $\hat{\sigma} = 0.0014$. In addition to STAT-SVD, we also apply S-HOOI and S-HOSVD for comparison. Suppose $\tilde{U}, \hat{V}, \hat{W}$ are the resulting singular subspaces of $\tilde{\Y}$. Then we transform $\tilde{U}$ (singular subspace of Mode age) back via $\hat{U} = D^{-1}\tilde{U}$. 

\begin{figure} \centering
\includegraphics[width=.85\columnwidth,height=6.8cm]{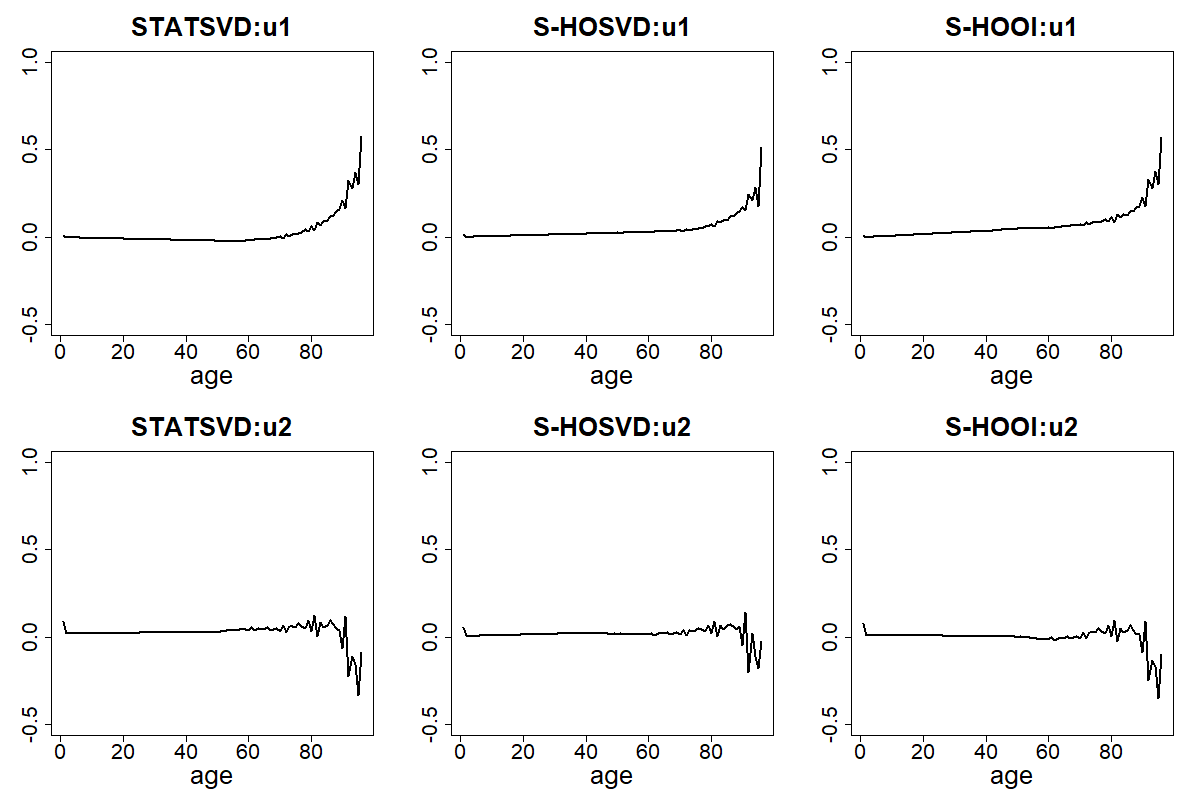}
\caption{Mortality rate data: $\hat{u}_{1}$, $\hat{u}_2$}
\label{fig:uplot}
\includegraphics[width=.85\columnwidth,height=6.8cm]{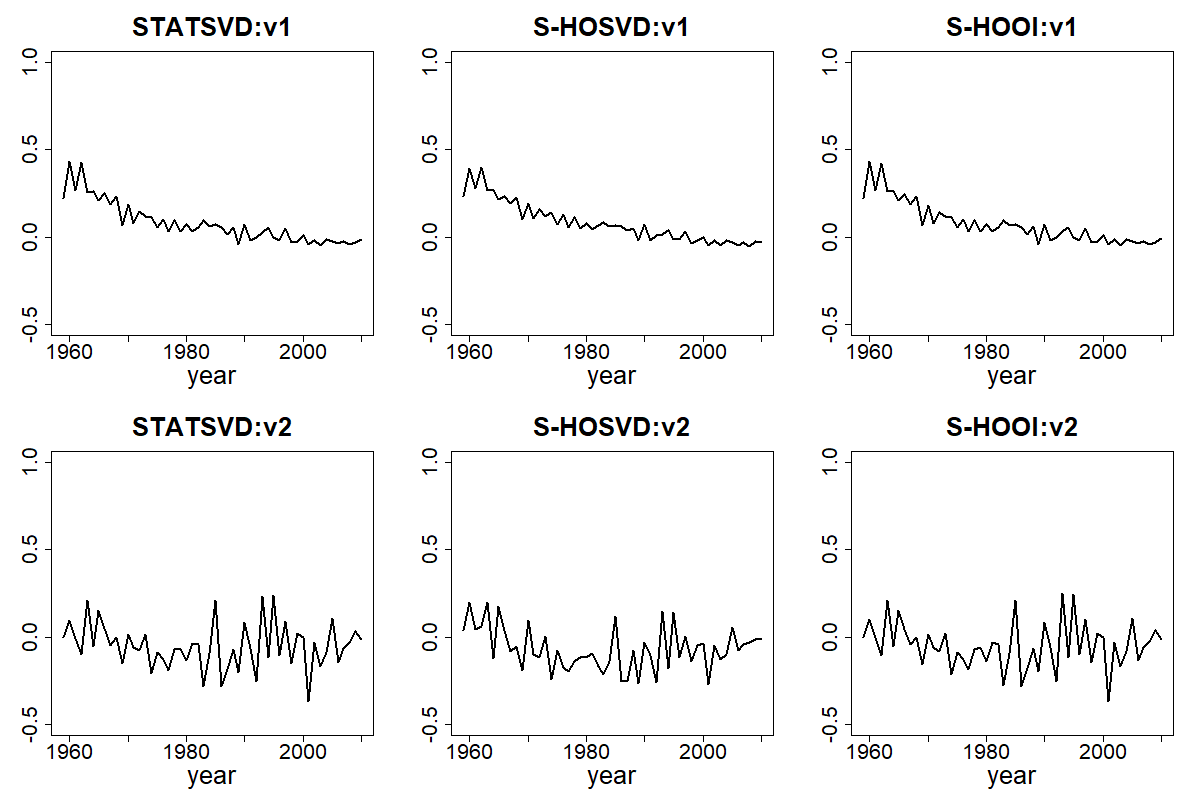}
\caption{Mortality rate data: $\hat{v}_1$, $\hat{v}_2$}
\label{fig:vplot}
\vskip.5cm
	\includegraphics[width=.7\columnwidth,height=4.8cm]{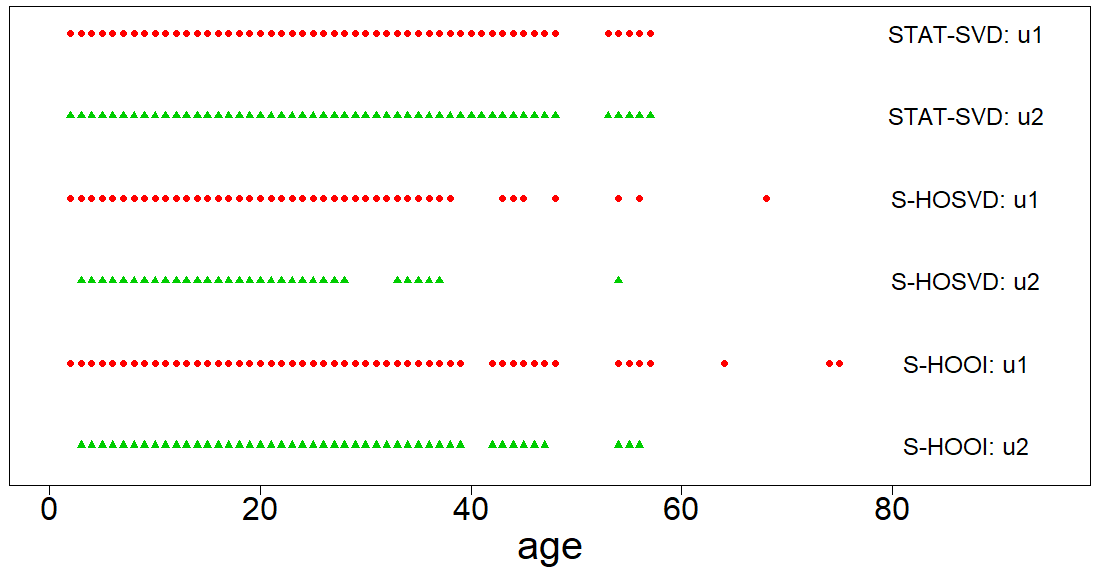}
	\caption{Sparse indexes of Mode age after secondary differential transformation}
	\label{fig:support}
\end{figure}

We first compare the first two singular vectors of Modes age and year, $\hat{u}_1, \hat{u}_2, \hat{v}_1,$ and $\hat{v}_2$, in Figures \ref{fig:uplot} and \ref{fig:vplot}. We can see from $\hat{u}_{1}$ that infants aged at 0-2 and old people aged over 55 have a higher risk of dying, and people aged from 2 to 50 have a relatively steady death rate. In addition, the shape of $u_{2}$ further suggests additional factors of mortality rate in certain periods. For example, there may be more complicated patterns in mortality rate for the infants and elderly aged over 80 due to the high risk of death in these two periods. From $\hat{v}_1$, we can see the mortality rate in these European countries declines significantly from the year 1955 to 2010, while $\hat{v}_2$ does not give any significant pattern. 
Since the three methods produce similar plots on values of singular vectors, we further compare the estimation support indexes in Figure \ref{fig:support}. Since the mortality rate is steady from childhood to adults, one expects that the zero indexes would cover such an age span for $\hat{u}_1$ and $\hat{u}_2$. The outcome of STAT-SVD matches this phenomenon. In comparison, the sparsity patterns given by S-HOSVD and S-HOOI are less interpretable.

\begin{figure} \centering
	\includegraphics[height=2.0in]{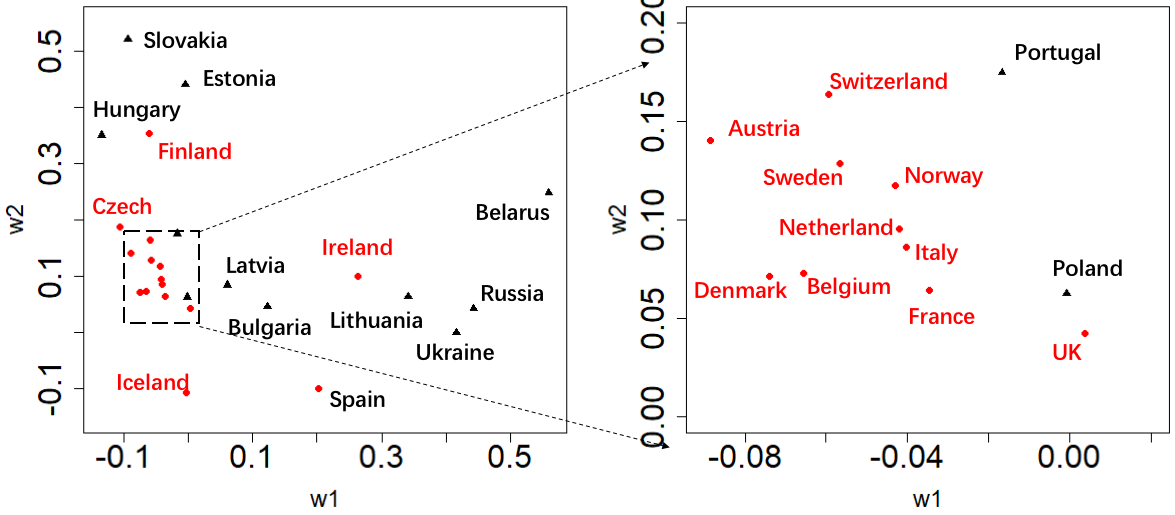}
	\caption{Mortality rate data: $\hat{w}_{1}$ and $\hat{w}_2$ of different countries. Red circles and black triangles represent the countries with top 50\% GDP per capita in Europe and the other countries, respectively.}
	\label{fig:ncountry}
	\vskip.5cm
	\includegraphics[height=2.0in]{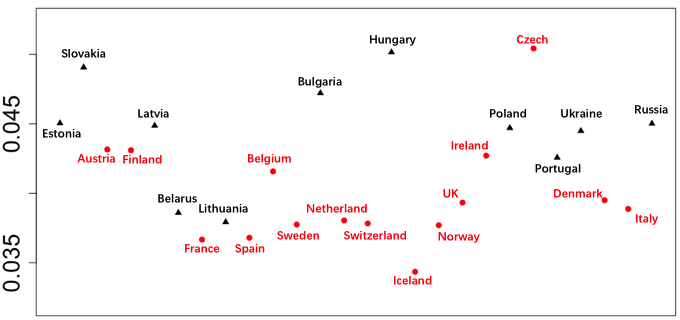}
	\caption{Mortality rate means for countries in Europe. Red circles and black triangles represent the countries with top 50\% GDP per capita in Europe and the other countries, respectively.}
	\label{fig:ncountry-naive}
\end{figure}

Next, we consider the mortality rate pattern for countries. In particular, we plot $\hat{w}_{1}$ against $\hat{w}_{2}$ for each country in Figure \ref{fig:ncountry}, refer to the 2014 European GDP per capita ranking data\footnote{Link: \url{http://statisticstimes.com/economy/countries-by-gdp-capita.php}}, search for the countries whose GDP per capita is ranked as top $50\%$ among all European counties, then highlight these countries with red circle markers and the other with black triangles. We can clearly see two clusters in Figure \ref{fig:ncountry} that countries with more GDP per capita have smaller values of $\hat{w}_{1}$ and $\hat{w}_2$, which implies a lower mortality rate. Also, wealthier countries are highly clustered in the graph, which indicates the common pattern of the death rate for these countries. In contrast, we also calculate the mean mortality rates of these countries and mark them with different colors in Figure \ref{fig:ncountry-naive}. By comparing Figures \ref{fig:ncountry} and \ref{fig:ncountry-naive}, the mortality data clustering performance of STAT-SVD is significantly better than the one by mean estimations.

%% file: stat-svd6-proof.tex
\section{Proofs}\label{sec:proofs}

\subsection{Proof of Theorem \ref{th:upper_bound}} Since the proof for this upper bound result is fairly complicated, we divide the proof into steps. Note that although we do not have refinement for non-sparse modes, to make the proof consistent, we first extend the refinement step to non-sparse algorithm. Specifically, we also apply algorithm 2 on $\hat{U}_{k}^{(t)}$ for $k \not \in J_s$, with $\eta_k=\bar{\eta}_k=0$ i.e. no truncation. This modification does not change the algorithm at all but makes the following analysis consistent for both sparse and non-sparse modes.\\
Without loss of generality, we assume $\sigma = 1$ throughout this proof. The idea of proving this theorem is that, we first impose a series of conditions, then prove the statement under these conditions, and finally prove these conditions hold with high probability.
\begin{enumerate}[leftmargin=*]
	\item[Step 1] (Introduction of Notations and Conditions) We introduce or rephrase the following list of notations.
	\begin{enumerate}[label=(\subscript{N}{\arabic*})]
		\item Matricizations, for $k=1,\ldots, d$,
		\begin{equation*}
			X_k = \mathcal{M}_k\left(\X\right),\quad Y_k = \mathcal{M}_k\left(\Y\right),\quad Z_k = \mathcal{M}_k\left(\Z\right).
		\end{equation*}
		\item Index sets for sparse mode $k\in J_{s}$: define
		\begin{equation}\label{eq:step-1-index}
			\begin{split}
				\text{(True support)}\quad & I_k = \left\{i: (X_k)_{[i, :]} \neq 0\right\}, \\
				\text{(Significant support)}\quad & I_k^{(0)} = \left\{i: \left\|(X_k)_{[i, :]}\right\|_2^2 \geq
				16\sigma s_{-k}\log p\right\},\\
				\text{(Selected support)} \quad & \hat{I}_k^{(0)} = \Bigg\{i: \left\|(Y_k)_{[i, :]}\right\|_2^2 \geq \sigma^2 \left(p_{-k} + 2\sqrt{p_{-k}\log p} + 2 \log p\right)\\
				& \quad \text{or } \max_{j} |(Y_k)_{[i,j]}| \geq 2\sigma\sqrt{\log p}\Bigg\},\quad k=1,\ldots, d.
			\end{split}
		\end{equation}
		For $t\geq 1$, we also define $\hat{I}_k^{(t)}$ and $\hat{J}_k^{(t)}$ with $k \in J_{s}$:
		\begin{equation}\label{eq:support-I-J}
			\begin{split}
				(\text{selected support for $A_k$ in $t$-th step})\quad & \hat{I}_k^{(t)} = \left\{i: \left\|(A_k^{(t)})_{[i, :]}\right\|_2^2 \geq \eta_k\right\},\\	
				(\text{selected support for $\bar{A}_k$ in $t$-th step})\quad & \hat{J}_k^{(t)} = \left\{i: \left\|(\bar{A}_k^{(t)})_{[i, :]}\right\|_2^2 \geq \bar{\eta}_k\right\},\\	
			\end{split}
		\end{equation}
		where $\eta_k = \sigma^2 \left(r_{-k} + 2\sqrt{r_{-k}\log p} + 2\log p \right)$, $\bar{\eta}_k = \sigma^2 \left(r_k + 2\sqrt{r_k\log p} + 2\log p\right)$.\\
		For non-sparse mode $k \not \in J_{s}$, define $I_{k}=I_{k}^{(t)}=\hat{I}_k^{(t)}=\hat{J}_k^{(t)}=[1:p_k]$. It is easy to see that for any $k$, $I_k, I_k^{(0)}\hat{I}_k^{(t)}, \hat{J}_k^{(t)}$ are all subsets of $\{1,\ldots, p_k\}$. We also define
		\begin{equation*}
			\begin{split}
				& I_{-k} = I_{k+1}\otimes \cdots \otimes I_d\otimes I_1\otimes\cdots \otimes I_{k-1}.
			\end{split}
		\end{equation*}
		Moreover, $\hat{I}_{-k}^{(t)}$, $\hat{J}_{-k}^{(t)}$, $I_{-k}^{(t)}$, $J_{-k}^{(t)}$ are defined in the similar way.

		\item Index projections: for any $I_k \subseteq \{1,\ldots, p_k\}$, we define
		\begin{equation*}
			D_{I_k} = \diag\left(1_{I_k}\right) \in \mathbb{R}^{p_k\times p_k}, \quad \text{if $I_k\subseteq \{1,\ldots, p_k\}$ is any index subset};
		\end{equation*}
		$D_{I_k}$ can be interpreted as the projection matrix, which set all rows with index not in $I_k$ to zero.
		\item Loadings:
		\begin{equation*}
			\begin{split}
				& U_{-k} = U_{k+1}\otimes\cdots\otimes U_{d}\otimes U_1\otimes \cdots \otimes U_{k-1} \in \mathbb{O}_{p_{-k}\times r_{-k}}, \quad k=1,\ldots, d;\\
				& \hat{U}_{-k}^{(t)} = \hat{U}_{k+1}^{(t-1)}\otimes\cdots\otimes \hat{U}_{d}^{(t-1)}\otimes \hat{U}_1^{(t)}\otimes \cdots \otimes \hat{U}_{k-1}^{(t)} \in \mathbb{O}_{p_{-k}\times r_{-k}}, \quad k=1,\ldots, d.
			\end{split}
		\end{equation*}
		\begin{equation*}
			\begin{split}
				& V_k = \SVD_{r_k}\left((X_k U_{-k})^\top\right)\in \mathbb{O}_{r_{-k}, r_k},~~~~ \text{i.e., leading $r_k$ right singular vectors of $X_kU_{-k}$};\\
				& \hat{V}_k^{(t)} = \SVD_{r_k}\left((D_{\hat{I}_k^{(t)}}Y_k\hat{U}_{-k}^{(t)})^\top\right)\in \mathbb{O}_{r_{-k}, r_k}, ~~~~ \text{i.e., leading $r_k$ right singular vectors of $D_{\hat{I}_k^{(t)}}Y_k\hat{U}_{-k}^{(t)}$}.
			\end{split}
		\end{equation*}	
		Combining these definitions, we can rewrite $A_k^{(t)}$ and $\bar{A}_k^{(t)}$ as
		\begin{equation*}
			\begin{split}
				& A_k^{(t)} = Y_k\hat{U}_{-k}^{(t)}, \quad B_k^{(t)} = D_{\hat{I}_k^{(t)}} A_k^{(t)} = D_{\hat{I}_k^{(t)}}Y_k\hat{U}_{-k}^{(t)},\\
				& \bar{A}_k^{(t)} = Y_k\hat{U}_{-k}^{(t)}\hat{V}_k^{(t)}, \quad \bar{B}_k^{(t)}  = D_{\hat{J}_k^{(t)}}Y_k\hat{U}_{-k}^{(t)}\hat{V}_k^{(t)}.
			\end{split}
		\end{equation*}
		We can immediately see that $A_k^{(t)}, B_k^{(t)}, \bar{A}_k^{(t)}, \bar{B}_k^{(t)}$ are essentially projections of $Y_k$. Since $Y_k = X_k + Z_k$, $A_k^{(t)}, B_k^{(t)}, \bar{A}_k^{(t)}, \bar{B}_k^{(t)}$ can be decomposed accordingly as
		\begin{equation}\label{eq:A_k-B_k-bar_A_k-bar_B_k}
			\begin{split}
				A_k^{(t)} = A_k^{(X, t)} +  A_k^{(Z, t)} \quad  & A_k^{(X, t)} = X_k\hat{U}_{-k}^{(t)}, \quad A_k^{(Z, t)} = Z_k\hat{U}_{-k}^{(t)};\\
				B_k^{(t)} = B_k^{(X, t)} + B_k^{(Z, t)} \quad & B_k^{(X, t)} = D_{\hat{I}_k^{(t)}}X_k\hat{U}_{-k}^{(t)},\quad  B_k^{(Z, t)} = D_{\hat{I}_k^{(t)}}Z_k\hat{U}_{-k}^{(t)}\\
				\bar{A}_k^{(t)} = \bar{A}_k^{(X, t)} + \bar{A}_k^{(Z, t)} \quad & \bar{A}_k^{(X, t)} = X_k\hat{U}_{-k}^{(t)}\hat{V}_k^{(t)},\quad \bar{A}_k^{(Z, t)} = Z_k\hat{U}_{-k}^{(t)}\hat{V}_k^{(t)} \\
				\bar{B}_k^{(t)} = \bar{B}_k^{(X, t)} + \bar{B}_k^{(Z, t)}\quad & \bar{B}_k^{(X, t)} = D_{\hat{J}_k^{(t)}}X_k\hat{U}_{-k}^{(t)}\hat{V}_k^{(t)},\quad \bar{B}_k^{(Z, t)} = D_{\hat{J}_k^{(t)}}Z_k\hat{U}_{-k}^{(t)}\hat{V}_k^{(t)}.
			\end{split}
		\end{equation}
		
		\item Error Bounds: for $t = 0, 1, \ldots, k =1,\ldots, d$,
		\begin{equation*}
			\begin{split}
				& E_{k}^{(t)} = \left\|(\hat{U}_{k\perp}^{(t)})^\top X_k\right\|_F, \quad F_{k}^{(t)} = \left\|\sin\Theta\left(\hat{U}_k^{(t)}, U_k\right)\right\|_F.
			\end{split}
		\end{equation*}
		For $k = 1,\ldots, d, t = 1,\ldots$,  we further define
		\begin{equation*}
			K_k^{(t)} = \left\|\sin\Theta\left(\hat{U}_{-k}^{(t)}\hat{V}_k^{(t)}, U_{-k}V_k\right)\right\|.
		\end{equation*}
	\end{enumerate}
	
	We also introduce the following conditions for the proof of this theorem.
	\begin{enumerate}[label=(\subscript{A}{\arabic*})]
		\item
		\begin{equation}\label{ineq:A-1}
			\max_{(i_1,\ldots, i_d)\in I_1\times \cdots \times I_d} \left|Z_{i_1,\ldots, i_d}\right| \leq 2\sqrt{\log p}
		\end{equation}
		\item
		\begin{equation}\label{ineq:A-2}
			\begin{split}
				& \left\|\Z_{[I_1,\ldots, I_d]}\right\|_F^2 \leq s + 2\sqrt{s \log p} + 2 \log p\\
				& \left\|\Z \times_1 U_1^\top \times \cdots \times_d U_d^\top \right\|_F^2 \leq r + 2\sqrt{r\log p} + 2\log p.
			\end{split}
		\end{equation}
		\item
		\begin{equation}\label{ineq:A-3}
			\forall 1\leq k \leq d, \quad \left\|Z_{k, [I_k, I_{-k}]}\right\| \leq \sqrt{s_k} + \sqrt{s_{-k}} + 2\sqrt{\log p}.
		\end{equation}
		\item $\forall 1\leq k \leq d$,
		\begin{equation}\label{ineq:A-4}
			N_k := \left\|(Z_k)_{[I_k, :]}U_{-k}V_k\right\| \leq \sqrt{s_k}+\sqrt{r_k}+2\sqrt{\log p}.
		\end{equation}
		\item $\forall 1\leq k \leq d$,
		\begin{equation}\label{ineq:A-5}
			M_k := \sup_{W_l \in \mathbb{R}^{s_l\times r_l}, \|W_l\|\leq 1} \left\|Z_{k, [I_k, I_{-k}]} W_{-k}\right\| \leq 2\left(\sqrt{s_k} + \sqrt{r_{-k}} + \sum_{l\neq k} (\sqrt{s_lr_l} + \sqrt{\log p})\right),
		\end{equation}
		where $W_{-k} = W_{k+1}\otimes \cdots \otimes W_d \otimes W_1 \otimes \cdots \otimes W_{k-1} \in \mathbb{O}_{s_{-k}, r_{-k}}$.
		\item Support consistency condition:
		\begin{equation}\label{ineq:A-6}
			\begin{split}
				& I_k^{(0)} \subseteq \hat{I}_k^{(0)},\\
				& \hat{I}_k^{(t)} \subseteq I_k,\quad t=0,\ldots, t_{\max},\\
				&  \hat{J}_k^{(t)} \subseteq I_k,\quad t=1,\ldots, t_{\max}.
			\end{split}
		\end{equation}
	\end{enumerate}

	\item[Step 2] (Theoretical guarantees for initialization: $\hat{U}_k^{(0)}$) In this second step, we show that under ($A_1$) -- ($A_7$), the initialization estimator $\hat{U}_k^{(0)}$ satisfies the following two inequalities for any $k=1,2,\cdots,d$:
	\begin{equation}\label{ineq:step2-result-1}
		\left\|\sin\Theta\left(\hat{U}_k^{(0)}, U_k\right)\right\|_F \leq \frac{12\sqrt{s\log p}}{\lambda_k}.
	\end{equation}
	\begin{equation}\label{ineq:step2-result-2}
		\left\|(\hat{U}_{k\perp}^{(0)})^\top X_k \right\|_F \leq 12\sqrt{ds\log p}.
	\end{equation}
	Recall that $\hat{U}_k^{(0)} = \SVD_{r_k}\left(\mathcal{M}_k\left(\tilde{\Y}^{(0)}\right)\right) = \SVD_{r_k}\left(D_{\hat{I}_{k}^{(0)}}Y_{k} D_{\hat{I}_{-k}^{(0)}}\right)$. Since $D_{\hat{I}_{k}^{(0)}}Y_{k} D_{\hat{I}_{-k}^{(0)}} = D_{\hat{I}_{k}^{(0)}}X_{k} D_{\hat{I}_{-k}^{(0)}} + D_{\hat{I}_{k}^{(0)}}Z_{k} D_{\hat{I}_{-k}^{(0)}}$, by the unilateral perturbation bound result (Proposition 1 in \cite{cai2016rate}),
	\begin{equation}\label{ineq:sine-theta-U^0}
		\begin{split}
			\left\|\sin\Theta\left(\hat{U}_k^{(0)}, U_k\right)\right\|_F \leq \frac{\sigma_{r_k}(U_k^\top D_{\hat{I}_{k}^{(0)}}Y_{k} D_{\hat{I}_{-k}^{(0)}}) \|U_{k\perp}^\top D_{\hat{I}_{k}^{(0)}}Y_{k} D_{\hat{I}_{-k}^{(0)}}\|_F}{\sigma_{r_k}^2(U_k^\top D_{\hat{I}_{k}^{(0)}}Y_{k} D_{\hat{I}_{-k}^{(0)}}) - \sigma_{r_k+1}^2(D_{\hat{I}_{k}^{(0)}}Y_{k} D_{\hat{I}_{-k}^{(0)}})}.
		\end{split}
	\end{equation}
	We set $C_0 \geq 80\sqrt{d}$, recall that
	\begin{equation}\label{ineq:lambda_k_min}
		\lambda_k \geq 80\sqrt{ds\log p}, \quad k=1,\ldots, d,
	\end{equation}
	we analyze each of the three key terms in the right hand side of \eqref{ineq:sine-theta-U^0} as follows,
	\begin{enumerate}[leftmargin=*]
		\item Since \eqref{ineq:A-6} holds, i.e. $I_k^{(0)}\subseteq \hat{I}_k^{(0)}$, we know the non-zero part of $D_{I_{k}^{(0)}}Y_{k} D_{I_{-k}^{(0)}}$ is a submatrix of $D_{\hat{I}_{k}^{(0)}}Y_{k} D_{\hat{I}_{-k}^{(0)}}$, then
		\begin{equation*}
			\begin{split}
				& \sigma_{r_k}(U_k^\top D_{\hat{I}_{k}^{(0)}}Y_{k} D_{\hat{I}_{-k}^{(0)}}) \geq \sigma_{r_k}\left(U_k^\top  D_{\hat{I}_{k}^{(0)}}X_{k} D_{\hat{I}_{-k}^{(0)}}\right) - \left\|U_k^\top D_{\hat{I}_{k}^{(0)}} Z_{k} D_{\hat{I}_{-k}^{(0)}}\right\|\\
				\geq & \sigma_{r_k}\left(U_k^\top X_k\right) - \left\|U_k^\top \left(X_k - D_{\hat{I}_k^{(0)}} X_k D_{\hat{I}_{-k}^{(0)}}\right)\right\| - \left\| D_{\hat{I}_{k}^{(0)}}Z_{k} D_{\hat{I}_{-k}^{(0)}}\right\|_F\\
				\overset{\eqref{ineq:A-6}}{\geq} & \sigma_{r_k}\left(X_k\right) - \left\|X_k - D_{I_k^{(0)}} X_k D_{I_{-k}^{(0)}}\right\| - \left\| D_{I_{k}}Z_{k} D_{I_{-k}}\right\|_F\\
				\overset{\eqref{ineq:A-2}}{\geq} & \lambda_k - \left\|\X - \X \times_1 D_{I_1^{(0)}} \times \cdots \times_d D_{I_d^{(0)}}\right\|_F - \left(s + 2\sqrt{s\log p} + 2\log p\right)^{1/2}\\
				\geq & \lambda_k - \left(\sum_{k=1}^d \|D_{I_k \backslash I_k^{(0)}}X_k\|_F^2 \right)^{1/2} - \left(s + 2\sqrt{s\log p} + 2\log p\right)^{1/2}\\
				\overset{\eqref{eq:step-1-index}}{\geq} & \lambda_k - \left(16ds\log p\right)^{1/2} - \left(\sqrt{s} + \sqrt{2\log p}\right)\\
				\geq  & \lambda_{k} - (4+1+\sqrt{2})\sqrt{ds\log p}\\
				\geq & \lambda_k - .1\lambda_k = .9\lambda_k.
			\end{split}
		\end{equation*}
		\item Note that $D_{\hat{I}_{k}^{(0)}}Y_{k} D_{\hat{I}_{-k}^{(0)}} = D_{\hat{I}_{k}^{(0)}}X_{k} D_{\hat{I}_{-k}^{(0)}} + D_{\hat{I}_{k}^{(0)}}Z_{k} D_{\hat{I}_{-k}^{(0)}}$, $\rank(D_{\hat{I}_{k}^{(0)}}X_{k} D_{\hat{I}_{-k}^{(0)}})\leq r_k$, we know
		\begin{equation*}
			\begin{split}
				& \sigma_{r_k+1} \left(D_{\hat{I}_{k}^{(0)}}Y_{k} D_{\hat{I}_{-k}^{(0)}}\right) \overset{\text{Lemma \ref{lm:projection-X-residual}}}{\leq} \sigma_{1}\left(D_{\hat{I}_{k}^{(0)}}Z_{k} D_{\hat{I}_{-k}^{(0)}}\right) \overset{\eqref{ineq:A-6}}{\leq} \sigma_1\left(D_{I_{k}}Z_{k} D_{I_{-k}}\right) \\
				\overset{\eqref{ineq:A-3}}{\leq} & \sqrt{s_k} + \sqrt{s_{-k}} + 2\sqrt{\log p} \leq 4\sqrt{s\log p} \overset{\eqref{ineq:lambda_k_min}}{\leq} .1\lambda_k.
			\end{split}
		\end{equation*}
		\item Since the left singular space of $D_{I_k}X_{k} D_{I_{-k}}$ is $U_k$, we have $U_{k\perp}^\top D_{I_k}X_kD_{I_{-k}}=O$, thus,
		\begin{equation*}
			\begin{split}
				& \left\|U_{k\perp}^\top D_{\hat{I}_k^{(0)}}Y_{k} D_{\hat{I}_{-k}^{(0)}}\right\|_F \overset{\eqref{ineq:A-6}}{\leq} \left\|U_{k\perp}^\top D_{\hat{I}_k^{(0)}}Y_{k} D_{I_{-k}}\right\|_F \\
				\leq &  \left\|U_{k\perp}^\top D_{I_k}Y_{k} D_{I_{-k}}\right\|_F + \left\|U_{k\perp}^\top D_{I_k\backslash \hat{I}_k^{(0)}}Y_{k} D_{I_{-k}}\right\|_F\\
				\leq & \left\|U_{k\perp}^\top D_{I_k}X_{k} D_{I_{-k}} + U_{k\perp}^\top D_{I_k}Z_{k} D_{I_{-k}}\right\|_F + \left\|U_{k\perp}^\top Y_{k, [I_k\backslash \hat{I}_{k}^{(0)}, I_{-k}]}\right\|_F\\
				\leq & \left\| U_{k\perp}^\top D_{I_k}Z_{k} D_{I_{-k}}\right\|_F + \left\|Y_{k, [I_k\backslash \hat{I}_{k}^{(0)}, I_{-k}]}\right\|_F\\
				\overset{\eqref{ineq:A-6}}{\leq} & \left\|Z_{k, [I_k, I_{-k}]}\right\|_F + \left\|X_{k, [I_k\backslash I_{k}^{(0)}, I_{-k}]}\right\|_F + \left\|Z_{k, [I_k\backslash I_{k}^{(0)}, I_{-k}]}\right\|_F\\
				\leq & 2\left\|Z_{k, [I_k, I_{-k}]}\right\|_F  + \sqrt{s_{k}16s_{-k}\log p}\\
				\overset{\eqref{ineq:A-2}}{\leq} & 2\sqrt{s+2\sqrt{s\log p}+2\log p} + 4\sqrt{s\log p} \\
				\leq & 2(\sqrt{s}+\sqrt{2\log p}) + 4\sqrt{s\log p} \\
				\leq & 9\sqrt{s\log p}.
			\end{split}
		\end{equation*}
	\end{enumerate}
	Summarizing (a), (b), (c), and \eqref{ineq:sine-theta-U^0}, we must have \eqref{ineq:step2-result-1}, i.e.
	\begin{equation*}
		\begin{split}
			\left\|\sin\Theta\left(\hat{U}_k^{(0)}, U_k\right)\right\|_F\leq \frac{12\sqrt{s\log p}}{\lambda_k}, \quad \forall k=1,\ldots, d,
		\end{split}
	\end{equation*}
	provided that ($A_1$) -- ($A_7$) hold. Next we consider the bound for $\|(\hat{U}_{k\perp}^{(0)})^\top X_k\|$. Since $\hat{U}_k^{(0)}$ is the leading $r_k$ left singular vectors of $\mathcal{M}_k(\tilde{Y}^{(0)}) = D_{I_k^{(0)}}Y_kD_{I_{-k}^{(0)}}$, $\rank(D_{I_k^{(0)}}X_{k}D_{I_{-k}^{(0)}})\leq r_k$, and
	$$D_{I_k^{(0)}}Y_{k}D_{I_{-k}^{(0)}} = D_{I_k^{(0)}}X_{k}D_{I_{-k}^{(0)}} + D_{I_k^{(0)}}Z_{k}D_{I_{-k}^{(0)}}.$$
	Then by Lemma \ref{lm:projection-X-residual},
	\begin{equation}\label{ineq:step1-inter-1}
		\begin{split}
			& \left\|(\hat{U}_k^{(0)})^\top  D_{I_k^{(0)}}X_{k}D_{I_{-k}^{(0)}}\right\|_F \leq 2\sqrt{r_k} \left\| D_{I_k^{(0)}}Z_{k}D_{I_{-k}^{(0)}}\right\|\leq 2\sqrt{r_k}\left\|Z_{[I_k, I_{-k}]}\right\|\\
			\overset{\eqref{ineq:A-3}}{\leq} & 2\sqrt{r_{k}}(\sqrt{s_k} + \sqrt{s_{-k}} + 2\sqrt{\log p}) \leq 8\sqrt{s\log p}.
		\end{split}
	\end{equation}
	The last inequality comes from the fact that $r_ks_k \leq s$, $r_ks_{-k} \leq s$ and $r_k \leq s_k$. Meanwhile,
	\begin{equation}\label{ineq:step1-inter-2}
		\begin{split}
			& \left\|(\hat{U}_k^{(0)})^\top \left(D_{I_k^{(0)}}X_{k}D_{I_{-k}^{(0)}} - X_k \right)\right\|_F\\
			\leq & \left\|X_k -  D_{I_k^{(0)}}X_{k}D_{I_{-k}^{(0)}}\right\|_F = \left\|\X - \X_{[I_1^{(0)}, \ldots, I_d^{(0)}]}\right\|_F  \leq \sqrt{\sum_{k=1}^d \|X_{[I_k \backslash I_k^{(0)}, :]}\|_F^2}\\
			= & \sqrt{\sum_{k=1}^d \sum_{i \in I_k \backslash I_k^{(0)}}\left\|X_{k, [i, :]}\right\|_2^2} \leq \sqrt{\sum_{k=1}^d s_k \cdot 16s_{-k} \log p} \\
			= &  4 \sqrt{ds \log p}.
		\end{split}
	\end{equation}
	Therefore,
	\begin{equation*}
		\begin{split}
			\left\|(\hat{U}_k^{(0)})^\top X_k \right\|_F \leq & \left\|(\hat{U}_k^{(0)})^\top  D_{I_k^{(0)}}X_{k}D_{I_{-k}^{(0)}}\right\|_F + \left\|(\hat{U}_k^{(0)})^\top \left(X_k -  D_{I_k^{(0)}}X_{k}D_{I_{-k}^{(0)}}\right)\right\|_F\\
			\overset{\eqref{ineq:step1-inter-1}\eqref{ineq:step1-inter-2}}{\leq} & 8\sqrt{s\log p} + 4\sqrt{ds \log p} \leq 12\sqrt{ds\log p},
		\end{split}
	\end{equation*}
	which has proved \eqref{ineq:step2-result-2}.
	\item[Step 3] Next we consider the refinement of each iteration. To be specific, we try to study the performance of $\hat{U}_k^{(t)}$ based on $\hat{U}_k^{(t-1)}$. We still assume \eqref{ineq:A-1} -- \eqref{ineq:A-6} all hold. Based on the result in Step 2, we have
	\begin{equation*}
		E_k^{(0)} \leq 12\sqrt{ds\log p} \quad \text{and}\quad F_k^{(0)} \leq \frac{12\sqrt{s\log p}}{\lambda_k},
	\end{equation*}
	provided that ($A_1$)--($A_7$) hold. We particularly provides the following upper bound for $E_k^{(t)}=\left\|(\hat{U}_{k\perp}^{(t)})^\top X_k\right\|_F$ and $F_{k}^{(t)} = \left\| \sin\Theta\left(\hat{U}_k^{(k)}, U_k\right)\right\|_F$
	\begin{equation}\label{ineq:E_new_upper}
	E_k^{(t)} \leq \frac{\sqrt{s_k\bar{\eta}_k} + 3\sqrt{r_k}N_k + 3\sqrt{r_k}M_k\left(K_k^{(t)} + \frac{E_{k+1}^{(t-1)}}{\lambda_{k+1}} + \cdots + \frac{E_{d}^{(t-1)}}{\lambda_{d}} + \frac{E_{1}^{(t)}}{\lambda_{1}} + \cdots + \frac{E_{k-1}^{(t)}}{\lambda_{k-1}}\right)}{1 - K_k^{(t)}} ,
	\end{equation}
	\begin{equation}\label{ineq:F_new_upper}
	F_{k}^{(t)} \leq \frac{E_k^{(t)}}{\lambda_k}.
	\end{equation}
	The detailed proof is collected in Section \ref{sec:step-3} in the supplementary materials.

	\item[Step 4] In this step, we combine the results in Step 4 and provide a upper bound for $E_k^{(t_{\max})}, F_k^{(t_{\max})}$ for $t_{\max} = C\left(\log(ds\log p)\right)$. We first show, for all sparse and non-sparse modes, we have the following upper bounds for $E_k^{(t)}$:
	\begin{equation}\label{ineq:step-5-induction}
		E_k^{(t)} \leq 30\sqrt{s_kr_k} + 30\sqrt{s_k\log p} + \frac{12\sqrt{ds\log p}}{2^t},\quad k=1,\ldots, d, \quad t=0,1,2\ldots.
	\end{equation}
 First, \eqref{ineq:step-5-induction} holds for $t=0$ due to \eqref{ineq:step2-result-2}. If \eqref{ineq:step-5-induction} holds for $t-1$ for some $t\geq 1$, we aim to prove \eqref{ineq:step-5-induction} for $t$. First, based on the basic principle of Tucker rank, we have
	\begin{equation}\label{ineq:r_k-r_-k}
		r_k \leq r_{-k}, \quad r_k \leq s_k \quad \Rightarrow \quad s_l r_l \leq s_l r_{-l} \leq s_l s_{-l} = s.
	\end{equation}
	Also, we shall recall that
	\begin{equation}\label{ineq:eta_N_M}
		\begin{split}
			&  \sqrt{\eta_k} = \sqrt{r_{-k} + 2\sqrt{r_{-k}\log p} + 2\log p} \leq \sqrt{r_{-k}} + \sqrt{2\log p};\\
			& \sqrt{\bar{\eta}_k} = \sqrt{r_{k} + 2\sqrt{r_{k}\log p} + 2\log p} \leq \sqrt{r_{k}} + \sqrt{2\log p};\\
			& N_l \leq \sqrt{s_l} + \sqrt{r_l} + 2\sqrt{\log p} \leq 2\sqrt{s_l} + 2\sqrt{\log p};\\
			& M_l \leq 2\left(\sqrt{s_l} + \sqrt{r_{-l}} + \sum_{l=2}^d \sqrt{s_lr_l} + \sqrt{\log p}\right),
		\end{split}
	\end{equation}
	and
	\begin{equation}\label{ineq:E_k^{(t-1)}-inequality}
		\begin{split}
			& E_k^{(t-1)} \leq 30\sqrt{s_kr_k} + 30\sqrt{s_k\log p} + \frac{12\sqrt{ds\log p}}{2^{t-1}} \\
			\leq & 30\sqrt{s} + 30\sqrt{s_k\log p} + \frac{12\sqrt{ds\log p}}{2^{t-1}}\\
			\leq & 60 \sqrt{s\log p} + \frac{12\sqrt{ds\log p}}{2^{t-1}}.
		\end{split}
	\end{equation}
	Thus, the following upper bound holds for $K_1^{(t)}$, if we set $C_{gap} \geq 200d$.
	\begin{equation*}
		\begin{split}
			K_1^{(t)} \overset{\eqref{ineq:K_t}}{\leq} & \left(\frac{(E_{2}^{(t-1)})^2 + \cdots + (E_{d}^{(t-1)})^2 + 2s_1\eta_1 + 6r_1M_1^2}{\lambda_1^2}\right)^{1/2}\\
			\overset{\eqref{ineq:eta_N_M}\eqref{ineq:E_k^{(t-1)}-inequality}}{\leq} & \frac{\left((d-1)\left(60\sqrt{s\log p} + \frac{12\sqrt{ds\log p}}{2^{t-1}}\right)^2 + 2s_1\eta_1 + 6r_1M_1^2\right)^{1/2}}{\lambda_1}\\
			\leq & \frac{60\sqrt{d-1}\sqrt{s\log p} + 12\sqrt{d-1}\sqrt{ds\log p} + \sqrt{2s_1r_{-1}} + 2\sqrt{s_1\log p}}{\lambda_1}\\
			& + \frac{2\sqrt{6}\sqrt{r_1}\left(\sqrt{s_1}+\sqrt{r_{-1}}+\sum_{k=2}^{d}\sqrt{s_kr_k}+\sqrt{\log{p}}\right)}{\lambda_1}\\
			\leq & \frac{100d \lambda_1}{C_{gap}\lambda_1} \leq \frac{1}{2}.
		\end{split}
	\end{equation*}
	Thus,
	\begin{equation}\label{ineq:E_1^t}
		\begin{split}
			E_1^{(t)} \leq & \frac{\sqrt{s_1\bar{\eta}_1} + 3\sqrt{r_1}N_1 + \sqrt{r_1}M_1\left(K_1^{(t)} + \frac{E_2^{(t-1)}}{\lambda_2} + \cdots + \frac{E_d^{(t-1)}}{\lambda_d}\right)}{1 - K_1^{(t)}}\\
			\leq & 2\sqrt{s_1\bar{\eta}_1} + 6\sqrt{r_1}N_1 + 2\sqrt{r_1}M_1\left(K_1^{(t)} + \frac{E_2^{(t-1)}}{\lambda_2} + \cdots + \frac{E_d^{(t-1)}}{\lambda_d}\right)\\
			\overset{\eqref{ineq:eta_N_M}}{\leq} & 14\sqrt{s_1r_1} + 15\sqrt{s_1\log p} + 2\sqrt{r_1}M_1\left(K_1^{(t)} + \frac{E_2^{(t-1)}}{\lambda_2} + \cdots + \frac{E_d^{(t-1)}}{\lambda_d}\right).
		\end{split}
	\end{equation}
	With $C_{gap}>20(d+2)(d-1)$, by Lemma \ref{lm:rate-analysis}, we have:
	\begin{equation}\label{ineq:r_1M_1K_1}
		\begin{split}
			 \sqrt{r_1}M_1K_1^{(t)} & \leq \frac{\sqrt{r_1}M_1\sum_{k=2}^dE_k^{(t-1)}}{\lambda_1} + \frac{M_1\sqrt{2s_1r_1\eta_1}}{\lambda_1} + \frac{\sqrt{6}r_1M_1^2}{\lambda_1},\\
			& \leq 4\sqrt{s_1r_1}+4\sqrt{s_1\log p} + \frac{12\sqrt{ds\log p}}{2^{(t+2)}}.
		\end{split}
	\end{equation}
	\begin{equation}\label{ineq:r_1M_1E...E}
		\begin{split}
			 \sqrt{r_1}M_1\left(\frac{E_2^{(t-1)}}{\lambda_2} + \cdots + \frac{E_d^{(t-1)}}{\lambda_d}\right) & \leq 3\sqrt{s_1r_1}+3\sqrt{s_1\log p} + \frac{12\sqrt{ds\log p}}{2^{(t+2)}}
		\end{split}
	\end{equation}
	
	Combining \eqref{ineq:E_1^t}, \eqref{ineq:r_1M_1K_1}, and \eqref{ineq:r_1M_1E...E}, we have proved
	\begin{equation*}
		\begin{split}
			E_k^{(t)} \leq 30\sqrt{s_kr_k} + 30\sqrt{s_k\log p} + \frac{12\sqrt{ds\log p}}{2^t},
		\end{split}
	\end{equation*}
	i.e. \eqref{ineq:step-5-induction} for $t$ and $k=1$. Similarly, one can also prove the upper bounds for $E_k^{(t)}$ for $k=2,\ldots, d$, which implies the claim \eqref{ineq:step-5-induction} holds for $t \geq 0$.
	
	Then we further provide another upper bound for non-sparse modes. Note that we assume $\log {p_1} \asymp \log{p_2} \asymp \cdots \asymp \log{p_d}$, thus we can find some constant $C$, such that $\sqrt{\log p} \leq C\sqrt{p_k},~k=1,2,\ldots,d$. Now we want to show:
	\begin{equation}\label{nonsparse-bound}
	E_k^{(t)} \leq (30+15\sqrt{2}C)\sqrt{p_kr_k}  + \frac{12\sqrt{ds\log p}}{2^t},~~~k \not \in J_s
	\end{equation}
	 Again, we assume mode 1 is non-sparse and only prove the bound for $E_1^{(t)}$, the similar bounds for other non-sparse modes essentially follow. \\
	Return to \eqref{ineq:E_1^t}, note that we particularly have $\bar{\eta}_1=0$ since it is a non-sparse mode, we have:
	\begin{equation}\label{reboundE}
		\begin{split}
		E_1^{(t)} & \leq  6\sqrt{r_1}N_1 + 2\sqrt{r_1}M_1\left(K_1^{(t)} + \frac{E_2^{(t-1)}}{\lambda_2} + \cdots + \frac{E_d^{(t-1)}}{\lambda_d}\right) \\
		& \overset{\eqref{ineq:eta_N_M}}{\leq} 12\sqrt{p_1r_1} + 6\sqrt{2r_1\log p} + 2\sqrt{r_1}M_1\left(K_1^{(t)} + \frac{E_2^{(t-1)}}{\lambda_2} + \cdots + \frac{E_d^{(t-1)}}{\lambda_d}\right)\\
		& \leq (12+6\sqrt{2}C)\sqrt{p_1r_1} + 2\sqrt{r_1}M_1\left(K_1^{(t)} + \frac{E_2^{(t-1)}}{\lambda_2} + \cdots + \frac{E_d^{(t-1)}}{\lambda_d}\right).
		\end{split}
	\end{equation}
	Also, we have the following bound for $K_1^{(t)}$ as $\eta_1=0$,
	\begin{equation}\label{reboundK}
		\begin{split}
		K_1^{(t)} &\leq \frac{(E_2^{(t-1)})+\ldots+(E_d^{(t-1)})+\sqrt{6r_1}M_1}{\lambda_1}. \\
		\end{split}
	\end{equation}
	Besides, we could rebound $M_1$ like following:
\begin{equation}\label{reboundM}
		\begin{split}
		M_1 &\leq 2\left((1+C)\sqrt{p_1}+\sqrt{r_{-1}}+\sum_{l\geq2}\sqrt{s_lr_l}\right). \\
		\end{split}
	\end{equation}
	If we set $C_{gap}>20(d+2)(d-1)$, then together with \eqref{ineq:step-5-induction}, \eqref{reboundE}, \eqref{reboundK} and \eqref{reboundM}, we can proved the following result as similar as we proved in Lemma \ref{lm:rate-analysis}:
	\begin{equation*}
	2\sqrt{r_1}M_1\left(K_1^{(t)} + \frac{E_2^{(t-1)}}{\sqrt{r_2}\lambda_2} + \cdots + \frac{E_d^{(t-1)}}{\sqrt{r_d}\lambda_d}\right) \leq (18+9\sqrt{2}C)\sqrt{p_1r_1} + \frac{12\sqrt{ds\log p}}{2^t}.
	\end{equation*}
	Then (63) follows.\\
	Now for $t_{\max} = C_1\log (ds\log p)$, for large constant $C_1>0$, we have
	\begin{equation}\label{ineq:E_k-upper-bound}
		\begin{split}
			E_k^{(t_{\max})} = \left\|(\hat{U}_k^{(t_{\max})})^\top X_k\right\|_F \leq \left\{\begin{array}{ll}
	C\left(\sqrt{s_kr_k} + \sqrt{s_k \log p}\right), &  k \in J_s,\\
	C\sqrt{p_kr_k}, & k\notin J_s,
	\end{array}\right.		\end{split}
	\end{equation}
	\begin{equation}\label{ineq:F_k-upper-bound}
		\begin{split}
			F_k^{(t_{\max})} = \left\|\sin\Theta\left(\hat{U}_k^{(t_{\max})}, U_k\right)\right\|_F \overset{\eqref{ineq:U-U-sin-theta}}{\leq}  \left\{\begin{array}{ll}
	C\left(\sqrt{s_kr_k} + \sqrt{s_k \log p}\right)/\lambda_k, &  k \in J_s,\\
	C\sqrt{p_kr_k}/\lambda_k, & k\notin J_s.
	\end{array}\right.
		\end{split}
	\end{equation}
	
	\item[Step 5] In this step, we consider the estimation error for $\hat{\X}$. We first have the following decomposition,
	\begin{equation*}
		\begin{split}
			& \left\|\hat{\X} - \X \right\|_F = \left\|\Y \times_1 P_{\hat{U}_1} \cdots \times_d P_{\hat{U}_d} - \X\right\|_F\\
			\leq & \left\|\Z \times_1 P_{\hat{U}_1} \cdots \times_d P_{\hat{U}_d}\right\|_F + \left\|\X \times_1 P_{\hat{U}_1} \cdots \times_d P_{\hat{U}_d} - \X\right\|_F.\\
		\end{split}
	\end{equation*}
	Here,
	\begin{equation*}
		\begin{split}
			& \left\|\X \times_1 P_{\hat{U}_1} \cdots \times_d P_{\hat{U}_d} - \X\right\|_F \\
			\leq & \left\|\X\times_1 P_{\hat{U}_{1\perp}} + \X\times_1 P_{\hat{U}_{1}} \times_2 P_{\hat{U}_{2\perp}} + \cdots  + \X\times_1 P_{\hat{U}_{1}} \times_2 P_{\hat{U}_{2}} \times \cdots \times_{(d-1)} P_{\hat{U}_{d-1}}\times_d P_{\hat{U}_{d}} \right\|_F\\
			\leq & \sum_{l=1}^d \left\|\X\times_l (\hat{U}_{l\perp})^\top \right\|_F = \sum_{l=1}^d \left\|\hat{U}_{l\perp}^\top X_l\right\|_F \overset{\eqref{ineq:E_k-upper-bound}}{\leq} C\left(\sum_{l\in J_s} (\sqrt{s_lr_l} + \sqrt{s_l\log p})+\sum_{l \not \in J_s}\sqrt{p_lr_l}\right);
		\end{split}
	\end{equation*}
	\begin{equation*}
		\begin{split}
			& \left\|\Z\times_1 P_{\hat{U}_1} \times \cdots \times \times_d P_{\hat{U}_d}\right\|_F\\
			\leq & \left\|\Z\times_1 \left(U_1U_1^\top \hat{U}_1\right)^\top \times_2 \hat{U}_2^\top \cdots \times \times_d \hat{U}_d^\top \right\|_F + \left\|\Z\times_1 \left(U_{1\perp}U_{1\perp}^\top \hat{U}_1\right)^\top \times_2 \hat{U}_2^\top \cdots \times \times_d \hat{U}_d^\top \right\|_F\\
			\leq & \left\|\Z\times_1 U_1^\top \times_2 \hat{U}_2^\top \cdots \times \times_d \hat{U}_d^\top \right\|_F + \left\|\left(U_{1\perp}U_{1\perp}^\top \hat{U}_1\right)^\top Z_1 \left(\hat{U}_2 \otimes \cdots \otimes \hat{U}_d^\top \right)\right\|_F\\
			\leq & \left\|\Z\times_1 U_1^\top \times_2 \hat{U}_2^\top \cdots \times \times_d \hat{U}_d^\top \right\|_F + M_1 \left\|\sin\Theta\left(\hat{U}_1, U_1\right)\right\|_F\\
			\leq & \cdots \\
			\leq & \left\|\Z\times_1 U_1^\top \times_2 U_2^\top \cdots \times \times_d \hat{U}_d^\top \right\|_F + M_1 \left\|\sin\Theta\left(\hat{U}_1, U_1\right)\right\|_F + M_2 \left\|\sin\Theta\left(\hat{U}_2, U_2\right)\right\|_F\\
			\leq & \cdots\\
			\leq & \left\|\Z\times_1 U_1^\top \times_2 \cdots \times_d U_d^\top\right\|_F + \sum_{l=1}^d M_l \left\|\sin\Theta(\hat{U}_l, U_l)\right\|_F.
		\end{split}
	\end{equation*}
	\begin{equation*}
		\left\|\Z \times_1 U_1^\top \times \cdots \times_d U_d^\top \right\| \overset{\eqref{ineq:A-2}}{\leq} \sqrt{r + 2\sqrt{r\log p} + 2\log p} = \sqrt{r} + \sqrt{2\log p};
	\end{equation*}
	\begin{equation*}
		\begin{split}
			\sum_{l=1}^{d}M_l \left\|\sin\Theta\left(\hat{U}_l, U_l\right)\right\|_F \overset{\eqref{ineq:A-5}\eqref{ineq:F_k-upper-bound}}{\leq} & C\sum_{l\in J_s}\frac{M_l}{\lambda_l}(\sqrt{s_lr_l}+ \sqrt{s_l\log p}) + C\sum_{l\not \in J_s}\frac{M_l}{\lambda_l}(\sqrt{p_lr_l})\\
			\leq & C\sum_{l\in J_s}(\sqrt{s_lr_l}+ \sqrt{s_l\log p}) + C\sum_{l\not \in J_s}(\sqrt{p_lr_l}).
		\end{split}
	\end{equation*}
	In summary,
	\begin{equation*}
		\left\|\hat{\X} - \X \right\|_F \leq C\left(\sqrt{r} + \sum_{l=1}^d \sqrt{r_ls_l} + \sum_{l \in J_s}\sqrt{s_l \log p}\right).
	\end{equation*}
	under the assumptions \eqref{ineq:A-1} -- \eqref{ineq:A-6} that we introduced in Step 1. Note that $\log p = \log p_1 + \cdots + \log p_d$ and we have the assumption that $\log{p_1} \asymp \cdots \asymp \log{p_d}$, thus the result is equitant to the one stated in Theorem 1.
	
	\item[Step 6] Finally, we prove that all imposed conditions in Step 1, i.e. \eqref{ineq:A-1}-\eqref{ineq:A-6}, hold with probability at least $1-O(\frac{(p_1+\cdots + p_d) \log{(ds\log p)}}{p})$. The detailed proof is collected in Section \ref{sec:step-6} in the supplementary materials. 
\end{enumerate}	
Combing all results from Steps 1 -- 6, we have finished the proof for this theorem. \quad $\square$

%% file: stat-svd7-appendix.tex

\begin{abstract}
	In this supplement, we provide discussions for the more straightforward single thresholding \& projection scheme and additional proofs of the main theorems. The key technical tools used in the proofs of the main technical results are also introduced and proved.
\end{abstract}

\section{Comparison with Single Thresholding \& Projection Scheme}

As illustrated in Section \ref{sec:method}, the update of $U_k$ along sparse mode is the core step in the proposed STAT-SVD algorithm. In order to allow more accurate thresholding, we propose a novel double thresholding \& projection scheme (Algorithm \ref{al:refinement-sparse}) in the main body of the paper. Compared with the proposed scheme, the following single thresholding \& projection scheme would be a more straightforward idea, as it can be seen as a simpler extension from matrix sparse SVD \citep{lee2010biclustering,yang2016rate}.
\begin{algorithm}
	\caption{An Alternative Sparse Mode Update Scheme -- Single Thresholding \& Projection}
	\begin{algorithmic}[1]
		\State Input: $\Y$, $\hat{U}_{k+1}^{(t)}, \ldots, \hat{U}_d^{(t)}, \hat{U}_1^{(t+1)},\ldots, \hat{U}_{k-1}^{(t+1)}$, rank $r_k$.
		\State (Projection) Calculate
		$$A_k^{(t+1)} = \mathcal{M}_k\left(\llbracket \Y;  (\hat{U}_1^{(t+1)})^\top,  \ldots, (\hat{U}_{k-1}^{(t+1)})^\top, I_{p_k}, (\hat{U}_{k+1}^{(t)})^\top, \ldots, (\hat{U}_{d}^{(t)})^\top \rrbracket \right).$$
		\State (Thresholding) Perform row-wise thresholding for $A_k^{(t+1)}$:
		$$B_k^{(t+1)} \in \mathbb{R}^{p_k\times r_{-k}}, \quad B^{(t+1)}_{k, [i, :]} = A^{(t+1)}_{k, [i, :]} 1_{\left\{\|A^{(t+1)}_{k, [i, :]}\|_2^2 \geq \eta_k\right\}}, \quad 1\leq i \leq p_k,$$
		\quad\quad \quad where $\eta_k = r_{-k}+ 2\left(\sqrt{r_{-k} \log p_k} + \log p_k\right)$.
		\State Return $\hat{U}_k^{(t+1)} = \SVD_{r_k}(B^{(t+1)})$.
	\end{algorithmic}
	\label{al:refinement-sparse-single-thresholding}
\end{algorithm}

Suppose $\tilde{\X}$ is the denoising estimator for $\X$ if one applies Algorithm \ref{al:procedure} with single thresholding \& projection (Algorithm \ref{al:refinement-sparse-single-thresholding}) instead of the double one (Algorithm \ref{al:refinement-sparse}), $\tilde{U}_k$ and $\tilde{\X}$ may yield a higher rate of convergence in general compared with the proposed STAT-SVD method as illustrated by the following theorem.
\begin{Theorem}\label{th:single-thresholding-projection}
	Suppose $s_k \geq 2r_k$ and $s = o(p^{1/4}/log p)$. There exists a constant $c>0$ and a parameter tensor $\X\in \mathcal{F}_{\bp, \br, \bs}(J_s)$, such that even after sufficient number of iterations that is required in Theorem \ref{th:upper_bound}, the single thresholding \& projection estimator $\tilde{\X}$ yields the following lower bound
	\begin{equation*}
	\mathbb{E}\|\tilde{\X} - \X\|_F^2 \geq c\sigma^2 \left(\prod_{k}r_k + \sum_{k} s_kr_k + \sum_{k\in J_s} s_k\log p_k + \sum_{k\in J_s} s_k (r_{-k}\log p)^{1/2}\right).
	\end{equation*}
\end{Theorem}

{\bf\noindent Proof of Theorem \ref{th:single-thresholding-projection}.} Without loss of generality we can assume that $\sigma^2 = 1$. By the lower bound argument in Theorem \ref{th:lower_bound_recovery}, we only need to show
\begin{equation*}
\mathbb{E}\|\tilde{\X} - \X\|_F^2 \geq c\sigma^2 s_k\sqrt{r_{-k}\log p}, \quad \forall k \in J_s.
\end{equation*}
for some $\X \in \mathcal{F}_{\bp, \br, \bs}(J_s)$ and constant $c>0$. Without loss of generality, for the rest of the proof, we focus on Mode-1 and aim to show $\mathbb{E}\|\tilde{\X} - \X\|_F^2 \geq c\sigma^2 s_1\sqrt{r_{-1}\log p}$ when $1 \in J_s$. Based on similar argument as the one of Theorem \ref{th:lower_bound_subspace}, we can construct $\tilde{\S}\in \mathbb{R}^{r_1\times \cdots\times r_d}$ with i.i.d. Gaussian entries. Then there exists constants $c_1, C_1>0$ such that
\begin{equation*}
\begin{split}
c_1 \sqrt{r_{-k}} \leq \sigma_{\min}(\mathcal{M}_k(\tilde{\S})) \leq \sigma_{\max}(\mathcal{M}_k(\tilde{\S})) \leq C_1 \sqrt{r_{-k}} \quad \text{and} \quad \|\tilde{\S}\|_\infty \leq C_1\sqrt{\log p}.
\end{split}
\end{equation*}
happens with a positive probability. Suppose $\tilde{\S}$ satisfies the previous condition, we rescale as $\S = \tilde{\S}\cdot \lambda / \min_{k}\sigma_{r_k} (\mathcal{M}_k(\tilde{\S}))$ to ensure that $\sigma_{\min}(\mathcal{M}_k(\S)) \geq \lambda$. By such construction, we can see $\S$ satisfies the following inequality,
\begin{equation}\label{ineq:th4-cond1}
\|\S\|_\infty\cdot \sqrt{\frac{r_{-k}}{\log p}} \leq C_2 \sigma_{\min}(\mathcal{M}_k(\tilde{\S})) \leq C_2\sigma_{\max}\left(\mathcal{M}_k(\S)\right) \leq C_2^2\sigma_{\min}(\mathcal{M}_k(\tilde{\S})), \quad 1\leq k \leq d.
\end{equation}
Define
\begin{equation*}
U_1 = \begin{bmatrix}
I_{r_1-1} & 0_{(r_1-1)\times 1} \\
0_{1\times (r_1-1)} & \sqrt{1 - (s_1-r_1)\tau^2} \\
0_{(s_1-r_1)\times (r_1-1)} & \tau \cdot 1_{(s_1-r_1)\times 1} \\
0_{(p_1-s_1) \times (r_1-1)} & 0_{(p_1-s_1)\times 1}
\end{bmatrix},
\end{equation*}
\begin{equation*}
\begin{split}
U_k = \begin{bmatrix}
I_{r_k}\\
0_{p_k - r_k, r_k}
\end{bmatrix}, \quad k=2,\ldots, d.
\end{split}
\end{equation*}
Here, $1_{a\times b}$ and $0_{a\times b}$ are the $a$-by-$b$ matrices with all ones and zeros, respectively;
\begin{equation}\label{ineq:th4-intermediate-0}
\tau^2 = \frac{\sqrt{r_{-1}\log p}}{(9C_2 \vee 3) \cdot\|\mathcal{M}_1(\S)\|^2},
\end{equation}
is a scaling factor. Under such a configuration, we have
\begin{equation}\label{ineq:th4-intermediate-1}
\begin{split}
\forall r_1+1\leq i \leq s_1, \quad \|\X_{1,[i, :]}\|^2 = & \| (U_1)_{[i, :]} \cdot \mathcal{M}_1(\S) \cdot (U_2\otimes \cdots \otimes U_d)^\top\|_2^2 \\
\leq & \tau^2 \cdot \|\mathcal{M}_1(\S)\|^2 \overset{\eqref{ineq:th4-intermediate-0}}{\leq} \frac{\sqrt{r_{-1}\log p}}{3};
\end{split}
\end{equation}
\begin{equation*}
\begin{split}
\forall r_1+1\leq i \leq s_1, \quad \|\X_{1,[i, :]}\|_\infty = & \left\|(U_1)_{[i, :]}[\mathcal{M}_1(\S)]\cdot (U_2\otimes \cdots \otimes U_d)^\top\right\|_\infty\\
= & \tau\left\|[\mathcal{M}_1(\S)]_{[i, :]}\cdot (U_2\otimes \cdots \otimes U_d)^\top\right\|_\infty\\
\leq & \tau\|\S\|_\infty \overset{\eqref{ineq:th4-intermediate-0}}{\leq} \sqrt{\frac{\sqrt{r_{-1}\log p}}{9C_2}}\cdot \frac{\|\S\|_\infty}{\|\mathcal{M}_1(\S)\|} \overset{\eqref{ineq:th4-cond1}}{\leq} \frac{1}{3}\sqrt{\log p}.
\end{split}
\end{equation*}
Recall $\tilde{\X}$ is the single thresholding \& projection estimator after $t_{\max}$ iterations. Next we show that with probability at least $1 - Cst_{\max}/p$, the Mode-1 support of $\mathcal{M}_1(\tilde{\X})$ does not cover the third block of $U_1$, i.e.,
\begin{equation}\label{eq:bad-event}
\supp(\mathcal{M}_1(\tilde{\X})) \cap \{r_1 + 1,\ldots, s_1\} = \emptyset.
\end{equation}
In order to prove \eqref{eq:bad-event}, we analyze the procedure of STAT-SVD with single thresholding \& projection, and show that with high probability, all Mode-1 indices in $\{r_1+1,\ldots, s_1\}$ will be thresholded in each round of iteration, i.e.,
\begin{equation}\label{eq:hat_I_1_bad}
\hat{I}_1^{(t)} \cap \supp(\mathcal{M}_1(\tilde{\X})) = \emptyset, \quad \forall 1\leq t \leq t_{\max}.
\end{equation}
Similarly to the Step 6 in the proof of Theorem \ref{th:upper_bound}, we construct another parallel sequence of estimations as follows.
\begin{enumerate}[(a)]
	\item Initialize
\begin{equation}
\hat{\cI}_k^{(0)} = \left\{\begin{array}{ll}
\Bigg\{i \in I_k: \left\|(Y_k)_{[i, :]}\right\|_2^2 \geq \sigma^2 \left(p_{-k} + 2\sqrt{p_{-k}\log p} + 2 \log p\right)\\
\quad\quad\quad \text{or } \max_{j} |(Y_k)_{[i,j]}| \geq 2\sigma\sqrt{\log p}\Bigg\} \backslash \{r_1+1,\ldots, s_1\}, & k = 1;\\
\Bigg\{i \in I_k: \left\|(Y_k)_{[i, :]}\right\|_2^2 \geq \sigma^2 \left(p_{-k} + 2\sqrt{p_{-k}\log p} + 2 \log p\right)\\
\quad\quad\quad \text{or } \max_{j} |(Y_k)_{[i,j]}| \geq 2\sigma\sqrt{\log p}\Bigg\}, & k\neq 1, k \in J_s;\\
\{1,\ldots, p_k\}, & k \notin J_s.
\end{array}
\right.
\end{equation}
$$\tilde{\Y} =\left\{\begin{array}{ll}
\Y_{i_1,\ldots,i_d}, & (i_1,\ldots, i_d)\in \hat{\cI}_1^{(0)}\otimes \cdots \otimes \hat{\cI}_d^{(0)};\\
0, & \text{otherwise};
\end{array}\right. \quad \hat{\cU}^{(0)}_k = \SVD_{r_k}\left(\mathcal{M}_k(\tilde{\Y})\right) \in \mathbb{R}^{p_k\times r_k}.$$
\item For $t = 1,\ldots, t_{\max}$, $k=1,\ldots, d$, let
$$\cA_k^{(t)} = \mathcal{M}_k\left(\Y \times_1 (\hat{\cU}_1^{(t)})^\top \times \cdots \times_{k-1} (\hat{\cU}_{k-1}^{(t)})^\top \times_{k+1} (\hat{\cU}_{k+1}^{(t-1)})^\top \times \cdots \times_d (\hat{\cU}_{d}^{(t-1)})^\top \right)\in \mathbb{R}^{p_k\times r_{-k}};$$
\begin{equation*}
\cB_k^{(t)} \in \mathbb{R}^{p_k\times r_{-k}}, \quad \cB^{(t)}_{k, [i, :]} = \left\{\begin{array}{ll}
\cA^{(t)}_{k, [i, :]} 1_{\left\{\|\cA^{(t)}_{k, [i, :]}\|_2^2 \geq \eta_k \text{ and } i\notin \{r_1,\ldots, s_1\}\right\}}, & k=1;\\
\cA^{(t)}_{k, [i, :]} 1_{\left\{\|\cA^{(t)}_{k, [i, :]}\|_2^2 \geq \eta_k \right\}}, & k \neq 1, k\in J_s;\\
\cA^{(t)}_{k, [i, :]}, & k \notin J_s,
\end{array}\right.
\end{equation*}
where $\eta_k = \sigma^2\left(r_{-k} + 2\left(r_{-k} \log p\right)^{1/2} + \log p\right)$. Finally, $U_k$ is updated via
$$\hat{\cU}_k^{(t+1)} = \SVD_{r_k}\left(\cB_k^{(t)}\right) \in \mathbb{O}_{p_k, r_k}.$$
\end{enumerate}
Essentially, $\hat{\cU}_k$ can be seen as the outcome of single projection \& thresholding algorithm performed on $\Y$ without $\{r_1+1, \ldots, s_1\}$-th indices on Mode-1. Next, we show the path of $\hat{\cU}_k^{(t)}$ and $\hat{\cU}^{(t)}$ exactly match with high probability. By comparing the evolution of $\hat{\cU}^{(t)}$ and $\hat{U}^{(t)}$, we only need to evaluate the probability that the following events happen,
$$\hat{\cB}_1^{(t)} = \hat{B}_1^{(t)}, \quad \text{ and } \quad \hat{\cI}_1^{(0)} = \hat{I}_1^{(0)}, $$
which is equivalent to
\begin{equation}\label{ineq:th4-wts2}
\forall 1\leq t\leq t_{\max}, r_1+1\leq i \leq s_1, \quad \|\cA^{(t)}_{1, [i, :]}\|_2^2 \geq r_{-1}+2(r_{-1}\log p)^{1/2} + \log p,
\end{equation}
\begin{equation}\label{ineq:th4-wts1}
\begin{split}
 & \left\|(Y_1)_{[i, :]}\right\|_2^2 < \sigma^2 \left(p_{-1} + 2\sqrt{p_{-1}\log p} + 2 \log p\right), ~~  \|(Y_1)_{[i,:]}\|_\infty < 2\sigma\sqrt{\log p}.
\end{split}
\end{equation}
Similarly as Step 6 of the proof for Theorem \ref{th:upper_bound}, next we evaluate the probability that \eqref{ineq:th4-wts2} and \eqref{ineq:th4-wts1} hold. First, based on the procedure, we know $\cU_k^{(t)}$ is independent of $(Z_1)_{[i, :]}$ for any $i\in \{r_1+1,\ldots, s_1\}$. Thus $\|\mathcal{A}_{1, [i, :]}^{(t)}\|_2^2 \sim \chi^2_{r_{-1}}(\theta)$, where
$$\theta = \left\|(X_1)_{[i, :]}\cdot (\hat{\cU}_2^{(t)}\otimes \cdots \otimes \hat{\cU}_d^{(t)})^\top \right\|_2^2 \leq \left\|(X_1)_{[i, :]}\right\|_2^2\leq \frac{1}{3}\sqrt{r_{-1}\log p}.$$
For specific $t$ and $i \in \{r_1+1,\ldots, s_1 \}$, by Lemma \ref{lm:chi-square-tail-bound}, one has
\begin{equation*}
\begin{split}
& \left\|\cA_{1, [i, :]}^{(t)}\right\|_2^2 < r_{-1} + \theta + 2\left((r_{-1} + 2\theta)\frac{\log p}{4}\right)^{1/2} + \frac{2\log p}{4}\\
\leq & r_{-1} + \frac{\sqrt{r_{-1}\log p}}{3} + 2\sqrt{\left(r_{-1} + \frac{\sqrt{r_{-1}\log p}}{3}\right)\frac{\log p}{4}} + \frac{2\log p}{4}\\
\leq & r_{-1} + \frac{\sqrt{r_{-1}\log p}}{3} + 2\sqrt{r_{-1}\frac{\log p}{4}} + \frac{1}{\sqrt{3}}\left(\sqrt{r_{-1}\log p}\cdot \log p\right)^{1/2} + \frac{2\log p}{4}\\
\leq & r_{-1} + 2\sqrt{r_{-1}\log p} + \log p - \frac{1}{2} \log p - 2/3\sqrt{r_{-1}\log p} + \frac{1}{\sqrt{3}}(\sqrt{r_{-1}\log p}\cdot\log p)^{1/2}\\
\leq & r_{-1}+2\sqrt{r_{-1}\log p} + \log p - \left(\frac{2}{\sqrt{3}} - \frac{1}{\sqrt{3}}\right)(\sqrt{r_{-1}\log p}\cdot\log p)^{1/2}\\
< & r_{-1} + 2(r_{-1}\log p)^{1/2} + \log p,
\end{split}
\end{equation*}
with probability at least $1 - e^{-\frac{\log p}{4}} = 1 - p^{1/4}$. So \eqref{ineq:th4-wts2} holds with probability at least $1 - Cs_1t_{\max}/p^{1/4}$. Similarly one can also show that \eqref{ineq:th4-wts1} holds with probability at least $1 - Cs_1t_{max}p^{1/4}$.

Therefore, within $t_{\max}$ iterations, the paths of $\hat{\cU}^{(t)}$ and $\hat{U}^{(t)}$ exactly match with probability at least $1 - Cs_1 t_{\max}/p^{1/4}$. Finally, if $\hat{\cU}^{(t)}$ and $\hat{U}^{(t)}$ are exactly the same, one has \eqref{eq:bad-event}, \eqref{eq:hat_I_1_bad}, and $\{i: \mathcal{M}_1(\tilde{\X})_{[i, :]}\neq 0\} \cap \{r_1+1,\ldots, s_1\} = \emptyset$. Then, with probability at least $1 - Cs_1t_{\max}/p^{1/4}$, one has
\begin{equation*}
\begin{split}
\|\tilde{\X} - \X\|_F^2 = & \left\|\mathcal{M}_1(\tilde{\X}) - \mathcal{M}_1(\X)\right\|_F^2 \geq \left\|\left[\mathcal{M}_1(\tilde{\X}) - \mathcal{M}_1(\X)\right]_{[(r_1+1): s_1, :]}\right\|_F^2 \\
= & \left\|\left[\mathcal{M}_1(\X)\right]_{[(r_1+1): s_1, :]}\right\|_F^2 = \sum_{i=r_1+1}^{s_1} \|(X_1)_{[i, :]}\|_2^2 \geq (s_1-r_1) \sigma_{r_1}^2(\mathcal{M}_{r_1}(X_1))\\
\geq & \sigma_{r_1}^2(\S)\cdot \tau^2 \cdot (s_1 - r_1) \overset{\eqref{ineq:th4-cond1}}{\geq} cs_1\sqrt{r_{-1}\log p}.
\end{split}
\end{equation*}
Provided that $t_{\max} = O(ds\log p)$ and $s = o(p^{1/4}/\log p)$, we have finished the proof of this lemma. \quad $\square$

\section{Implementation Details of S-HOSVD and S-HOOI}\label{baseline-algorithm}
Here, we summarize the sparse high-order SVD (S-HOSVD) and sparse high-order orthogonal iteration (S-HOOI), which served as two baseline methods in numerical comparison in Section \ref{sec:simu}. Intuitively speaking, S-HOSVD and S-HOOI are sparse modifications of HOSVD and HOOI, with regular SVD replaced by the Sparse SVD in each step. Specifically, let ${\rm SSVD}()$ be the Algorithm 1 in \cite{yang2014sparse}. 
\begin{algorithm}
	\caption{Sparse High-order Singular Value Decomposition (S-HOSVD)}
	\begin{algorithmic}[1]
		\State Input: order-$d$ tensor data $\Y\in \mathbb{R}^{\bp}$, rank $\br$, set of sparse modes $J_s$.
		\For{$k = 1,\ldots, d$}
		\If{$k\in J_s$} 				
		\State $\hat U_k = \SSVD\left(\mathcal M_k(\Y),r_k\right)$
		\Else
		\State $\hat U_k = \SVD_{r_k}\left(\mathcal M_k(\Y)\right)$
		\EndIf
		\EndFor
	\end{algorithmic}\label{al:S-HOSVD}
\end{algorithm}

\begin{algorithm}
	\caption{Sparse High-order Orthogonal Iteration (S-HOOI)}
	\begin{algorithmic}[1]
		\State Input: order-$d$ tensor data $\Y\in \mathbb{R}^{\bp}$, rank $\br$, set of sparse modes $J_s$.
		\State Initialize $\hat U_1^{(0)},\ldots \hat U_d^{(0)}$ as the output of S-HOSVD
		\While{$t < t_{\max}$ or convergence criterion not satisfied}
		\For{$k = 1,\ldots, d$}
		\State $ A_k^{(t)} = \mathcal{M}_k\left(\llbracket \Y;  (\hat{U}_1^{(t)})^\top,  \ldots, (\hat{U}_{k-1}^{(t)})^\top, I_{p_k}, (\hat{U}_{k+1}^{(t-1)})^\top, \ldots, (\hat{U}_{d}^{(t-1)})^\top \rrbracket \right). $
		\If{$k\in J_s$} 				
		\State $\hat U_k^{(t+1)} = \SSVD\left(A_k^{(t)},r_k\right)$
		\Else
		\State $\hat U_k^{(t+1)} = \SVD_{r_k}\left(A_k^{(t)}\right)$
		\EndIf
		\EndFor
		\State $t = t+1$.
		\EndWhile
		
	\end{algorithmic}\label{al:S-HOOI}
\end{algorithm}

\section{Additional Proofs}

\subsection{Proof of Theorem \ref{th:lower_bound_subspace}}
Without loss of generality, we assume $\sigma=1$. Since $\sqrt{2}\left\|\sin \Theta \left(\hat U_k, U_k\right)\right\|_F = \left\|\hat U_k \hat U_k^\top - U_k U_k^\top\right\|_F$, to prove this theorem, it suffices to show that for each $k=1,2,\ldots,d$, we have following inequalities:
\begin{equation}\label{proof_subspace_skrk}
\inf_{\hat{U}_k}\sup_{\X \in \mathcal{F}_{\bp, \bs, \br}^{(k)}} \mathbb{E}\left\|\hat{U}_k\hat{U}_k^\top - U_kU_k^\top\right\|_F^2 \geq c \left(\frac{s_kr_k}{\lambda_k^2} \wedge r_k \right),
\end{equation}
\begin{equation}\label{proof_subspace_sklogp}
\inf_{\hat{U}_k}\sup_{\X \in \mathcal{F}_{\bp, \bs, \br}^{(k)}} \mathbb{E}\left\|\hat{U}_k\hat{U}_k^\top - U_kU_k^\top\right\|_F^2 \geq c \left(\frac{s_k\log{(p_k/s_k)}}{\lambda_k^2} \wedge r_k\right), \quad \text{if} ~k \in J_s.
\end{equation}
We consider the first mode, $k=1$. 
By the proof of Theorem 3 in \cite{zhang2017tensor}, we can construct a core tensor $\S \in \mathbb R^{r_1\times r_2 \times \cdots \times r_d}$ that satisfies the following property:
\begin{equation}\label{good_core_tensor}
\lambda_1 \leq \sigma_{min}(\mathcal M_1(\S)) \leq \sigma_{max}(\mathcal M_1(\S)) \leq C\lambda_1.
\end{equation}
Define the metric space $\mathcal B_{s_1-r_1,r_1} := \{U, U \in \mathbb O_{s_1-r_1,r_1}\}$ equipped with the metric $d(U_1,U_2) := \left\|U_1U_1^\top - U_2U_2^\top\right\|_F$. Recall that the $(\varepsilon\sqrt{r}_1)$-packing number of the metric space $(\mathcal B_{s_1-r_1,r_1},d)$ is defined as:
\begin{equation*}
D(\mathcal B_{s_1-r_1,r_1},d,\varepsilon\sqrt{r_1}) := \max \left\{n:~\text{there are}~t_1,\ldots,t_n \in \mathcal B_{s_1-r_1,r_1},~\text{such that}~\min_{i\neq j}d(t_i,t_j)>\varepsilon\sqrt{r_1}\right\}.
\end{equation*}
By Lemma 5 in \cite{koltchinskii2015optimal}, we have the following control on this packing number if $r_1 \leq s_1-2r_1$:
\begin{equation*}
\left(\frac{c}{\varepsilon}\right)^{r_1(s_1-2r_1)} \leq D(B_{s_1-r_1,r_1},d,\varepsilon\sqrt{r_1}) \leq \left(\frac{C}{\varepsilon}\right)^{r_1(s_1-2r_1)},
\end{equation*}
where $c, C$ are some absolute constants. By choosing $\varepsilon=\frac{c}{2}$, we can find a subset $\mathcal V_{s_1-r_1,r_1} \subset \mathcal B_{s_1-r_1,r_1}$, with Card$\left(\mathcal V_{s_1-r_1,r_1}\right)\geq 2^{r_1(s_1-2r_1)}$ such that for any $V_i \neq V_j \in \mathcal V_{s_1-r_1,r_1}$, 
\begin{equation*}
\left\|V_iV_i^\top - V_jV_j^\top\right\|_F \geq \frac{c}{2}\sqrt{r_1}.
\end{equation*}
Now for each $V_i \in \mathcal V_{s_1-r_1,r_1}$, we define $\tilde{V}_i\in \mathbb O_{p_1,r_1}$ as:
\begin{equation*}
\tilde{V}_i = \begin{bmatrix}
0_{p_1-s_1, r_1}\\
\sqrt{1-\delta}I_{r_1}\\
\sqrt{\delta}V_i
\end{bmatrix}.
\end{equation*}
Now for any $\tilde V_i \neq \tilde V_j \in \mathbb O_{p_1,r_1}$,
\begin{equation*}
\begin{split}
\left\|\tilde V_i \tilde V_i^\top - \tilde V_j\tilde V_j^\top\right\|_F &\geq \sqrt{2\delta(1-\delta)}\left\|V_i - V_j\right\|_F \geq \sqrt{2\delta(1-\delta)}\inf_{O\in \mathbb O_{r_1}} \left\|V_i - V_jO\right\|_F\\
& \geq \sqrt{\delta(1-\delta)}\left\|V_iV_i^\top-V_jV_j^\top\right\|_F\\
&\geq \frac{c}{2}\sqrt{\delta(1-\delta)r_1}.
\end{split}
\end{equation*}
Now we construct a series of fixed signal tensors: $\X_i = \S \times_1 \tilde V_i \times_2 U_2 \times \ldots \times_d U_d, i=1,\ldots, m$, where $U_k \in \mathbb O_{p_k,r_k}(s_k), m=1,\ldots,2^{r_1(s_1-2r_1)}$. By \eqref{good_core_tensor}, we have $\X_1,\ldots, \X_m \in \mathcal{F}_{\bp, \bs, \br}^{(k)}$. Let $\Y_i = \X_i + \Z_i$, where $\Z_i$ are tensors with i.i.d. standard normal distributed entries. Then, $\Y_i \sim N(\X_i, I_{p_1\times p_2 \times \ldots \times p_d})$ and we have:
\begin{equation*}
\begin{split}
D_{KL}(\Y_i\|\Y_j) &= \frac{1}{2}\|\X_i-\X_j\|_F^2 = \frac{1}{2}\left\|\S \times_1 (\tilde V_i - \tilde V_j) \times_2 U_2 \times \ldots \times_d U_d\right\|_F^2 \\
&= \frac{1}{2}\left\|(\tilde V_i - \tilde V_j) \cdot \mathcal M_1(\S) \cdot (U_2 \otimes \ldots \otimes U_d)^\top\right\|_F^2 \\
&\leq C\lambda_1^2 \left\|\tilde V_i - \tilde V_j\right\|_F^2 \leq C\lambda_1^2\delta\left(\|V_i\|_F + \|V_j\|_F\right)^2\\
& = 4C\lambda_1^2\delta r_1.
\end{split}
\end{equation*}
By the generalized Fano's lemma, we have the following lower bound
\begin{equation*}
\inf_{\hat U_1} \sup_{U_1 \in \{\tilde V_i\}_{i=1}^{m}} \mathbb E \|\hat U_1 \hat U_1^\top - U_1U_1^\top\|_F^2 \geq  c\delta(1-\delta) r_1\left(1 - \frac{C\lambda_1^2\delta r_1 + \log 2}{r_1(s_1-2r_1)\log 2}\right).
\end{equation*}
By setting $\delta = c_1\frac{(s_1-2r_1)}{\lambda^2} \wedge \frac{1}{2}$ for a sufficient small constant $c_1 >0$, then under the condition that $s_1>3r_1$, we get \eqref{proof_subspace_skrk}.

To prove \eqref{proof_subspace_sklogp}, we first construct a fixed orthogonal matrix $V^* \in \mathbb R^{(s_1-r_1)\times r_1}$. Denote $t = \lfloor \frac{s_1-r_1}{r_1} \rfloor$, $q = s_1-tr_1$. Since $s_1 \geq 3r_1$, we have $t \geq 1$. Now let
\begin{equation*}
V= \begin{bmatrix}
I_{r_1}\\
\vdots \\
I_{r_1} \\
I_{q} ~;~ 0_{q,r_1-q} 
\end{bmatrix}
\end{equation*}
and $V^* = VA$, where $A \in \mathbb R^{r\times r}$ is a diagonal matrix with first $q$ diagonal entries equal to $1/\sqrt{t+1}$ and the rest diagonal entries equal to $1/\sqrt{t}$. One can verify that $V^* \in \mathbb R^{(s_1-r_1)\times r_1}$ is an orthogonal matrix with the following incoherent constraint: 
\begin{equation*}
\begin{split}
& \frac{r_1}{s_1} \leq \frac{1}{\lfloor s_1/r_1 \rfloor} = \frac{1}{\lfloor (s_1-r_1)/r_1 \rfloor+1} \leq \frac{1}{t+1} \leq \|V^*_{[k,:]}\|_F^2; \\
& \|V^*_{[k,:]}\|_F^2 \leq \frac{1}{t} = \frac{1}{\lceil\frac{s_1-r_1}{r_1} \rceil} \leq \frac{1}{\frac{s_1-r_1}{r_1}-1} \leq \frac{r_1}{s_1-2r_1}, \quad \forall k=1,\ldots, s_1-r_1.
\end{split}
\end{equation*}
By the assumption that $s_1 \geq 3r_1$, we further have $\|V^*_{[k,:]}\|_F^2 \geq \frac{r_1}{s_1} \geq  \frac{r_1}{2(s_1-r_1)}$ for $k=1,\ldots,s_1-r_1$. \\
Let $\Omega_1,\ldots,\Omega_N$ be uniformly random subsets of $\{1,\ldots,p_1-r_1\}$ with ascending order and $|\Omega_1| = \cdots = |\Omega_N| = s_1-r_1$, where $N$ is specified later. Construct $V^{(1)},\ldots,V^{(N)} \in \mathbb R^{(p_1-r_1) \times r_1}$ such that $V^{(i)}_{[\Omega_i,:]} = V^{*}$, $V^{(i)}_{[\Omega_i^c,:]} = 0_{p_1-s_1,r_1}$. Now we define
\begin{equation*}
\tilde V^{(i)}= \begin{bmatrix}
\sqrt{1-\delta}I_{r_1}\\
\sqrt{\delta}V^{(i)}
\end{bmatrix}.
\end{equation*}
It is easy to see that 
$$\|\tilde V^{(i)}-\tilde V^{(j)}\|_F^2 \leq \delta\|V^{(i)} - V^{(j)}\|_F^2 \leq \delta\left(\|V^{(i)}\|_F^2  + \|V^{(j)}\|_F^2\right) \leq 4\delta r_1.$$ 
We further denote $p' = p_1-r_1$ and $s' = s_1-r_1$, note that if $\left|\Omega_i \cap \Omega_j\right| < \frac{s'}{2}$, we must have $\left\|V^{(i)} - V^{(j)}\right\|_F^2 > r_1/2$. Thus we have
\begin{equation*}
\begin{split}
& \mathbb P \left(\left\|\tilde V^{(i)}-\tilde V^{(j)}\right\|_F^2 \leq \frac{r_1\delta}{2}\right) \leq \mathbb P\left(|\Omega_i \cap \Omega_j|\geq\frac{s'}{2}\right) = \sum_{\frac{s'}{2}\leq t \leq s'} \mathbb P(|\Omega_i \cap \Omega_j|=t) \\
& = \sum_{\frac{s'}{2}\leq t \leq s'} \frac{\binom{s'}{t}\binom{p'-s'}{s'-t}}{\binom{p'}{s'}} =  \sum_{\frac{s'}{2}\leq t \leq s'} \binom{s'}{t} \frac{(p'-s')(p'-s'-1)\cdots(p'-2s'+t+1)/(s'-t)!}{p'(p'-1)\cdots(p'-s'+1)/s'!} \\
&  =  \sum_{\frac{s'}{2}\leq t \leq s'} \binom{s'}{t} \frac{(p'-s')(p'-s'-1)\cdots(p'-2s'+t+1)\cdot s' (s'-1) \cdots (s'-t+1)}{p'(p'-1)\cdots(p'-s'+1)} \\
& <  \sum_{\frac{s'}{2}\leq t \leq s'} \binom{s'}{t} \frac{s' (s'-1) \cdots (s'-t+1)}{(p'-s'+t)\cdots(p'-s'+1)}  <  \sum_{\frac{s'}{2}\leq t \leq s'} \binom{s'}{t} \left(\frac{s'}{p'-s'+1}\right)^t \\
& \leq 2^{s'-1}\left(\frac{s'}{p'-s'+1}\right)^{s'/2} = \frac{1}{2}\left(\frac{4s'}{p'-s'+1}\right)^{s'/2}.
\end{split}
\end{equation*}
Let $N = \lfloor 2 (\frac{p'-s'+1}{4s'})^{s'/4}\rfloor$, then 
\begin{equation*}
\begin{split}
	& \mathbb P \left(\left\|\tilde V^{(i)}-\tilde V^{(j)}\right\|_F^2 \geq \delta r_1/2,~\forall i \neq j\right) \geq 1-\frac{N(N-1)}{4}\left(\frac{4s'}{p'-s'+1}\right)^{s'/2} \\
	&  > 1-\frac{N^2}{4}\left(\frac{4s'}{p'-s'+1}\right)^{s'/2} \geq 0.
\end{split}
\end{equation*}
That implies there exist $\{\tilde V^{(i)}\}_{i=1}^{N}$, such that
\begin{equation*}
\delta r_1/2 \leq \left\|\tilde V^{(i)}-\tilde V^{(j)}\right\|_F^2 \leq 4\delta r_1,\quad \forall 1\leq i\neq j \leq N.
\end{equation*}
Then we construct $\X_i = \S \times_1 \tilde V^{(i)} \times_2 U_2 \times \ldots \times_d U_d$ and $\Y_i = \X_i + \Z_i$, where $\Z_i$ has i.i.d. standard normal entries and $i=1,\ldots, N$. Similarly to the proof of \eqref{proof_subspace_skrk}, we have
\begin{equation*}
\left\|\tilde V^{(i)} \tilde V^{(i)\top}-\tilde V^{(j)} \tilde V^{(j)\top}\right\|_F^2 \geq  c(1-\delta)\delta r_1,
\end{equation*}
\begin{equation*}
D_{KL}(\Y_i\|\Y_j) \leq C\lambda_1^2\|\tilde V_i-\tilde V_j\|_F^2 \leq C\lambda_1^2\delta r_1.
\end{equation*}
Note that in the case when $p' < 10s'$, we must have $\log(p_k/s_k) < c$ and \eqref{proof_subspace_sklogp} is directly implied by \eqref{proof_subspace_skrk}. Then we can assume $p' > 10s'$, so that $N\geq 3$. Under such the circumstance, the generalized Fano's lemma gives the following lower bound:
\begin{equation*}
\inf_{\hat U_1} \sup_{U_1 \in \{V_i\}_{i=1}^{m}} \mathbb E \|\hat U_1 \hat U_1^\top - U_1U_1^\top\|_F^2 \geq  c\delta(1-\delta)r_1\left(1 - \frac{C\delta\lambda_1^2r_1 + \log 2}{\frac{s'}{4}\log{((p'-s'+1)/4s'})}\right).
\end{equation*}
We set $\delta = c\frac{s'\log{((p'-s'+1)/4s')}}{\lambda_1^2} \wedge \frac{1}{2}$ with sufficient small constant $c$. Since $s' = s_1-r_1>2r_1$, we can obtain \eqref{proof_subspace_sklogp}. \quad $\square$

\subsection{Proof of Theorem \ref{th:lower_bound_recovery}.} Again, we assume $\sigma=1$. In order to prove this theorem, we only need to show the following inequalities, separately,
\begin{equation}\label{proof_tensor_skrk}
\inf_{\hat{\X}}\sup_{\X \in \mathcal{F}_{\bp, \bs, \br}} \mathbb{E}\left\|\hat{\X} - \X\right\|_F^2 \geq cs_kr_k, \quad \forall k =1,2,\ldots, k,
\end{equation}
\begin{equation}\label{proof_tensor_sklogp}
\inf_{\hat{\X}}\sup_{\X \in \mathcal{F}_{\bp, \bs, \br}} \mathbb{E}\left\|\hat{\X} - \X\right\|_F^2 \geq cs_k\log(p_k/s_k), \quad \forall k \in J_s,
\end{equation}
and
\begin{equation}\label{proof_tensor_r}
\inf_{\hat{\X}}\sup_{\X \in \mathcal{F}_{\bp, \bs, \br}} \mathbb{E}\left\|\hat{\X} - \X\right\|_{F}^2 \geq c r_1\cdots r_d.
\end{equation}
\begin{enumerate}

	\item In order to prove \eqref{proof_tensor_skrk}, it suffices to focus on $k=1$. We use the metric space $\mathcal B_{s_1-r_1,r_1} := \{U \in \mathbb O_{s_1-r_1,r_1}\}$ equipped with the metric $d(U_1,U_2) := \left\|U_1U_1^\top - U_2U_2^\top\right\|_F$. As we showed in the proof of Theorem 2, when $s_1 \geq 3r_1$, we can find a subset $\mathcal V_{s_1-r_1,r_1} \subset \mathcal B_{s_1-r_1,r_1}$, with Card$\left(\mathcal V_{s_1-r_1,r_1}\right)\geq 2^{r_1(s_1-2r_1)}$ such that for each $V_i \neq V_j \in \mathcal V_{s_1-r_1,r_1}$, 
\begin{equation*}
\|V_iV_i^\top - V_jV_j^\top\|_F \geq \frac{c}{2}\sqrt{r_1}.
\end{equation*}
Now for every $V \in \mathcal V_{s_1-r_1,r_1}$, we define $\tilde{V}\in \mathbb O_{p_1,r_1}$ as:
\begin{equation*}
\tilde{V} = \begin{bmatrix}
0_{p_1-s_1, r_1}\\
\sqrt{1-\delta}I_{r_1}\\
\sqrt{\delta}V
\end{bmatrix}.
\end{equation*}
Then For $\tilde V_i \neq \tilde V_j \in \mathbb O_{p_1,r_1}$,
\begin{equation*}
\begin{split}
\|\tilde V_i \tilde V_i^\top - \tilde V_j\tilde V_j^\top\|_F &\geq \sqrt{2\delta(1-\delta)}\|V_i - V_j\|_F \geq \sqrt{2\delta(1-\delta)}\inf_{O\in \mathbb O_{r_1}} \|V_i - V_jO\|_F\\
& \geq \sqrt{\delta(1-\delta)}\|V_iV_i^\top-V_jV_j^\top\|_F\\
&\geq \frac{c}{2}\sqrt{\delta(1-\delta)r_1}.
\end{split}
\end{equation*} Now we construct $\X_i = \S \times_1 \tilde V_i \times_2 U_2 \times \ldots \times_d U_d$ for $i=1,\ldots,2^{r_1(s_1-2r_1)}$. Similarly as the part 1 in the proof of Theorem \ref{th:lower_bound_subspace}, we set $U_2, \ldots U_d$ as fixed sparse orthogonal matrices and $\S\in \mathbb{R}^{r_1\times \cdots \times r_d}$ as a core tensor with the following property: 
\begin{equation*}
\lambda_1 \leq \sigma_{min}(\mathcal M_1(\S)) \leq \sigma_{max}(\mathcal M_1(\S)) \leq C\lambda_1,
\end{equation*}
where $\lambda_1 > C\sqrt{s_1r_1}$. Then we have:
	\begin{equation*}
	\begin{split}
	\left\|\X_i - \X_j\right\|_F^2 &= \left\|(\tilde V_i - \tilde V_j)\cdot \mathcal M_1(\S) \cdot (U_2\otimes \cdots \otimes U_d)\right\|_F^2 \\
	& \geq \lambda_1^2 \left\|\tilde V_i - \tilde V_j\right\|_F^2  \geq \lambda_1^2 \inf_{O \in \mathbb O_{r_1}} \left\|\tilde V_i - \tilde V_jO\right\|_F^2 \\
	&\geq \frac{\lambda_1^2}{2} \left\|\tilde V_i \tilde V_i^\top - \tilde V_j \tilde V_j^\top\right\|_F^2 \geq c\delta (1-\delta)\lambda_1^2 r_1.
	\end{split}
	\end{equation*}
	Now consider $\Y_i = \X_i + \Z_i$, where $\Z_i$ is a random matrix with $N(0,1)$ entries. The KL-divergence between $\Y_i$ and $\Y_j$ could be bounded as:
	\begin{equation*}
	\begin{split}
	D_{KL}(\Y_i\|\Y_j) &= \frac{1}{2}\|\X_i-\X_j\|_F^2 = \|(\tilde V_i - \tilde V_j)\cdot \mathcal M_1(\S) \cdot (U_2\otimes \cdots \otimes U_d)\|_F^2\\
	& \leq C\lambda^2\|\tilde V_i-\tilde V_j\|_F^2 \leq C\delta\lambda_1^2 r_1.
	\end{split}
	\end{equation*}
	By generalized Fano's lemma, we have:
	\begin{equation*}
	\inf_{\hat \X} \sup_{\X \in \{\X_i\}_{i=1}^{m}} \mathbb E\|\hat \X-\X\|_F^2 \geq c\delta(1-\delta)\lambda_1^2 r_1\left(1-\frac{C\delta\lambda_1^2r_1+\log 2}{r_1(s_1-r_1)\log 2}\right).
	\end{equation*}
	We set $\delta = c_1\frac{r_1(s_1-r_1)}{\lambda_1^2}$ with sufficient small constant $c_1$. Since $s_1r_1 = O(\lambda_1^2)$ by our construction, we could get the corresponding lower bound \eqref{proof_tensor_skrk}.
	
	\item The construction of packing sets and the proof of \eqref{proof_tensor_sklogp} is essentially the same as the proof of \eqref{proof_subspace_sklogp}.
	
	\item We finally prove \eqref{proof_tensor_r} by construct a packing set of core tensors. 
	Let $\S_i = \delta \W_i,~i=1,2,\cdots,N$, where $\W_i$ is a random gaussian tensor with i.i.d. standard normal entries and $0<\delta<1$. Note that $\|\S_i-\S_j\|_F^2 = \delta^2\|\W_i-\W_j\|_F^2 \sim 2\delta^2 \chi_{r}^2$, where $r=r_1r_2\cdots r_p$. Then by Lemma \ref{lm:chi-square-tail-bound}, we have:
	\begin{equation*}
	\mathbb P\left(\frac{\|\S_i-\S_j\|_F^2}{2\delta^2} \leq r - 2\sqrt{rx}\right) + \mathbb P\left(\frac{\|\S_i-\S_j\|_F^2}{2\delta^2} \geq r + 2\sqrt{rx} +2x\right) \leq 2e^{-x}.
	\end{equation*}
	Set $x=r/16$, we could obtain:
	\begin{equation*}
	\mathbb P(\|\S_i-\S_j\|_F^2 \leq r\delta^2) + \mathbb P(\|\S_i-\S_j\|_F^2 \geq 4r\delta^2) \leq 2e^{-r/16}.
	\end{equation*}
	Define the event $A_1 = \{r\delta^2\leq\|\S_i-\S_j\|_F^2\leq4r\delta^2,~~\forall i\neq j\}$. Then,
	\begin{equation*}
	\mathbb P (A_1) \geq 1 - N(N-1)e^{-r/16}.
	\end{equation*}
	We could specify $N = \lfloor e^{cr} \rfloor$ such that the above probability lower bound is positive, which implies we could find a packing sets $\{\S_i\}_{i=1}^N$ satisfying
	\begin{equation*}
	r\delta^2 \leq \|\S_i-\S_j\|_F^2 \leq 4r\delta^2,~~\forall i\neq j, \quad N = \lfloor e^{cr}\rfloor.
	\end{equation*}
	Now we let $\X_i = \S_i \times_1 U_1 \times \cdots \times_d U_d$ and $\Y_i = \X_i + \Z_i$, where $U_1,U_2,\ldots,U_d$ are fixed sparse orthogonal matrices and $\Z_i$ has i.i.d. standard normal entries. Then
	\begin{equation*}
	D_{KL}(\Y_i\|\Y_j) = \frac{1}{2}\|\X_i-\X_j\|_F^2 = \frac{1}{2}\|\S_i-\S_j\|_F^2 \leq 2r\delta^2.
	\end{equation*}
	When $cr \geq 2$, $\lfloor e^{cr} \rfloor  \geq e^{cr}- 1 \geq e^{cr-1}$, then by generalized Fano lemma we have:
	\begin{equation*}
		\begin{split}
			& \inf_{\hat \X} \sup_{\X \in \{\X_i\}_{i=1}^{m}} \mathbb E\left\|\hat \X-\X\right\|_F^2 \geq r\delta^2\left(1-\frac{2r\delta^2+\log 2}{\log{\lfloor e^{cr} \rfloor}}\right) \\
			& \geq r\delta^2\left(1-\frac{2r\delta^2+\log 2}{cr-1}\right).
		\end{split}
	\end{equation*}
	Let $\delta^2$ be a sufficient small constant then we obtain \eqref{proof_tensor_r}. 
	
	On the hand, note that by \eqref{proof_tensor_skrk}, we have
	\begin{equation*}
		\inf_{\hat \X} \sup_{\X \in \{\X_i\}_{i=1}^{m}} \mathbb E\left\|\hat \X-\X\right\|_F^2 \geq cs_1r_1 \geq c,
	\end{equation*}	 
	and when $cr < 2$, we can find some universal small constant $c_1$, such that $c_1r < c$, which also implies
	\begin{equation*}
	\inf_{\hat \X} \sup_{\X \in \{\X_i\}_{i=1}^{m}} \mathbb E\left\|\hat \X-\X\right\|_F^2 \geq c_1r.
	\end{equation*}

\end{enumerate}
In summary, we have finished the proof for this theorem.\quad $\square$

\subsection{Proofs to Theorem \ref{th:upper_bound} - Step 3}\label{sec:step-3}

Without loss of generality, we consider the case that $k=1$, while the proof for $k\geq 2$ essentially follows. The proof can be divided into the following two steps: (a) provides upper bounds for $\left\|\sin\Theta\left(\hat{U}_{-1}^{(t)}\hat{V}_1^{(t)} ,U_{-1}V_1\right)\right\|$, and (b) provide upper bounds for $E_k^{(t)}$ and $F_k^{(t)}$.
\begin{enumerate}[leftmargin=*]
	\item We first develop the following perturbation bound for $\hat{U}_{-1}^{(t)}\hat{V}^{(t)}_1$:
	\begin{equation}\label{ineq:step-3-result-1}
	\left\|\sin\Theta\left(\hat{U}_{-1}^{(t)}\hat{V}_1^{(t)}, U_{-1}V_1\right)\right\|^2 \leq \frac{(E_{2}^{(t-1)})^2 + \cdots + (E_{d}^{(t-1)})^2 + 2s_1\eta_1 + 6r_1M_1^2}{\lambda_1^2}.
	\end{equation}
	Recall that the right singular subspace of $X_1$ is $(U_{-1}V_1)$, we turn to consider the upper bound of $\left\|X_1 \left(\hat{U}_{-1}^{(t)} \hat{V}_1^{(t)}\right)_{\perp}\right\|_F^2$. By Lemma \ref{lm:property-orthogonal-matrix}, $\left(\hat{U}_{-1}^{(t)} \hat{V}_1^{(t)}\right)_{\perp} = \left[\hat{U}^{(t)}_{-1\perp} ~~ \hat{U}_{-1}^{(t)}\hat{V}_{1\perp}^{(t)}\right]$, thus
	\begin{equation*}
	\begin{split}
	& \left\|X_1 \left(\hat{U}_{-1}^{(t)} \hat{V}_1^{(t)}\right)_{\perp}\right\|_F^2 = \left\|X_1 \hat{U}_{-1\perp}^{(t)}\right\|_F^2 + \left\|X_1\hat{U}^{(t)}_{-1} \hat{V}_{1\perp}^{(t)} \right\|_F^2\\
	= & \left\|X_1\left(\hat{U}_{2}^{(t-1)}\otimes \cdots \otimes \hat{U}_d^{(t-1)}\right)_{\perp}\right\|_F^2 + \left\|X_1\hat{U}^{(t)}_{-1} \hat{V}_{1\perp}^{(t)} \right\|_F^2\\
	\overset{\text{Lemma \ref{lm:property-orthogonal-matrix}}}\leq & \sum_{l=2}^d \left\|(\hat{U}_{l\perp}^{(t-1)})^\top X_l\right\|_F^2 + \left\|X_1\hat{U}^{(t)}_{-1} \hat{V}_{1\perp}^{(t)} \right\|_F^2\\
	\leq & (E_{2}^{(t-1)})^2 + \cdots + (E_{d}^{(t-1)})^2 + \left\|X_1\hat{U}^{(t)}_{-1} \hat{V}_{1\perp}^{(t)} \right\|_F^2.
	\end{split}
	\end{equation*}
	For $\|X_1\hat{U}_{-1}^{(t)}\hat{V}_{1\perp}^{(t)}\|_F$, since all non-zero rows of $X_1$ are in $I_1$, we have the following decomposition,
	\begin{equation}\label{eq:key-term-decomposition}
	\begin{split}
	& \left\|X_1\hat{U}_{-1}^{(t)}\hat{V}_{1\perp}^{(t)}\right\|_F^2 = \left\|D_{I_1} X_1 \hat{U}_{-1}^{(t)}\hat{V}_{1\perp}^{(t)}\right\|_F^2\\
	= & \left\|D_{\hat{I}_1^{(t)}} X_1 \hat{U}_{-1}^{(t)}\hat{V}_{1\perp}^{(t)}\right\|_F^2 + \left\|D_{I_1\backslash \hat{I}_1^{(t)}} X_1 \hat{U}_{-1}^{(t)}\hat{V}_{1\perp}^{(t)}\right\|_F^2.
	\end{split}
	\end{equation}
	Recall that $\hat{V}_1^{(t)}$ is the leading $r_k$ right singular vectors of $B_1^{(t)} \overset{\eqref{eq:A_k-B_k-bar_A_k-bar_B_k}}{=} B_1^{(X, t)} + B_1^{(Z, t)}$, $\rank(B_1^{(X, t)}) \overset{\eqref{eq:A_k-B_k-bar_A_k-bar_B_k}}{=} \rank(D_{\hat{I}_1^{(t)}}X_1\hat{U}_{-1}^{(t)})\leq r_1$, by Lemma \ref{lm:projection-X-residual}, $\left\|B_1^{(X, t)}\hat{V}_{1\perp}^{(t)}\right\|_F^2 \leq 4r_1\|B_1^{(Z, t)}\|^2$, which yields
	\begin{equation}\label{ineq:step-4-inter1}
	\begin{split}
	& \left\|D_{\hat{I}_1^{(t)}} X_1 \hat{U}_{-1}^{(t)}\hat{V}_{1\perp}^{(t)}\right\|_F^2 \leq 4r_1\left\|D_{\hat{I}_1^{(t)}} Z_1 \hat{U}_{-1}^{(t)}\right\|^2 \overset{\eqref{ineq:A-6}}{\leq} 4r_1\left\|(Z_1)_{[I_1, :]}\hat{U}_{-1}^{(t)}\right\|^2 \leq 4r_1M_1^2.
	\end{split}
	\end{equation}
	For the last inequality, it is because that $\hat{U}_{-1}^{(t)}$ only have nonzero rows on index $\hat{I}^{(t)}_{-1} \subset I_{-1}$.
	Meanwhile,
	\begin{equation}\label{ineq:step-4-inter2}
	\begin{split}
	& \left\|D_{I_1\backslash \hat{I}_1^{(t)}}X_1\hat{U}_1^{(t)}\hat{V}_{1\perp}^{(t)}\right\|_F^2\\
	\overset{\text{Lemma \ref{lm:projection-X-residual}}}{\leq} & \left(\left\|D_{I_1\backslash \hat{I}_1^{(t)}}Y_1\hat{U}_{-1}^{(t)}\hat{V}_{1\perp}^{(t)}\right\|_F + \sqrt{r_1} \left\|D_{I_1\backslash \hat{I}_1^{(t)}}Z_1\hat{U}_{-1}^{(t)}\hat{V}_{1\perp}^{(t)}\right\| \right)^2\\
	\overset{\text{Cauchy-Schwarz}}{\leq} & 2\left\|D_{I_1\backslash \hat{I}_1^{(t)}}Y_1\hat{U}_{-1}^{(t)}\hat{V}_{1\perp}^{(t)}\right\|_F^2 + 2r_1\left\|Z_{1,[I_1,:]}\hat{U}_{-1}^{(t)}\right\|^2\\
	\overset{\eqref{ineq:A-5}}{\leq} & 2\left\|D_{I_1\backslash \hat{I}_1^{(t)}}Y_1\hat{U}_{-1}^{(t)}\right\|_F^2 + 2r_1M_1^2=  2\left\|D_{I_1\backslash \hat{I}_1^{(t)}}A_1^{(t)}\right\|_F^2 + 2r_{1}M_1^2\\
	\leq & 2\sum_{i\in I_1\backslash \hat{I}_1^{(t)}} \|A_{1, [i, :]}^{(t)}\|_2^2 + 2r_1M_1^2\\
	\overset{\eqref{eq:support-I-J}}{\leq} & 2s_1\eta_1 + 2r_1M_1^2.
	\end{split}
	\end{equation}
	Combining \eqref{eq:key-term-decomposition}, \eqref{ineq:step-4-inter1}, and \eqref{ineq:step-4-inter2}, we have $\left\|X_1\hat{U}_{-1}^{(t)}\hat{V}_{1\perp}^{(t)}\right\|_F^2 \leq 2s_1\eta_1 + 6r_1M_1^2$, and
	\begin{equation}
	\left\|X_1\left(\hat{U}_{-1}^{(t)}\hat{V}_1^{(t)}\right)_{\perp}\right\|_F^2 \leq (E_2^{(t-1)})^2 + \cdots + (E_d^{(t-1)})^2 + 2s_1\eta_1 + 6r_1M_1^2.
	\end{equation}

	On the other hand, since the right singular subspace of $X_k$ is $U_{-1}V_1$, we also have the following lower bound:
	\begin{equation*}
	\begin{split}
	& \left\|X_1 \left(\hat{U}_{-1}^{(t)} \hat{V}_1^{(t)}\right)_{\perp}\right\|_F^2 = \left\|X_k P_{U_{-1}V_{1}}\left(\hat{U}_{-1}^{(t)} \hat{V}_1^{(t)}\right)_{\perp}\right\|^2_F\\
	= & \left\|\left(X_1 U_{-1}V_{1}\right) \cdot \left(U_{-1} V_1\right)^\top \left(\hat{U}_{-1}^{(t)} \hat{V}_1^{(t)}\right)_{\perp} \right\|_F^2\\
	\geq & \sigma_{r_1}^2(X_1)  \left\|\left(U_{-1}V_1\right)^\top \left(\hat{U}_{-1}^{(t)} \hat{V}_1^{(t)}\right)_{\perp} \right\|^2 = \lambda_1^2 \left\|\sin\Theta\left(U_{-1}V_1, \hat{U}_{-1}^{(t)} \hat{V}_1^{(t)}\right)\right\|^2.
	\end{split}
	\end{equation*}
	Therefore, we have derived the following upper bound,
	\begin{equation*}
	\begin{split}
	& \left\|\sin\Theta\left(\hat{U}_{-1}^{(t)}\hat{V}_1^{(t)}, U_{-1} V_1\right)\right\|^2\\
	\leq & \frac{(E_{2}^{(t-1)})^2 + \cdots + (E_{d}^{(t-1)})^2 + 2s_1\eta_1 + 6r_1M_1^2}{\lambda_1^2},
	\end{split}
	\end{equation*}
	which has shown \eqref{ineq:step-3-result-1}.
	
	\item Next, we develop the upper bound for $E_1^{(t)}$.
	First, $E_1^{(t)}$ has the following decomposition,
	\begin{equation}\label{ineq:E_1^{(t)}-decomposition}
	\begin{split}
	E_{1}^{(t)} = & \left\|(\hat{U}_{1\perp}^{(t)})^\top X_1 \right\|_F = \left\|(\hat{U}_{1\perp}^{(t)})^\top   X_1\left(P_{\hat{U}_{-1}^{(t)}\hat{V}_1^{(t)}} + P_{(\hat{U}_{-1}^{(t)}\hat{V}_1^{(t)})_\perp}\right) \right\|_F\\
	\leq & \left\|(\hat{U}_{1\perp}^{(t)})^\top X_1 P_{\hat{U}_{-1}^{(t)} \hat{V}_1^{(t)}}\right\|_F + \left\|(\hat{U}_{1\perp}^{(t)})^\top  X_1 P_{\left(\hat{U}_{-1}^{(t)}\hat{V}_1^{(t)}\right)_{\perp}} \right\|_F\\
	& \quad\quad \text{(since the non-zero row index set is $I_1$)}\\
	\leq & \left\|(\hat{U}_{1\perp}^{(t)})^\top  D_{\hat{J}_1^{(t)}} X_1 \hat{U}_{-1}^{(t)}\hat{V}_1^{(t)}\right\|_F + \left\|(\hat{U}_{1\perp}^{(t)})^\top  D_{I_1\backslash \hat{J}_1^{(t)}} X_1 \hat{U}_{-1}^{(t)}\hat{V}_1^{(t)}\right\|_F\\
	& + \left\|(\hat{U}_{1\perp}^{(t)})^\top  X_1 \left(\hat{U}_{-1}^{(t)}\hat{V}_1^{(t)}\right)_{\perp}\right\|_F.
	\end{split}
	\end{equation}
	Next we analyze the three terms in the inequality above respectively.
	\begin{enumerate}
		\item Note that $ D_{\hat{J}_1^{(t)}} Y_1 \hat{U}_{-1}^{(t)}\hat{V}_1^{(t)} = \bar{B}_1^{(t)}$, and $\hat{U}_1^{(t)} = \SVD_{r_1}\left(\bar{B}_1^{(t)}\right) = \SVD_{r_1}\left(\bar{B}_1^{(X, t)} + \bar{B}_1^{(Z, t)}\right)$, then
		\begin{equation*}
		\begin{split}
		& \left\|(\hat{U}_{1\perp}^{(t)})^\top  D_{\hat{J}_1^{(t)}} X_1 \hat{U}_{-1}^{(t)}\hat{V}_1^{(t)}\right\|_F \overset{\text{Lemma \ref{lm:projection-X-residual}}}{\leq} 2\sqrt{r_1}\left\|\bar{B}_1^{(Z, t)}\right\| = 2\sqrt{r_1}\left\| D_{\hat{J}_1^{(t)}} Z_1 \hat{U}_{-1}^{(t)}\hat{V}_1^{(t)}\right\|_F\\
		\leq & 2\sqrt{r_1}\left\| Z_{1, [I_1, :]} \hat{U}_{-1}^{(t)}\hat{V}_1^{(t)}\right\|
		\end{split}
		\end{equation*}
		Meanwhile, by Lemma \ref{lm:key-perturbation-result},
		\begin{equation}\label{ineq:Z_1-hatU-hatV}
		\begin{split}
		& \left\|Z_{1, [I_1, :]}\hat{U}^{(t)}_{-1}\hat{V}_1^{(t)}\right\|\\
		\overset{\eqref{ineq:A-4}\eqref{ineq:A-5}}{\leq} & N_1 + M_1\left( \left\|\sin\Theta\left(\hat{U}_{-1}^{(t)}\hat{V}_1^{(t)}, U_{-1}V_1\right)\right\| + \frac{E_2^{(t-1)}}{\lambda_2} + \cdots + \frac{E_d^{(t-1)}}{\lambda_d}\right).
		\end{split}
		\end{equation}
		Therefore,
		\begin{equation*}
		\begin{split}
		& \left\|(\hat{U}_{1\perp}^{(t)})^\top  D_{\hat{J}_1^{(t)}} X_1 \hat{U}_{-1}^{(t)}\hat{V}_1^{(t)}\right\|_F \\
		\leq & 2\sqrt{r_1}N_1 + 2\sqrt{r_1}M_1\left( \left\|\sin\Theta\left(\hat{U}_{-1}^{(t)}\hat{V}_1^{(t)}, U_{-1}V_1\right)\right\| + \frac{E_2^{(t-1)}}{\lambda_2} + \cdots + \frac{E_d^{(t-1)}}{\lambda_d}\right).
		\end{split}
		\end{equation*}
		\item
		\begin{equation*}
		\begin{split}
		& \left\|(\hat{U}_{1\perp}^{(t)})^\top  D_{I_1\backslash \hat{J}_1^{(t)}} X_1 \hat{U}_{-1}^{(t)}\hat{V}_1^{(t)}\right\|_F \\
		\leq & \left\|(\hat{U}_{1\perp}^{(t)})^\top D_{I_1\backslash\hat{J}_1^{(t)}}Y_{1}\hat{U}_{-1}^{(t)}\hat{V}_1^{(t)}\right\|_F + \left\|(\hat{U}_{1\perp}^{(t)})^\top D_{I_1\backslash \hat{J}_1^{(t)}}Z_1\hat{U}_{-1}^{(t)}\hat{V}_1^{(t)}\right\|_F\\
		\leq & \left\|D_{I_1\backslash\hat{J}_1^{(t)}}\bar{A}_1^{(t)}\right\|_F + \sqrt{r_1} \left\|(\hat{U}_{1\perp}^{(t)})^\top D_{I_1\backslash \hat{J}_1^{(t)}}Z_1\hat{U}_{-1}^{(t)}\hat{V}_1^{(t)}\right\|\\
		\leq & \left(\sum_{i \in I_1\backslash \hat{J}_1^{(t)}}\left\|\bar{A}_{1,[i,:]}^{(t)}\right\|_2^2\right)^{1/2} + \sqrt{r_1}\left\|Z_{1, [I_1, :]}\hat{U}_{-1}^{(t)}\hat{V}_1^{(t)}\right\|\\
		\overset{\eqref{eq:support-I-J}\eqref{ineq:Z_1-hatU-hatV}}{\leq} & \sqrt{s_1\bar{\eta}_1} + \sqrt{r_1}N_1\\
		& + \sqrt{r_1}M_1\left( \left\|\sin\Theta\left(\hat{U}_{-1}^{(t)}\hat{V}_1^{(t)}, U_{-1}V_1\right)\right\| + \frac{E_2^{(t-1)}}{\lambda_2} + \cdots + \frac{E_d^{(t-1)}}{\lambda_d}\right)
		\end{split}
		\end{equation*}
		\item Since the right singular subspace of $X_1$ is $U_{-1}V_1$, one has
		\begin{equation*}
		\begin{split}
		& \left\|\left(\hat{U}_{1\perp}^{(t)}\right)^\top X_1 \left(\hat{U}_{-1}^{(t)}\hat{V}_1^{(t)}\right)_{\perp}\right\|_F \leq \left\|\left(\hat{U}_{1\perp}^{(t)}\right)^\top X_1 \right\|_F \cdot \left\|\sin\Theta\left(\hat{U}_{-1}^{(t)}\hat{V}_1^{(t)}, U_{-1}V_1\right)\right\|\\
		= & E_1^{(t)} \cdot \left\|\sin\Theta\left(\hat{U}_{-1}^{(t)}\hat{V}_1^{(t)}, U_{-1}V_1\right)\right\|.
		\end{split}
		\end{equation*}
	\end{enumerate}
	Combining i, ii, iii, and \eqref{ineq:E_1^{(t)}-decomposition}, we have
	\begin{equation*}
	\begin{split}
	E_1^{(t)} \leq & \sqrt{s_1\bar{\eta}_1} + 3\sqrt{r_1}N_1 + 3\sqrt{r_1}M_1\left( \left\|\sin\Theta\left(\hat{U}_{-1}^{(t)}\hat{V}_1^{(t)}, U_{-1}V_1\right)\right\| + \frac{E_2^{(t-1)}}{\lambda_2} + \cdots + \frac{E_d^{(t-1)}}{\lambda_d}\right)\\
	& + E_1^{(t)}\cdot \left\|\sin\Theta\left(\hat{U}_{-1}^{(t)}\hat{V}_1^{(t)}, U_{-1}V_1\right)\right\|.
	\end{split}
	\end{equation*}
	Thus when $\left\|\sin\Theta (\hat{U}_{-1}^{(t)}\hat{V}_1^{(t)}, U_{-1}V_1)\right\| <1$,
	\begin{equation}\label{ineq:E_1^{(t)}}
	\begin{split}
	E_1^{(t)} \leq & \frac{\sqrt{s_1\bar{\eta}_1} + 3\sqrt{r_1}N_1 + 3\sqrt{r_1}M_1\left( \left\|\sin\Theta\left(\hat{U}_{-1}^{(t)}\hat{V}_1^{(t)}, U_{-1}V_1\right)\right\| + \frac{E_2^{(t-1)}}{\lambda_2} + \cdots + \frac{E_d^{(t-1)}}{\lambda_d}\right)}{1 - \left\|\sin\Theta \left(\hat{U}_{-1}^{(t)}\hat{V}_1^{(t)}, U_{-1}V_1\right)\right\|}\\
	= & \frac{\sqrt{s_1\bar{\eta}_1} + 3\sqrt{r_1}N_1 + 3\sqrt{r_1}M_1\left(K_1^{(t)} + \frac{E_2^{(t-1)}}{\lambda_2} + \cdots + \frac{E_d^{(t-1)}}{\lambda_d}\right)}{1 - K_1^{(t)}},
	\end{split}
	\end{equation}
	where $K^{(t)}_1$ is denoted as
	\begin{equation*}
	\begin{split}
	K^{(t)}_1 := & \left\|\sin\Theta\left(\hat{U}_{-1}^{(t)}\hat{V}_1^{(t)}, U_{-1}V_{1}\right)\right\| \\
	\leq & \left(\frac{(E_{2}^{(t-1)})^2 + \cdots + (E_{d}^{(t-1)})^2 + 2s_1\eta_1 + 6r_1M_1^2}{\lambda_1^2}\right)^{1/2}.
	\end{split}
	\end{equation*}
	Similarly, one can also show for any general $1\leq k \leq d$,
	\begin{equation}\label{ineq:K_t}
	E_k^{(t)} \leq \frac{\sqrt{s_k\bar{\eta}_k} + 3\sqrt{r_k}N_k + 3\sqrt{r_k}M_k\left(K_k^{(t)} + \frac{E_{k+1}^{(t-1)}}{\lambda_{k+1}} + \cdots + \frac{E_{d}^{(t-1)}}{\lambda_{d}} + \frac{E_{1}^{(t)}}{\lambda_{1}} + \cdots + \frac{E_{k-1}^{(t)}}{\lambda_{k-1}}\right)}{1 - K_k^{(t)}},
	\end{equation}
	where
	\begin{equation}\label{ineq:K^t}
	\begin{split}
	K^{(t)}_k := & \left\|\sin\Theta\left(\hat{U}_{-k}^{(t)}\hat{V}_k^{(t)}, U_{-k}V_{k}\right)\right\|\\
	\leq & \left(\frac{(E_{k+1}^{(t-1)})^2 + \cdots + (E_{d}^{(t-1)})^2 + (E_{1}^{(t)})^2 + \cdots + (E_{k-1}^{(t)})^2 + 2s_k\eta_k + 6r_kM_k^2}{\lambda_k^2}\right)^{1/2}.
	\end{split}
	\end{equation}
	
	\item Then we derive the upper bound for $F_k^{(t)}$. Recall that $U_k$ is the left singular vectors of $X_k$, $\rank(X_k) = r_k$, then by Lemma \ref{lm:projection-perturbation-sine-theta},
	\begin{equation*}
	\begin{split}
	\left\|(\hat{U}_{k\perp}^{(t)})^\top X_k \right\|_F \geq \sigma_{r_k}(X_k) \left\|\sin\Theta\left(\hat{U}_k^{(t)}, U_k\right)\right\|_F,
	\end{split}
	\end{equation*}
	which yields
	\begin{equation}\label{ineq:U-U-sin-theta}
	\begin{split}
	F_k^{(t)} = & \left\|\sin\Theta\left(\hat{U}_k^{(t)}, U_k\right)\right\|_F \leq \frac{E_k^{(t)}}{\lambda_k}.
	\end{split}
	\end{equation}
\end{enumerate}

\subsection{Proofs to Theorem \ref{th:upper_bound} - Step 6}\label{sec:step-6}

First, by Lemma \ref{lm:tail-probability}, \eqref{ineq:A-1} holds with probability $1-O(1/p)$.
Note that $\Z \times_1 U_1^\top \times \cdots \times_d U_d^\top $ is a $r$-dimensional projection of i.i.d. Gaussian, and $\Z_{[I_1,\ldots, I_d]}$ is of $s$-dimensional, we have
\begin{equation*}
\begin{split}
& \|\Z_{[I_1,\ldots, I_d]}\|_F^2 \sim \chi_{s}^2, \\
& \left\|\Z \times_1 U_1^\top \times \cdots \times_d U_d^\top \right\|_F^2 \sim \chi^2_{r}.
\end{split}
\end{equation*}
By Lemma \ref{lm:chi-square-tail-bound}, \eqref{ineq:A-2} holds with probability at least $1 - O(1/p)$. Next, since $Z_{k, [I_k, I_{-k}]}$ is a $(s_k)$-by-$(s_{-k})$ matrix with i.i.d. Gaussian entries, by random matrix theory (c.f. Corollary 5.35 in \cite{vershynin2010introduction}), \eqref{ineq:A-3} holds with probability at least $1 - O(d/p)$. Since $U_{-k}\in \mathbb{O}_{p_{-k}, r_{-k}}$ and $V_k\in \mathbb{O}_{r_{-k}, r_k}$ are fixed orthogonal matrix, thus $(Z_k)_{[I_1, :]}U_{-k}V_{k}$ is a $s_k$-by-$r_k$ i.i.d. Gaussian matrix, then by Corollary 5.35 in \cite{vershynin2010introduction}, \eqref{ineq:A-4} holds with probability at least $1 - O(d/p)$. Then, by Lemma \ref{lm:random-projection-lemma}, \eqref{ineq:A-5} holds with probability at least $1 - O(1/p)$.

Here we want to emphasize that Conditions \eqref{ineq:A-1} -- \eqref{ineq:A-5} are actually only rely on $Z_{I_1,\ldots, I_k}$, i.e. the noise on the non-zero entries of $\X$.

The evaluation for the probability that \eqref{ineq:A-6} is fairly complicated. To this end, we construct a parallel sequence of estimations rather than working directly on $\hat{I}_k^{(t)}, \hat{J}_k^{(t)}, \hat{U}_k^{(t)}, \hat{V}_k^{(t)}$ as follows.
\begin{enumerate}
	\item Initialize
	\begin{equation}\label{eq:I_k^{(0)}-cal}
	\begin{split}
	\hat{\cI}_k^{(0)} = & I_k \cap \Bigg(I_k^{(0)} \cup \Bigg\{i \in I_k: \left\|(Y_k)_{[i, :]}\right\|_2^2 \geq \sigma^2 \left(p_{-k} + 2\sqrt{p_{-k}\log p} + 2 \log p\right)\\
	& \quad\quad\quad \text{or } \max_{j} |(Y_k)_{[i,j]}| \geq 2\sigma\sqrt{\log p}\Bigg\}\Bigg),\quad k\in J_s;
	\end{split}
	\end{equation}
	\begin{equation*}
	\hat{\mathcal{I}}^{(0)}_k = \{1,\ldots, p_k\}, \quad k \notin J_s;
	\end{equation*}
	$$\hat{\cU}^{(0)}_k = \SVD_{r_k}\left(\mathcal{M}_k(\tilde{\Y})\right) \in \mathbb{R}^{p_k\times r_k},\quad \tilde{\Y} =\left\{\begin{array}{ll}
	\Y_{i_1,\ldots,i_d}, & (i_1,\ldots, i_d)\in \hat{\cI}_1^{(0)}\otimes \cdots \otimes \hat{\cI}_d^{(0)};\\
	0, & \text{otherwise}.
	\end{array}\right.$$
	\item For $t = 1,\ldots, t_{\max}$, $k=1,\ldots, d$, let
	$$\cA_k^{(t)} = \mathcal{M}_k\left(\Y \times_1 (\hat{\cU}_1^{(t)})^\top \times \cdots \times_{k-1} (\hat{\cU}_{k-1}^{(t)})^\top \times_{k+1} (\hat{\cU}_{k+1}^{(t-1)})^\top \times \cdots \times_d (\hat{\cU}_{d}^{(t-1)})^\top \right)\in \mathbb{R}^{p_k\times r_{-k}};$$
	$$\cB_k^{(t)} \in \mathbb{R}^{p_k\times r_{-k}}, \quad \cB^{(t)}_{k, [i, :]} = \cA^{(t)}_{k, [i, :]} 1_{\left\{\|\cA^{(t)}_{k, [i, :]}\|_2^2 \geq \eta_k \text{ and } i\in I_k\right\}}, \quad 1\leq i \leq p_k,$$
	where $\eta_k = \sigma^2\left(r_{-k}+2\left(\left(r_{-k} \log p\right)^{1/2} + \log p\right)\right)$.
	$$\hat{\cV}_k^{(t)} = \SVD_{r_k}\left(\cB_k^{(t)\top}\right) \in \mathbb{O}_{p_k, r_k}.$$
	For each $t=1,\ldots, t_{\max}, k=1,\ldots, d$, after obtaining $\mathcal{V}_k^{(t)}$, we calculate the projection $\bar{\cA}_k^{(t)} = \cA_k^{(t)} \hat{\cV}_k^{(t)} \in\mathbb{R}^{p_k\times r_k}$, and
	$$\bar{\cB}_k^{(t)} \in \mathbb{R}^{p_k\times r_{k}}, \quad \bar{\cB}^{(t)}_{k, [i, :]} = \bar{\cA}^{(t)}_{k, [i, :]} 1_{\left\{\|\cA^{(t)}_{k, [i, :]}\|_2^2 \geq \bar{\eta}_k  \text{ and } i\in I_k\right\}}, \quad 1\leq i \leq p_k,$$
	where $\bar{\eta}_k = \sigma^2\left(r_k + 2\left((r_k\log p)^{1/2} + \log p\right)\right)$. Finally, apply QR decomposition to $\bar{\cB}_k^{(t)}$, and assign the $Q$ part to $\hat{\cU}_k^{(t)}\in \mathbb{O}_{p_k, r_k}$.
\end{enumerate}
Intuitively speaking, the above procedure is the counterpart of Algorithm \ref{al:procedure} restricted on index sets $I_1\times \cdots\times I_d$: particularly $\hat{\cI}_k^{(t)}, \hat{\cJ}_k^{(t)}\subseteq I_k$ always hold; meanwhile, by taking the union of $I^{(0)}$, it makes sure that $I^{(0)} \subseteq \hat{\cI}_k^{(0)}$. In order to show \eqref{ineq:A-6}, we aim to prove that the sequence of $\hat{\cI}_k^{(t)}, \hat{\cJ}_k^{(t)}, \hat{\cU}^{(t)}_k, \hat{\cV}_k^{(t)}$ coincide with the original sequence up to $t_{\max}$ steps with high probability, i.e.
\begin{equation}\label{eq:final-step-coincide-0}
I^{(0)} \subseteq \hat{I}^{(0)} = \hat{\cI}^{(0)}
\end{equation}
and
\begin{equation}\label{eq:final-step-coincide}
\hat{\cI}_k^{(t)} = I_k^{(t)}, \hat{\cJ}_k^{(t)} = J_k^{(t)}, \hat{\cU}^{(t)}_k = \hat{U}^{(t)}_k, \hat{\cV}_k^{(t)} = \hat{V}_k^{(t)}, \quad t = 1,\ldots, t_{\max};\quad k=1,\ldots, d,
\end{equation}
with probability at least $1 - O(\log p/p)$.

Next, we particularly prove \eqref{eq:final-step-coincide-0} and \eqref{eq:final-step-coincide} by conditioning on fixed $\Z_{[I_1,\ldots, I_d]}$ satisfying Conditions \eqref{ineq:A-1} -- \eqref{ineq:A-5}. By conditioning on fixed $\Z_{[I_1,\ldots, I_d]}$, the system achieves the following two important properties:
\begin{enumerate}
	\item the entries of $\Z$ outside of $[I_1, \ldots, I_d]$ are still i.i.d. Gaussian distributed, given the independence between $Z_{[I_1,\ldots, I_d]}$ and $\Z_{[I_1,\ldots, I_d]^c}$;
	\item $\hat{\cI}^{(t)}, \hat{\cJ}^{(t)}, \hat{\cU}^{(t)}, \hat{\cV}^{(t)}$ all becomes fixed, since they all only relies on $\X$ and $\Z_{[I_1,\ldots, I_d]}$.
\end{enumerate}
By comparing the procedures for $\hat{\cI}, \hat{\cJ}, \hat{\cU}, \hat{\cV}$ above and Algorithm \ref{al:procedure}, \eqref{eq:final-step-coincide-0} and \eqref{eq:final-step-coincide} are implied by the following statements: \eqref{ineq:final-coincide-to-prove-1} -- \eqref{ineq:final-coincide-to-prove-4}.
\begin{equation}\label{ineq:final-coincide-to-prove-1}
\begin{split}
\forall 1\leq k \leq d, \forall i\notin I_k, \quad & \left\|(Y_k)_{[i, :]}\right\|_2^2 < \sigma^2 \left(p_{-k} + 2\sqrt{p_{-k}\log p} + 2 \log p\right);\\
& \quad \text{and  } \max_{j} |(Y_k)_{[i,j]}| < 2\sigma\sqrt{\log p};
\end{split}
\end{equation}
\begin{equation}\label{ineq:final-coincide-to-prove-2}
\begin{split}
\forall 1\leq k \leq d, \forall i\in I_k^{(0)}, \quad & \left\|(Y_k)_{[i, :]}\right\|_2^2 \geq \sigma^2 \left(p_{-k} + 2\sqrt{p_{-k}\log p} + 2 \log p\right);\\
& \quad \text{or  } \max_{j} |(Y_k)_{[i,j]}| \geq 2\sigma\sqrt{\log p};
\end{split}
\end{equation}
\begin{equation}\label{ineq:final-coincide-to-prove-3}
\begin{split}
\forall 1\leq k \leq d, 1\leq t \leq t_{\max}, i\notin I_k, \quad & \left\|\cA_{k, [i, :]}^{(t)}\right\|_2^2 = \left\|\left(Y_k \hat{\cU}_{-k}^{(t)}\right)_{[i,:]}\right\|_2^2 < r_{-k} + \eta_k;\\
\end{split}
\end{equation}
\begin{equation}\label{ineq:final-coincide-to-prove-4}
\begin{split}
\forall 1\leq k \leq d, 1\leq t \leq t_{\max}, i\notin I_k, \quad & \left\|\bar{\cA}_{k, [i, :]}^{(t)}\right\|_2^2 = \left\|\left(Y_k \hat{\cU}_{-k}^{(t)}\hat{\cV}_k^{(t)}\right)_{[i,:]}\right\|_2^2 < r_k + \bar{\eta}_k.
\end{split}
\end{equation}
\begin{enumerate}
	\item Note that for any $k=1,\ldots, d$ and $i\notin I_k$,
	$$\left\|(Y_k)_{[i, :]}\right\|_2^2 = \|(Z_k)_{[i,:]}\|_2^2 \sim \chi^2_{p_{-k}},$$
	$$\max_{j} |(Y_k)_{[i, j]}| = \max_{j} |(Z_k)_{[i, j]}| \text{ is the maximum of $p_{-k}$ i.i.d. Gaussians}.$$
	Thus by Lemma \ref{lm:tail-probability}, \eqref{ineq:final-coincide-to-prove-1} holds with probability at least $1-O\left((p_1+\cdots+p_d)/p\right)$.
	\item Since $X_{k}$ is supported on $I_k \times I_{-k}$, then for any $1\leq k \leq d$, $i \in I_k^{(0)}$,
	\begin{equation*}
	\begin{split}
	\|X_{k, [i, :]}\|_2^2 \geq 16 s_{-k}\log p \quad \Rightarrow & \quad \max_{j\in I_{-k}} \left(X_{k, [i, j]}\right)^2 \geq \frac{16 s_{-k}\log p}{s_{-k}} = 16\log p\\
	\Rightarrow & \quad \max_{j\in I_{-k}} \left|X_{k, [i, j]}\right| \geq 4 \sqrt{\log p}\\
	\overset{\eqref{ineq:A-1}}{\Rightarrow} & \quad \max_{j\in I_{-k}} \left|Y_{k, [i,j]}\right| \geq 4\sqrt{\log p} - 2\sqrt{\log p} = 2\sqrt{\log p}\\
	\Rightarrow & \quad \max_{j} \left|Y_{k, [i,j]}\right| \geq 2\sqrt{\log p},
	\end{split}
	\end{equation*}
	which means \eqref{ineq:final-coincide-to-prove-2} definitely holds given \eqref{ineq:A-1} -- \eqref{ineq:A-5}.
	
	\item Note for any $1\leq k \leq d, 1\leq t \leq t_{\max}, i\notin I_k$, $X_{k, [i, :]} = 0$, $Z_{k, [i, :]}\in \mathbb{R}^{p_{-k}}$ is an i.i.d. Gaussian vector, $\hat{\cU}_{-k}^{(t)}$ is a fixed $p_{-k}$-by-$r_{-k}$ orthogonal matrix, then
	\begin{equation*}
	\left\|\mathcal{A}_{k, [i, :]}^{(t)}\right\|_2^2 = \left\|\left(Y_k\hat{\cU}_{-k}^{(t)}\right)_{[i, :]}\right\|_2^2 = \left\|\left(Z_{k, [i, :]}\hat{\cU}_{-k}^{(t)}\right)_{[i, :]}\right\|_2^2 \sim \chi_{r_{-k}}^2,
	\end{equation*}
	thus by Lemma \ref{lm:chi-square-tail-bound}, \eqref{ineq:final-coincide-to-prove-3} holds with probability $1 - O\left(t_{\max}(p_1+\cdots+p_d)/p\right)$.
	
	\item Similarly as the previous case, for any $1\leq k \leq d, 1\leq t \leq t_{\max}, i\notin I_k$, $X_{k, [i, :]} = 0$, $Z_{k, [i, :]}\in \mathbb{R}^{p_{-k}}$ is an i.i.d. Gaussian vector, $\hat{\cU}_{-k}^{(t)}\hat{\cV}_{k}^{(t)}$ is a fixed $p_{-k}$-by-$r_{k}$ orthogonal matrix, then
	\begin{equation*}
	\left\|\bar{\cA}_{k, [i, :]}^{(t)}\right\|_2^2 = \left\|\left(Y_k\hat{\cU}_{-k}^{(t)}\hat{\cV}_{k}^{(t)}\right)_{[i, :]}\right\|_2^2 = \left\|\left(Z_{k, [i, :]}\hat{\cU}_{-k}^{(t)}\hat{\cV}_{k}^{(t)}\right)_{[i, :]}\right\|_2^2 \sim \chi_{r_{k}}^2,
	\end{equation*}
	thus by Lemma \ref{lm:chi-square-tail-bound}, \eqref{ineq:final-coincide-to-prove-4} holds with probability $1 - O\left(t_{\max}(p_1+\cdots+p_d)/p\right)$.
\end{enumerate}
Therefore, conditioning fixed $\Z_{[I_1,\ldots, I_d]}$ satisfying \eqref{ineq:A-1} -- \eqref{ineq:A-5}, \eqref{eq:final-step-coincide-0}, \eqref{eq:final-step-coincide}, and \eqref{ineq:final-coincide-to-prove-1} -- \eqref{ineq:final-coincide-to-prove-4} hold with probability at least $1 - O\left(t_{\max} (p_1+\cdots+p_d)/ p\right)$. Finally, we conclude from the analysis in this step that,
\begin{equation*}
\begin{split}
& P\left(\text{\eqref{ineq:A-1} -- \eqref{ineq:A-6} all hold}\right)\\
= & P\left(\text{\eqref{ineq:A-1} -- \eqref{ineq:A-5} all hold}\right) \cdot P\left(\text{\eqref{ineq:A-6} holds}\big|\text{\eqref{ineq:A-1} -- \eqref{ineq:A-5} all hold}\right)\\
\geq & \left(1 - Cd/p\right) \left(1 - Ct_{\max}(p_1+\cdots+p_d)/p \right)\\
\geq & 1 - O\left(\frac{(p_1+\cdots+p_d)\log(ds\log(p))}{p}\right).
\end{split}
\end{equation*}

\subsection{Proof of Proposition \ref{prop:sigma_estimate}.}\label{tuning-parameter}

Let $Y, X, Z \in \mathbb{R}^p$ be vectorizations of $\Y, \X, \Z$. Since the sample median does not depend on the order of its entries, we can assume without loss of generality that the first $s$ elements of $X$ correspond to the non-sparse entries of $\X$. Then, $|Y_i| \sim |N(X_i,\sigma^2)|$ for any $1\leq i \leq s$ and $|Y_i| \sim |N(0,\sigma^2)|$ for $s+1 \leq i \leq p$. Denote the median of $|N(0,\sigma^2)|$ as $m = \sigma z_{0.75}$ and $b_i = 1_{\{Y_i > m + \varepsilon \sigma z_{0.75}\}}$. Then $b_i \sim {\rm Bernoulli}(q)$ for $i=s+1,\ldots,p$, where $q=2-2\Phi(z_{0.75}+\varepsilon)$ and $\Phi(x)$ is the distribution function of $N(0,1)$. Thus we have:
\begin{equation*}
	\begin{split}
		\mathbb{P} \left(\hat{\sigma} \geq \sigma + \sigma \varepsilon \right) & = \mathbb{P} \left({\rm median}(|Y|) \geq m + \sigma\varepsilon z_{0.75}\right) \\
		& \leq \mathbb{P} \left(\sum_{i=1}^p b_i \geq p/2  \right) \leq \mathbb{P} \left(\sum_{i=s+1}^p b_i \geq  p/2 - s \right) \\
		& = \mathbb{P}\left( \sum_{i=s+1}^p b_i - (p-s)q \geq p(0.5 - q) - (1-q)s\right)\\
		& \leq  e^{-2\left(p(0.5-q) - (1-q)s\right)^2/(p-s)},
	\end{split}
\end{equation*}
where the last inequality comes from Hoeffding's inequality. Now we set $\varepsilon = c\sqrt{(\log p)/p}$. For some small constant $c>0$, we additionally have $q = 2 - 2\Phi(z_{0.75}+\varepsilon) = 2 - 2\Phi(z_{0.75}+c\sqrt{(\log p)/p})\geq 1/2 - \sqrt{(\log p)/p}$. Since $s = o(\sqrt{p\log p})$, we have $2\left(p(0.5-q) - (1-q)s\right)^2/(p-s) \gtrsim \log p$. Thus,
$$\mathbb{P}\left(\hat \sigma > \sigma + c\sigma\sqrt{(\log p)/p}\right) \leq e^{\log p} = 1/p.$$ 
We can similarly prove that $\mathbb{P} \left(\hat{\sigma} \leq \sigma - c\sigma \sqrt{(\log p)/p}\right) \leq 1/p$ for the same $c>0$, which has finished the proof of this proposition. \quad $\square$

\subsection{Proof of Proposition \ref{prop:rank_estimate}.}
Without loss of generality, we assume $\sigma = 1$. For each $r_k$, we aim to show that $\hat r_k \leq r_k$ and $\hat r_k \geq r_k$ both happen with high probability. Note that the estimation procedure relies on $\hat \sigma$ instead of $\sigma$ now, we introduce similar notations with some constants changed as the one introduced in the proof of Theorem \ref{th:upper_bound}. Define the index sets similarly as \eqref{eq:step-1-index} and \eqref{eq:support-I-J},
\begin{equation}\label{eq:new-step-1-index}
			\begin{split}
				\text{(True support)}\quad & I_k = \left\{i: (X_k)_{[i, :]} \neq 0\right\}, \\
				\text{(Significant support)}\quad & I_k^{(0)} = \left\{i: \left\|(X_k)_{[i, :]}\right\|_2^2 \geq 
				25\hat \sigma s_{-k}\log p\right\},\\
				\text{(Selected support)} \quad & \hat{I}_k^{(0)} = \Bigg\{i: \left\|(Y_k)_{[i, :]}\right\|_2^2 \geq \hat\sigma^2 \left(p_{-k} + 2\sqrt{p_{-k}\log p} + 2 \log p\right)\\
				& \quad \text{or } \max_{j} |(Y_k)_{[i,j]}| \geq 2\hat\sigma\sqrt{\log p}\Bigg\},\quad k=1,\ldots, d.
			\end{split}
		\end{equation}
		For $t\geq 1$, we also define $\hat{I}_k^{(t)}$ and $\hat{J}_k^{(t)}$ with $k \in J_{s}$:
		\begin{equation}\label{eq:new-support-I-J}
			\begin{split}
				(\text{selected support for $A_k$ in $t$-th step})\quad & \hat{I}_k^{(t)} = \left\{i: \left\|(A_k^{(t)})_{[i, :]}\right\|_2^2 \geq \eta_k\right\},\\	
				(\text{selected support for $\bar{A}_k$ in $t$-th step})\quad & \hat{J}_k^{(t)} = \left\{i: \left\|(\bar{A}_k^{(t)})_{[i, :]}\right\|_2^2 \geq \bar{\eta}_k\right\},\\	
			\end{split}
		\end{equation}
where $\eta_k = \hat \sigma^2\left(\hat r_{-k} + 2\sqrt{\hat r_{-k}\log p} + 2\log p\right)$, $\bar{\eta}_k = \hat \sigma^2 \left(\hat r_k + 2\sqrt{\hat r_k\log p} + 2\log p\right)$.

We also list the following assumptions:
\begin{enumerate}[label=(\subscript{A}{\arabic*})]
		\item
		\begin{equation}\label{new-ineq:A-1}
			\max_{(i_1,\ldots, i_d)\in I_1\times \cdots \times I_d} \left|Z_{i_1,\ldots, i_d}\right| \leq 2\sqrt{\log p}.
		\end{equation}
		\item
		\begin{equation}\label{new-ineq:A-2}
			\begin{split}
				& \left\|\Z_{[I_1,\ldots, I_d]}\right\|_F^2 \leq s + 2\sqrt{s \log p} + 2 \log p,\\
				& \left\|\Z \times_1 U_1^\top \times \cdots \times_d U_d^\top \right\|_F^2 \leq r + 2\sqrt{r\log p} + 2\log p.
			\end{split}
		\end{equation}
		\item
		\begin{equation}\label{new-ineq:A-3}
			\forall 1\leq k \leq d, \quad \left\|Z_{k, [I_k, I_{-k}]}\right\| \leq \sqrt{s_k} + \sqrt{s_{-k}} + 2\sqrt{\log p}.
		\end{equation}
		\item $\forall 1\leq k \leq d$,
		\begin{equation}\label{new-ineq:A-4}
			N_k := \left\|(Z_k)_{[I_k, :]}U_{-k}V_k\right\| \leq \sqrt{s_k}+\sqrt{r_k}+2\sqrt{\log p}.
		\end{equation}
		\item $\forall 1\leq k \leq d$,
		\begin{equation}\label{new-ineq:A-5}
			M_k := \sup_{W_l \in \mathbb{R}^{s_l\times r_l}, \|W_l\|\leq 1} \left\|Z_{k, [I_k, I_{-k}]} W_{-k}\right\| \leq 2\left(\sqrt{s_k} + \sqrt{r_{-k}} + \sum_{l\neq k} (\sqrt{s_lr_l} + \sqrt{\log p})\right),
		\end{equation}
		where $W_{-k} = W_{k+1}\otimes \cdots \otimes W_d \otimes W_1 \otimes \cdots \otimes W_{k-1} \in \mathbb{O}_{s_{-k}, r_{-k}}$.
		\item Support consistency for initialization step:
		\begin{equation}\label{new-ineq:A-6}
			I_k^{(0)} \subseteq \hat I_k^{(0)} \subseteq I_k.
		\end{equation}
		\item Estimation error of $\hat \sigma$:
		\begin{equation}\label{new-ineq:A-7}
			\begin{split}
				1 - c\sqrt{\frac{\log p}{p}} \leq \hat \sigma \leq 1 + c\sqrt{\frac{\log p}{p}}.
			\end{split}
		\end{equation}
	\end{enumerate}

We first show $\hat r_k \leq r_k$ and $\hat r_k \geq r_k$ both happen with high probability under Assumptions $(A_6)$ and $(A_7)$. The proof of the first part follows the idea of the proof of Proposition 7 in \cite{yang2016rate}. Specifically, if $p \geq C$ for large constant $C>0$, we have
\begin{equation*}
	\begin{split}
		& \Proj (\hat r_k > r_k) = \Proj \left(\sigma_{r_k+1}(Y_{k,[\hat{I}_k^0,\hat{J}_k^0]}) > \hat\sigma \delta_{|\hat{I}_k^0||\hat{J}_k^0|}\right) \leq \Proj \left(\max_{|A|=|\hat{I}_k^0|,|B|=|\hat{J}_k^0|} \sigma_{r_k+1}(Y_{k,[A,B]}) > \hat\sigma \delta_{|A||B|}\right) \\
		& \leq  \sum_{i=r_k+1}^{p_k} \sum_{j=r_k+1}^{p_{-k}}\Proj\left(\max_{|A|=i,|B|=j} \sigma_{r+1}(Y_{k,[A,B]}) > \hat\sigma \delta_{ij}\right) \\
		& \leq  \sum_{i=r_k+1}^{p_k} \sum_{j=r_k+1}^{p_{-k}}\Proj\left(\max_{|A|=i,|B|=j} \sigma_1(Z_{k,[A,B]}) > \hat\sigma \delta_{ij}\right)  \\
		& \leq  \sum_{i=r_k+1}^{p_k} \sum_{j=r_k+1}^{p_{-k}}\Proj\left(\max_{|A|=i,|B|=j} \sigma_1(Z_{k,[A,B]}) > \hat\sigma \delta_{ij}\right) \\
		& \leq \sum_{i=r_k+1}^{p_k} \sum_{j=r_k+1}^{p_{-k}} {p_k \choose i}  {p_{-k} \choose j} \Proj \left(\sigma_1 (Z_{k,[A, B]})>\left(1-c\sqrt{\frac{\log p}{p}}\right)\delta_{ij}\right) \\
		& \leq \sum_{i=r_k+1}^{p_k} \sum_{j=r_k+1}^{p_{-k}} \left(\frac{ep_k}{i}\right)^i  \left(\frac{ep_{-k}}{j}\right)^j \Proj \left(\sigma_1 (Z_{k,[A, B]})>0.99\delta_{ij}\right)  \\
		& \leq \sum_{i=r_k+1}^{p_k} \sum_{j=r_k+1}^{p_{-k}} \left(\frac{ep_k}{i}\right)^i  \left(\frac{ep_{-k}}{j}\right)^j \exp\left(-i \log \frac{ep_k}{i} - j \log \frac{ep_{-k}}{j} - 2\log p\right) \\
		& \leq \sum_{i=r_k+1}^{p_k} \sum_{j=r_k+1}^{p_{-k}} p^{-2} \leq p^{-1} .
	\end{split}
\end{equation*}
The second part relies on the intermediate result in the proof of Theorem \ref{th:upper_bound}. When $(A_6)$ and $(A_7)$ hold, $|\hat I_{k}^{(0)}| \leq |I_{k}| = s_k$, $|\hat I_{-k}^{(0)}| \leq |I_{-k}| = s_{-k}$, and
\begin{equation*}
	\begin{split}
		\sigma_{r_k} (\tilde{Y}_k) &= \sigma_{r_k} (D_{\hat{I}_k^{(0)}}Y_k D_{\hat{I}_{-k}^{(0)}}) \\
		& \geq \sigma_{r_k} (D_{\hat{I}_k^{(0)}}X_k D_{ \hat{I}_{-k}^{(0)}}) - \|D_{\hat{I}_k^{(0)}}Z_k D_{\hat{I}_{-k}^{(0)}}\| \\
		& \geq \sigma_{r_k} (X_k) - \|X_k - D_{\hat{I}_k^{(0)}}X_k D_{\hat{I}_{-k}^{(0)}}\| - \|D_{\hat{I}_k^{(0)}}Z_k D_{\hat{I}_{-k}^{(0)}}\| \\
		& \geq 0.9 \lambda_k \gtrsim \left(1+\sqrt{\frac{\log p}{p}}\right) \delta_{s_ks_{-k}} \geq \hat \sigma\delta_{|\hat{I}_k^{(0)}||\hat I_{-k}^{(0)}|}.
	\end{split}
\end{equation*}
The last but two inequality is the same as Step 2(a) in the proof of Theorem \ref{th:upper_bound}; the last but one inequality comes from the signal-noise condition and the fact that $\delta_{s_ks_{-k}} = O(\sqrt{s\log p})$. Since the above inequality implies that $\hat r_k \geq r_k$, we know that under assumptions $(A_1)-(A_7)$, $\hat r_k = r_k$ with high probability. \\
Now it suffices to show Assumptions $(A_1)-(A_7)$ hold with high probability. Note that assumptions $(A_1)-(A_5)$ have nothing to do with hyper-parameters, thus step 6 in Theorem \ref{th:upper_bound} has already shown $(A_1)-(A_5)$ happen with probability at least $1 - O((p_1+\cdots+p_d)/p)$. Also, Proposition \ref{prop:sigma_estimate} showed $(A_7)$ happens with probability at least $1 - O(1/p)$. So we only need to show $(A_6)$ happens with high probability. To this end, we redefine $\hat{\cI}_k^{(t)}, \hat{\cJ}_k^{(t)}, \hat{\cU}^{(t)}_k, \hat{\cV}_k^{(t)}$ as the outcome of the oracle version of Algorithm \ref{al:procedure} with $\sigma$ replaced by $\hat{\sigma}$, $r_k$ replaced by $\hat r_k$. 
\begin{enumerate}
	\item Initialize
	\begin{equation*}
	\begin{split}
	\hat{\cI}_k^{(0)} = & I_k \cap \Bigg(I_k^{(0)} \cup \Bigg\{i \in I_k: \left\|(Y_k)_{[i, :]}\right\|_2^2 \geq \hat{\sigma}^2 \left(p_{-k} + 2\sqrt{p_{-k}\log p} + 2 \log p\right)\\
	& \quad\quad\quad \text{or } \max_{j} |(Y_k)_{[i,j]}| \geq 2\hat{\sigma}\sqrt{\log p}\Bigg\}\Bigg),\quad k\in J_s;
	\end{split}
	\end{equation*}
	\begin{equation*}
	\hat{\mathcal{I}}^{(0)}_k = \{1,\ldots, p_k\}, \quad k \notin J_s;
	\end{equation*}
	$$\hat{\cU}^{(0)}_k = \SVD_{r_k}\left(\mathcal{M}_k(\tilde{\Y})\right) \in \mathbb{R}^{p_k\times r_k},\quad \tilde{\Y} =\left\{\begin{array}{ll}
	\Y_{i_1,\ldots,i_d}, & (i_1,\ldots, i_d)\in \hat{\cI}_1^{(0)}\otimes \cdots \otimes \hat{\cI}_d^{(0)};\\
	0, & \text{otherwise}.
	\end{array}\right.$$
	\item For $t = 1,\ldots, t_{\max}$, $k=1,\ldots, d$, let
	$$\cA_k^{(t)} = \mathcal{M}_k\left(\Y \times_1 (\hat{\cU}_1^{(t)})^\top \times \cdots \times_{k-1} (\hat{\cU}_{k-1}^{(t)})^\top \times_{k+1} (\hat{\cU}_{k+1}^{(t-1)})^\top \times \cdots \times_d (\hat{\cU}_{d}^{(t-1)})^\top \right)\in \mathbb{R}^{p_k\times \hat r_{-k}};$$
	$$\cB_k^{(t)} \in \mathbb{R}^{p_k\times \hat r_{-k}}, \quad \cB^{(t)}_{k, [i, :]} = \cA^{(t)}_{k, [i, :]} 1_{\left\{\|\cA^{(t)}_{k, [i, :]}\|_2^2 \geq \eta_k \text{ and } i\in I_k\right\}}, \quad 1\leq i \leq p_k,$$
	where $\eta_k = \hat{\sigma}^2\left(\hat r_{-k}+2\left(\hat r_{-k} \log p\right)^{1/2} + 2\log p\right)$.
	$$\hat{\cV}_k^{(t)} = \SVD_{\hat r_k}\left(\cB_k^{(t)\top}\right) \in \mathbb{O}_{p_k, \hat r_k}.$$
	For each $t=1,\ldots, t_{\max}, k=1,\ldots, d$, after obtaining $\mathcal{V}_k^{(t)}$, we calculate the projection $\bar{\cA}_k^{(t)} = \cA_k^{(t)} \hat{\cV}_k^{(t)} \in\mathbb{R}^{p_k\times \hat r_k}$, and
	$$\bar{\cB}_k^{(t)} \in \mathbb{R}^{p_k\times \hat r_{k}}, \quad \bar{\cB}^{(t)}_{k, [i, :]} = \bar{\cA}^{(t)}_{k, [i, :]} 1_{\left\{\|\cA^{(t)}_{k, [i, :]}\|_2^2 \geq \bar{\eta}_k  \text{ and } i\in I_k\right\}}, \quad 1\leq i \leq p_k,$$
	where $\bar{\eta}_k = \hat{\sigma}^2\left(\hat{r}_k + 2(\hat r_k\log p)^{1/2} + 2\log p\right)$. 
	Finally, apply QR decomposition to $\bar{\cB}_k^{(t)}$, and assign the $Q$ part to $\hat{\cU}_k^{(t)}\in \mathbb{O}_{p_k, \hat r_k}$.
\end{enumerate}
Since $(A_7)$ is implied by the following two conditions,
\begin{equation}\label{ineq:new-final-coincide-to-prove-1}
\begin{split}
\forall 1\leq k \leq d, \forall i\notin I_k, \quad & \left\|(Y_k)_{[i, :]}\right\|_2^2 < \hat\sigma^2 \left(p_{-k} + 2\sqrt{p_{-k}\log p} + 2 \log p\right)\\
& \quad \text{and  } \max_{j} |(Y_k)_{[i,j]}| < 2\hat\sigma\sqrt{\log p};
\end{split}
\end{equation}
\begin{equation}\label{ineq:new-final-coincide-to-prove-2}
\begin{split}
\forall 1\leq k \leq d, \forall i\in I_k^{(0)}, \quad  \max_{j} |(Y_k)_{[i,j]}| \geq \hat\sigma\sqrt{\log p};
\end{split}
\end{equation}
Similar to Step 6 of the proof of Theorem \ref{th:upper_bound}, in order to show $(A_7)$ holds with high probability, it suffices to show \eqref{ineq:new-final-coincide-to-prove-1} and \eqref{ineq:new-final-coincide-to-prove-2} hold with probability at least $1 - O( (p_1+\ldots + p_d)/p^{1-\delta})$ when conditioning on fixed $\Z_{[I_1,\ldots,I_d]}$. By Lemmas \ref{lm:chi-square-tail-bound} and \ref{lm:tail-probability}, $\forall 1\leq k \leq d, i\notin I_k$, $\left\|(Y_k)_{[i, :]}\right\|_2^2 < \left(p_{-k} + 2\sqrt{(1-\delta)p_{-k}\log p} + 2(1-\delta) \log p\right)$ and $\max_{j} |(Y_k)_{[i,j]}| < \sqrt{(4-2\delta)\log p}$ hold with probability $1-O((p_1+\ldots+p_d)/p^{1-\delta})$. Thus for large constant $C>0$ and any $p \geq C$, we have
\begin{equation*}
	\begin{split}
		\left\|(Y_k)_{[i, :]}\right\|_2^2 < & \frac{\hat{\sigma}^2}{\left(1 - c\sqrt{(\log p)/p}\right)}\left(p_{-k} + 2\sqrt{(1-\delta)p_{-k}\log p} + 2(1-\delta) \log p\right) \\
		< & \hat \sigma^2 (p_{-k} + 2\sqrt{1-\delta}\sqrt{p_{-k}\log p} + 2(1-\delta)\log p) \\
		& + c \sqrt{\frac{p_{-k} \log p}{p_k}} + 2c\sqrt{1-\delta}\cdot \frac{\log p}{\sqrt p_k} + o\left(\frac{\log p}{\sqrt{p}}\right) \\
		\leq & \hat \sigma^2 \left(p_{-k} + 2\sqrt{p_{-k}\log p} + 2\log p\right),
	\end{split}
\end{equation*}
\begin{equation*}
		\max_{j} |(Y_k)_{[i,j]}|  < \sqrt{(4-2\delta)\log p}  \leq \frac{\hat{\sigma}}{\left(1 - c\sqrt{\frac{\log p}{p}}\right)}\cdot \sqrt{(4-2\delta)\log p} \leq 2\hat{\sigma}\sqrt{\log p}
\end{equation*}
happen with probability at least $1 - O((p_1+\cdots p_d)/p^{1-\delta})$. Moreover, similarly as the proof of \eqref{ineq:final-coincide-to-prove-2}, we have for any $1\leq k \leq d$, $i \in I_k^{(0)}$, if $p\geq C$ for some large constant $C>0$,
	\begin{equation*}
	\begin{split}
	\|X_{k, [i, :]}\|_2^2 \geq 25 s_{-k}\log p \quad \Rightarrow & \quad \max_{j\in I_{-k}} \left(X_{k, [i, j]}\right)^2 \geq \frac{C_0 s_{-k}\log p}{s_{-k}} = 25\log p\\
	\Rightarrow & \quad \max_{j\in I_{-k}} \left|X_{k, [i, j]}\right| \geq 5 \sqrt{\log p}\\
	\overset{\eqref{ineq:A-1}}{\Rightarrow} & \quad \max_{j\in I_{-k}} \left|Y_{k, [i,j]}\right| \geq 5\sqrt{\log p} - 2\sqrt{\log p} \\
	\Rightarrow & \quad \max_{j\in I_{-k}} \left|Y_{k, [i,j]}\right| \geq 3\left(\hat \sigma - c\sqrt{\frac{\log p}{p}}\right)\sqrt{\log p} \geq 2\hat \sigma \sqrt{\log p}.
	\end{split}
	\end{equation*}
	Thus, we have shown that $(A_6)$ holds with probability at least $1 - O((p_1+\cdots p_d)/p^{1-\delta})$ conditioning on fixed $\Z_{[I_1,\ldots, I_d]}$. This implies that $(A_1)-(A_7)$ happen with probability at least $1 - O(p_1+\cdots+p_d)/p^{1-\delta}$ (note that this is even true if $p\leq C$ for the large constant $C>0$), which has finished the proof of this proposition. \quad $\square$

\subsection{Proof of Theorem \ref{th:emperical-upper-bound}.}

Again, we assume $\sigma = 1$. We inherit the notations and assumptions from the proof of Proposition \ref{prop:rank_estimate}. We added two additional assumptions:
\begin{enumerate}[label=(\subscript{A}{\arabic*})]
		\setcounter{enumi}{7}
		\item
		\begin{equation}\label{new-ineq:A-8}
			\hat r_k = r_k, \quad \forall k = 1,\ldots,d.
		\end{equation}
		\item
		\begin{equation}\label{new-ineq:A-9}
			\begin{split}
				& \hat{I}_k^{(t)} \subseteq I_k,\quad t=1,\ldots, t_{\max},\\
&  \hat{J}_k^{(t)} \subseteq I_k,\quad t=1,\ldots, t_{\max}.
			\end{split}
		\end{equation}
	\end{enumerate}
As long as the assumptions$(A_1)-(A_9)$ are established, we can directly establish the error bound similarly as we did in the Steps 2-5 in the proof of Theorem \ref{th:upper_bound}. By Proposition \ref{prop:rank_estimate}, $(A_8)$ holds with high probability. To show $(A_9)$ happens with high probability, it suffices to show the following condition holds with probability at least $1 - O\left(t_{\max}(p_1+\cdots+p_d)/p\right)$ conditioning on fixed $\Z_{[I_1,\ldots,I_d]}$:
\begin{equation}\label{ineq:new-final-coincide-to-prove-3}
\begin{split}
\forall 1\leq k \leq d, 1\leq t \leq t_{\max}, i\notin I_k, \quad & \left\|\cA_{k, [i, :]}^{(t)}\right\|_2^2 < \hat \sigma^2 \left(r_{-k} + (\sqrt{r_{-k}\log p} + \log p)\right);\\
\end{split}
\end{equation}
\begin{equation}\label{ineq:new-final-coincide-to-prove-4}
\begin{split}
\forall 1\leq k \leq d, 1\leq t \leq t_{\max}, i\notin I_k, \quad & \left\|\bar{\cA}_{k, [i, :]}^{(t)}\right\|_2^2 < \hat \sigma^2 \left(r_k + (\sqrt{r_{k}\log p} + \log p)\right).
\end{split}
\end{equation}

For \eqref{ineq:new-final-coincide-to-prove-3}, we know that $\left\|\cA_{k, [i, :]}^{(t)}\right\|_2^2 < r_{-k} + 2\sqrt{(1-\delta)r_{-k}\log p} + 2(1-\delta)\log p)$ holds with probability at least $1-O(t_{max}(p_1+\cdots+p_d)/p^{1-\delta})$, thus \eqref{ineq:new-final-coincide-to-prove-3} follows from the following inequality if $p\geq C$ for large constant $C>0$,
	\begin{equation*}
		\begin{split}
			& \hat \sigma^2 \left(r_{-k} + 2\sqrt{r_{-k}\log p} + 2\log p\right) - \left(r_{-k} + 2\sqrt{(1-\delta)r_{-k}\log p} + 2(1-\delta)\log p\right) \\
			& \geq \left(1-c_1\sqrt{\frac{\log p}{p}}\right) \left(r_{-k} + 2\sqrt{r_{-k}\log p} + 2\log p\right) - \left(r_{-k} + 2\sqrt{(1-\delta)r_{-k}\log p} + 2(1-\delta)\log p\right) \\
			& \geq 2(1-\sqrt{1-\delta})\sqrt{r_{-k} \log p} + 2\delta\log p - c_1\left(r_{-k}\sqrt{\frac{\log p}{p}} +  \frac{r_{-k}\log p}{p} + o\left(\frac{\log p}{\sqrt{p}}\right)\right) \geq 0.
		\end{split}
	\end{equation*}
The last inequality is due to the assumption that $r \leq s = o(p)$. The proof of \eqref{ineq:new-final-coincide-to-prove-4} is essentially the same, we omit it here. \quad $\square$

\section{Technical Lemmas}

We collect the technical lemmas which will be used in the proof of the main results.

\begin{Lemma}[Properties of Orthogonal Matrices]\label{lm:property-orthogonal-matrix} ~
	\begin{itemize}
		\item Suppose $U\in \mathbb{O}_{p, m}, V\in \mathbb{O}_{m, r}$, recall $U_{\perp}\in \mathbb{O}_{p, p-m}, V_{\perp}\in \mathbb{O}_{m, m-r}$ are the orthogonal complements, then
		$$\left(UV\right)_{\perp} = \left[U_{\perp} ~~ UV_{\perp}\right].$$
		\item Suppose $U_1\in \mathbb{O}_{p_1,r_1}, \ldots, U_d\in \mathbb{O}_{p_d, r_d}$ are orthogonal matrices (not necessarily of the same dimension), then
		\begin{equation}\label{eq:outer-product-perpendicular}
			\begin{split}
				& \left(U_1\otimes \cdots \otimes U_d\right)_{\perp} \\
				= & \left[U_{1\perp}\otimes I_{p_2}\otimes \cdots \otimes I_{d}  \quad U_1\otimes U_{2\perp} \otimes I_3\otimes \cdots \otimes I_d \quad \cdots \quad  U_1\otimes \cdots \otimes U_{d-1}\otimes U_{d\perp}\right].
			\end{split}
		\end{equation}
		\item Suppose $\X \in \mathbb{R}^{p_1\times \cdots \times p_d}$ is a tensor, $X_k = \mathcal{M}_k(\X)$,
		$$\left\|X_k \left(U_{k+1}\otimes\cdots \otimes U_d\otimes U_1 \otimes \cdots \otimes U_d\right)_{\perp} \right\|_F^2\leq \sum_{l=1, l\neq k}^d \|U_{l\perp}^{\top}X_l\|_F^2. $$
	\end{itemize}
\end{Lemma}

{\noindent\bf Proof of Lemma \ref{lm:property-orthogonal-matrix}.}

\begin{itemize}
	\item Since
	\begin{equation*}
		\begin{split}
			& \left[UV ~~ UV_{\perp} ~~ U_{\perp}\right]^\top \cdot \left[UV ~~ UV_{\perp} ~~ U_{\perp}\right] =
			\begin{bmatrix}
			V^\top U^\top UV & O & O\\
			O & V^\top_{\perp} U^\top UV_{\perp} & O\\
			O & O & U_{\perp}^\top U_{\perp}
			\end{bmatrix}
			= I_p,
		\end{split}
	\end{equation*}
	
	$[UV_{\perp} ~~ U_{\perp}]$ is the orthogonal complement of $UV$.
	
	\item We only need to show $(U_1 \otimes U_2)_\perp = [U_{1\perp} \otimes I_{p2} ~~ U_1 \otimes U_{2\perp}]$, and the result follows by induction. Since $(A\otimes B)^\top = A^\top \otimes B^\top$ and $(A\otimes B)\cdot(C \otimes D) = (AC)\otimes (BD)$,  we could easily verify the following:
	\begin{equation*}
		\begin{split}
			& \left[U_1 \otimes U_2 ~~ U_{1\perp} \otimes I_{p2} ~~ U_1 \otimes U_{2\perp}\right]^\top \cdot \left[U_1 \otimes U_2 ~~ U_{1\perp} \otimes I_{p2} ~~ U_1 \otimes U_{2\perp}\right]
			= I_{p_1p_2}.
		\end{split}
	\end{equation*}
	
	\item Without loss of generality, we only need to prove the situation when $k=1$. Based on \eqref{eq:outer-product-perpendicular}, we have
	\begin{equation*}
		\begin{split}
			& X_1 \left(U_2\otimes\cdots \otimes U_d\right)_{\perp} \\
			= & \left[X_1\left(U_{2\perp}\otimes I_{p_3}\otimes \cdots \otimes I_{d}\right)  \quad X_1\left(U_2\otimes U_{3\perp} \otimes I_3\otimes \cdots \otimes I_d\right) \quad \cdots \quad X_1\left(U_2\otimes \cdots \otimes U_{d-1}\otimes U_{d\perp}\right)\right].
		\end{split}
	\end{equation*}
	Thus,
	\begin{equation*}
		\begin{split}
			& \left\|X_1 \left(U_2\otimes\cdots \otimes U_d\right)_{\perp}\right\|_F^2 = \sum_{l=2}^d \left\|X_1\left(U_{p_2}\otimes \cdots\otimes U_{l-1}\otimes U_{l\perp}\otimes I_{p_{l+1}}\otimes \cdots \otimes I_{p_{d}}\right) \right\|^2_F\\
			\leq & \sum_{l=2}^d\left\|X_1\left(I_{p_2}\otimes \cdots\otimes I_{l-1}\otimes U_{l\perp}\otimes I_{p_{l+1}}\otimes \cdots \otimes I_{p_{d}}\right)\right\|^2_F\\
			= & \sum_{l=2}^d\left\|\X \times_l U_{l\perp}\right\|_F^2 = \sum_{l=2}^d\left\|U_{l\perp}^\top X_l\right\|_F^2.
		\end{split}
	\end{equation*}
\end{itemize}
\quad $\square$

\begin{Lemma}[Projections and Sine-Theta distances]\label{lm:projection-perturbation-sine-theta}~
	\begin{itemize}
		\item Suppose $X \in \mathbb{R}^{m\times n}$ is any matrix, then for any $V, \hat{V}\in \mathbb{O}_{n, r}$, we have
		\begin{equation*}
			\left\|X\hat{V}\right\| \leq \left\|XV\right\|+ \|X\|\cdot \left\|\sin\Theta\left(\hat{V}, V\right)\right\|.
		\end{equation*}
		\item Suppose $X\in \mathbb{R}^{m\times n}$ is any matrix, $\rank(X)=r$, the left singular vectors of $X$ is $U\in \mathbb{O}_{m, r}$. Recall $\hat{U}_\perp\in\mathbb{O}_{m, m-r}$ is the orthogonal complement of $\hat{U}$, then for any $\hat{U}\in \mathbb{O}_{m, r}$,
		\begin{equation*}
			\begin{split}
				& \left\|(\hat{U}_\perp)^\top X\right\| \geq \sigma_r(X) \left\|\sin\Theta\left(\hat{U}, U\right)\right\|;\\
				& \left\|(\hat{U}_\perp)^\top X\right\|_F \geq \sigma_r(X) \left\|\sin\Theta\left(\hat{U}, U\right)\right\|_F.
			\end{split}
		\end{equation*}
	\end{itemize}
\end{Lemma}

{\bf\noindent Proof of Lemma \ref{lm:projection-perturbation-sine-theta}.}
\begin{itemize}
	\item Since $I = P_{V} + P_{V_\perp}$,
	\begin{equation*}
		\begin{split}
			\left\|X\hat{V}\right\| = & \left\|X\left(P_V + P_{V_\perp}\right)\hat{V} \right\| \leq \left\|XP_V\hat{V} \right\| + \left\|XP_{V_\perp}\hat{V} \right\|\\
			= & \left\|XVV^\top\hat{V}\right\| + \left\|XV_{\perp}\left(V_{\perp}^\top \hat{V}\right)\right\|\\
			\leq & \left\|XV\right\| + \|X\|\cdot \left\|V_{\perp}^\top \hat{V}\right\|.
		\end{split}
	\end{equation*}
	By Lemma 1 in \cite{cai2016rate}, $\left\|V_{\perp}^\top \hat{V}\right\| = \left\|\sin\Theta\left(\hat{V}, V\right)\right\|$, which has finished the proof of the first part of this lemma.
	\item Since the left singular vectors of $X$ is $U\in \mathbb{O}_{m, r}$, the projection satisfies $P_U X = UU^\top X = X$, thus,
	\begin{equation*}
		\begin{split}
			\left\|\hat{U}_\perp^\top X\right\| = \left\|\hat{U}_\perp^\top UU^\top X\right\| \geq \left\|\hat{U}_\perp^\top U\right\| \sigma_{\min}(U^\top X) = \left\|\sin\Theta\left(\hat{U}, U\right)\right\| \sigma_{r}(X);
		\end{split}
	\end{equation*}
	\begin{equation*}
		\begin{split}
			\left\|\hat{U}_\perp^\top X\right\|_F = \left\|\hat{U}_\perp^\top UU^\top X\right\|_F \geq \left\|\hat{U}_\perp^\top U\right\|_F \sigma_{\min}(U^\top X) = \left\|\sin\Theta\left(\hat{U}, U\right)\right\|_F \sigma_{r}(X).
		\end{split}
	\end{equation*}
	
\end{itemize}
\quad $\square$

\begin{Lemma}\label{lm:random-projection-lemma}
	Suppose $\Z \in \mathbb{R}^{s_1\cdots s_d}$ is an order-$d$ tensor with i.i.d. standard normal entries. Let $Z_k = \mathcal{M}_k(\Z)$ be the matricizations for $k=1,\ldots, d$, then
	\begin{equation*}
		\sup_{\substack{V_k \in \mathbb{R}^{s_k \times r_k}, \\ \|V_k\|\leq 1, 1\leq k \leq d}} \left\|Z_{k,[I_k,I_{-k}]} \cdot \left(V_{k+1} \otimes \cdots \otimes V_{d} \otimes V_1\otimes \cdots \otimes V_{k-1}\right)\right\| \leq C\left(\sqrt{s_k}+\sqrt{r_{-k}} + \sqrt{1+t}\sum_{l\neq k}\sqrt{s_lr_l}\right)
	\end{equation*}
	with probability at least $1 - C\exp\left(-Ct\sum_{l\neq k} s_lr_l\right)$.
\end{Lemma}

{\bf\noindent Proof of Lemma \ref{lm:random-projection-lemma}.} Without loss of generality we let $k=1$. By Lemma 7 in \cite{zhang2017tensor}, for each $m=2,\ldots,d$, there exists $\varepsilon$-net $\{V_m^{(1)},\ldots, V_m^{(N_m)}\}$ for $\{V_m\in \mathbb R^{s_m \times r_m}: \left\|V_m\right\| \leq 1\}$ with $\left|N_m\right| \leq ((4+\varepsilon)/\varepsilon)^{s_kr_k}$.\\
Consider the (d-1) dimensional index set $\mathcal I = [1:N_2] \times [1:N_3] \times \ldots \times [1:N_d]$, for a fixed index $i = (i_2,i_3,\ldots,i_d)\in \mathcal I$, consider the following matrix:
$$
Z_{1}^{(i)} = Z_{1,[I_1,I_{-1}]} \cdot \left(V_{2}^{(i_2)} \otimes \cdots \otimes V_{d}^{(i_d)} \right)
$$
Since $Z_{1}^{(i)}$ has independent rows, and each row follows a joint Gaussian distribution:
$$
N\left(0,\left(V_2^{(i_2)\top}V_2^{(i_2)}\right) \otimes \ldots \otimes \left(V_d^{(i_d)\top}V_d^{(i_d)}\right) \right)
$$
Note that $\left\| \left(V_2^{(i_2)\top}V_2^{(i_2)}\right) \otimes \ldots \otimes \left(V_d^{(i_d)\top}V_d^{(i_d)}\right) \right\| \leq 1$, by random matrix theory, we have:
\begin{equation*}
	\begin{split}
		P\left(\left\|Z_1^{(i)}\right\| \leq \sqrt{s_1} + \sqrt{r_{-1}} + x \right) \geq 1-2\exp(-cx^2/2).
	\end{split}
\end{equation*}
Then we further have:
\begin{equation*}
	\begin{split}
		& P\left(\max_{i \in \mathcal I}{\left\|Z_1^{(i)}\right\|} \leq \sqrt{s_1} + \sqrt{r_{-1}} + x \right) \\
		\geq & 1-2\left|\mathcal I\right|\exp(-cx^2/2) = 1-2((4+\varepsilon)/\varepsilon)^{\sum_{l=2}^{d}s_lr_l}\exp(-cx^2/2).
	\end{split}
\end{equation*}
Now, let
\begin{equation*}
	\begin{split}
		(V_2^*, V_3^*, \cdots, V_d^*) = \argmax_{\substack{V_m \in \mathbb{R}^{s_m \times r_m}, \\ \|V_m\|\leq 1, 2\leq m \leq d}} \left\|Z_{1,[I_1,I_{-1}]} \cdot \left(V_2 \otimes \cdots \otimes V_{d} \right)\right\|, \\
		M = \max_{\substack{V_m \in \mathbb{R}^{s_m \times r_m}, \\ \|V_m\|\leq 1, 2\leq m \leq d}} \left\|Z_{1,[I_1,I_{-1}]} \cdot \left(V_2 \otimes \cdots \otimes V_{d} \right)\right\|.
		\end{split}
\end{equation*}
Then by the definition of $\varepsilon$-net, we can find a index $i = (i_2,i_3,\cdots,i_d)$, such that $\left\|V_m^{i_m} - V_m^*\right\| \leq \varepsilon$ for any $2 \leq m \leq d$. Provided $\varepsilon \leq 1$, we have:
\begin{equation*}
	\begin{split}
		M & =  \left\|Z_{1,[I_1,I_{-1}]} \cdot \left(V_2^* \otimes \cdots \otimes V_{d}^* \right)\right\|	 = \left\|Z_{1,[I_1,I_{-1}]} \cdot \left((V_2^*-V_2^{(i_2)}+V_2^{(i_2)}) \otimes \cdots \otimes (V_{d}^* - V_d^{(i_d)} + V_d^{(i_d)}) \right)\right\| \\
		& \leq \left\| Z_{1,[I_1,I_{-1}]} \cdot \left(V_2^{(i_2)} \otimes \cdots \otimes V_{d}^{(i_d)}\right) \right\| + (2^{d-1}-1)\varepsilon M \\
		& \leq \sqrt{s_1}+\sqrt{r_{-1}} + x + (2^{d-1}-1)\varepsilon M
	\end{split}
\end{equation*}
with probability at least $1-2((4+\varepsilon)/\varepsilon)^{\sum_{l=2}^{d}s_lr_l}\exp(-cx^2/2)$, thus we have:
\begin{equation*}
	\begin{split}
		P\left(M \leq \frac{\sqrt{s_1}+\sqrt{r_{-1}}+x}{1-(2^{d-1}-1)\varepsilon}\right) \geq 1-2((4+\varepsilon)/\varepsilon)^{\sum_{l=2}^{d}s_lr_l}\exp(-cx^2/2).
	\end{split}
\end{equation*}
We set $\varepsilon = \frac{1}{2^d-2}$, $x^2 = C\sum_{l=2}^d s_lr_l(1+t)$ with sufficient big constant $C$, then we obtain the result for $k=1$, the bound for other modes follows similarly.

\begin{Lemma}[Probability Tail Bound for Non-central $\chi^2$ Distribution]\label{lm:chi-square-tail-bound} If $W$ satisfies the non-central Chi-square distribution $\chi_m^2(\lambda)$ with degrees of freedom $m$ and non-centrality parameter $\lambda$, so that $W = \sum_{i=1}^m X_i^2$, where $X_i\sim N(\mu_i, 1)$ with $\sum_i \mu_i^2 = \lambda$. Then for any $x>0$,
	\begin{equation*}
		P\left\{X \geq m + \lambda + 2 \sqrt{(m+2\lambda)x} + 2x\right\} \leq e^{-x},
	\end{equation*}
	\begin{equation*}
		P\left\{X \leq m + \lambda - 2\sqrt{(m+2\lambda)x} \right\} \leq e^{-x}.
	\end{equation*}
\end{Lemma}
The proof of this lemma is provided in Lemma 8.1 in \cite{birge2001alternative}, we thus omit it here.

\begin{Lemma}[Probability Tail Bound for Gaussian Extreme Values]\label{lm:tail-probability} If $u \in \mathbb{R}^p$, $u\overset{iid}{\sim}N(0, 1)$, then for any $x>0$, $p\geq 2$,
	\begin{equation*}
		P\left(\|u\|_\infty = \max_i |u_i| \geq \sqrt{x\log p} \right) \leq \sqrt{\frac{2}{\pi \log p}}p^{-x/2+1}.
	\end{equation*}
\end{Lemma}

{\noindent\bf Proof of Lemma \ref{lm:tail-probability}.} By the tail bound probability for Gaussian random variables,
\begin{equation*}
	P\left(|u_i| \geq \sqrt{x\log p}\right) \leq \frac{2e^{-(\sqrt{x\log p})^2/2}}{\sqrt{2\pi x\log p}} = \sqrt{\frac{2}{\pi \log p}} p^{-x/2}.
\end{equation*}
Thus,
\begin{equation*}
	P\left(\|u\|_\infty = \max_i |u_i| \geq \sqrt{x\log p}\right) \leq \sqrt{\frac{2}{\pi \log p}} p^{-x/2+1},
\end{equation*}
which has finished the proof for this lemma. \quad $\square$

\begin{Lemma}[Properties related to low-rank matrix perturbation]\label{lm:projection-X-residual}~
	\begin{itemize}
		\item Suppose $X, Z\in \mathbb{R}^{m\times n}$ and $Y=X+Z$, $\rank(X) = r$. If the leading $r$ left and right singular vector of $Y$ are $\hat{U}\in \mathbb{O}_{m, r}$ and $\hat{V}_{n, r}$, then
		\begin{equation*}
			\begin{split}
				& \max\left\{\left\|\hat{U}_\perp^\top X\right\|, \left\|X \hat{V}_{\perp}\right\|\right\} \leq 2\|Z\|,\\
				& \max\left\{\left\|\hat{U}_\perp^\top X\right\|_F, \left\|X \hat{V}_{\perp}\right\|_F\right\} \leq \min\left\{2\sqrt{r}\|Z\|, 2\|Z\|_F\right\}.
			\end{split}
		\end{equation*}
		\item Suppose $X, Z\in\mathbb{R}^{m\times n}$ and $Y=X+Z$, $\rank(X)\leq r$. Then
		\begin{equation*}
			\sigma_{r+1} \left(Y\right) \leq \|Z\|, \quad \left(\sum_{i=r+1}^{m\wedge n}\sigma_i^2(Y)\right)^{1/2} \leq \|Z\|_F.
		\end{equation*}
		\item Suppose $Y = X + Z$, $\rank(X)\leq r$, then
		\begin{equation*}
			\|X\|_F \leq \|Y\|_F + \sqrt{r}\|Z\|.
		\end{equation*}
	\end{itemize}
\end{Lemma}
{\noindent\bf Proof of Lemma \ref{lm:projection-X-residual}.} The proof is essentially the same as Lemma 7 in \cite{zhang2017tensor}. For completeness of the presentation we provide the proof here.

\begin{equation*}
	\begin{split}
		\left\|\hat{U}_{\perp}^\top X\right\| \leq & \left\| \hat{U}_{\perp}^\top (X+Z)\right\| + \|Z\| = \sigma_{r+1}(Y) + \|Z\| = \min_{\substack{\tilde{X}\in \mathbb{R}^{m\times n}\\\rank(\tilde{X})\leq r}} \|Y - \tilde{X}\| + \|Z\| \\
		\leq & \|Y-X\| + \|Z\| =  2\|Z\|.
	\end{split}
\end{equation*}
Since $\rank\left(\hat{U}_\perp^\top X \right) \leq  \rank(X) \leq r$, it is clear that
\begin{equation*}
	\left\|\hat{U}_\perp^\top X\right\|_{\rm F} \leq \sqrt{r}\left\|\hat{U}_\perp^\top X\right\|  \leq 2\sqrt{r}\|Z\|;
\end{equation*}
meanwhile,
\begin{equation*}
	\begin{split}
		\left\|\hat{U}_\perp^\top X\right\|_{\rm F} \leq & \left\|\hat{U}_\perp^\top (X+Z)\right\|_{\rm F} + \|Z\|_{\rm F} = \left(\sum_{i=r+1}^{p_1\wedge p_2}\sigma_{i}^2(Y)\right)^{1/2} + \|Z\|_{\rm F}\\
		\leq & \min_{\substack{\tilde{X}\in \mathbb{R}^{p_1\times p_2}\\\rank(\tilde{X})\leq r}} \|Y - \tilde{X}\|_{\rm F} + \|Z\|_{\rm F} \leq \|Y-X\|_{\rm F} + \|Z\|_{\rm F} \leq 2\|Z\|_{\rm F}.
	\end{split}
\end{equation*}

Furthermore,
\begin{equation*}
	\sigma_{r+1}(Y) = \min_{\rank(\tilde{X})\leq r} \|Y-\tilde{X}\| \leq \|Y-X\| = \|Z\|;
\end{equation*}
\begin{equation*}
	\left(\sum_{i=r+1}^{m\wedge n} \sigma_i^2(Y)\right)^{1/2} \leq \min_{\rank(\tilde{X})\leq r} \|Y-\tilde{X}\|_F \leq \|Y-X\|_F = \|Z\|_F,
\end{equation*}

Finally,
\begin{equation*}
	\begin{split}
		\|X\|_F = \|P_X X\|_F \leq \|P_X Y \|_F + \|P_X Z\|_F \leq \|Y\|_F + \sqrt{r}\|P_{X}Z\| \leq \|Y\|_F + \sqrt{r}\|Z\|,
	\end{split}
\end{equation*}
which has proved this lemma. \quad $\square$

\begin{Lemma}\label{lm:key-perturbation-result}Following the notations from the proof of Theorem \ref{th:upper_bound}, we have
	\begin{equation*}
		\left\|Z_{1, [I_1, :]} \hat{U}_{-1} \hat{V}_1\right\| \leq N_1 + M_1\left(\frac{E_2}{\lambda_2} + \cdots + \frac{E_d}{\lambda_d} + K_1\right).
	\end{equation*}
\end{Lemma}

{\noindent\bf Proof of Lemma \ref{lm:key-perturbation-result}.} By Lemma \ref{lm:property-orthogonal-matrix},
\begin{equation*}
	\begin{split}
		& \left(U_{-1}V_1\right)_{\perp} = \left[U_{-1}V_{1\perp} ~~ U_{-1\perp}\right]\\
		= & \left[U_{-1} V_{1\perp} \quad U_{2\perp} \otimes I \otimes \cdots \otimes I \quad U_2\otimes U_{3\perp} \otimes I \otimes \cdots \otimes I \quad \cdots \quad U_2\otimes \cdots \otimes U_{d-1}\otimes U_{d\perp}\right].
	\end{split}
\end{equation*}
Thus, we have the following decomposition.
\begin{equation*}
	\begin{split}
		Z_{1, [I_1, :]}\hat{U}_{-1} \hat{V}_1 = & Z_{1, [I_1, :]}P_{U_{-1}V_1} \hat{U}_{-1}\hat{V}_1 + Z_{1, [I_1, :]}P_{U_{-1}V_{1\perp}}\hat{U}_{-1}\hat{V}_1 + Z_{1, [I_1, :]}P_{U_{2\perp} \otimes I \otimes \cdots \otimes I} \hat{U}_{-1}\hat{V}_1\\
		& + \cdots  + Z_{1, [I_1, :]}P_{U_2\otimes \cdots \otimes U_{d-1}\otimes U_{d\perp}} \hat{U}_{-1}\hat{V}_1
	\end{split}
\end{equation*}
We analyze the spectral norm of the terms in the equation above respectively as follows.
\begin{equation*}
	\left\|Z_{1, [I_1, :]}P_{U_{-1}V_1} \hat{U}_{-1}\hat{V}\right\| = \left\|Z_{1, [I_1, :]}U_{-1}V_1 (U_{-1}V_1)^\top \hat{U}_{-1}\hat{V}\right\| \leq \|Z_{1, [I_1, :]}U_{-1}V_1\| \leq N_1;
\end{equation*}
\begin{equation*}
	\begin{split}
		&\left\|Z_{1, [I_1, :]}P_{U_{-1}V_{1\perp}}\hat{U}_{-1}\hat{V}_1 \right\| = \left\|Z_{1, [I_1, :]}U_{-1}V_{1\perp} \left(U_{-1}V_{1\perp}\right)^\top \hat{U}_{-1}\hat{V}_1 \right\| \leq \left\|Z_{1, [I_1, :]}U_{-1}\right\| \cdot \left\|\left(U_{-1}V_{1\perp}\right)^\top\hat{U}_{-1}\hat{V}_1\right\|\\
		\leq & M_1 \cdot \left\|\left(U_{-1}V_1\right)_{\perp}^\top \hat{U}_{-1}\hat{V}_1\right\| \leq M_1 \cdot \left\|\sin\Theta\left(\hat{U}_{-1}\hat{V}_1, U_{-1}V_1\right)\right\| = M_1K_1;
	\end{split}
\end{equation*}
\begin{equation*}
	\begin{split}
		& \left\|Z_{1, [I_1, :]} P_{U_{2\perp} \otimes I \otimes \cdots \otimes I} \hat{U}_{-1}\hat{V}_1\right\| = \left\|Z_{1, [I_1, :]} \left(U_{2\perp} \otimes I \otimes \cdots \otimes I\right)\left(U_{2\perp} \otimes I \otimes \cdots \otimes I\right)^\top\hat{U}_{-1}\hat{V}_1\right\|\\
		= & \left\|Z_{1, [I_1, :]} \left((U_{2\perp}U_{2\perp}^\top \hat{U}_2) \otimes \hat{U}_3\otimes \cdots \otimes \hat{U}_d\right)\hat{V}_1\right\| \leq \left\|Z_{1, [I_1, :]} \left((U_{2\perp}U_{2\perp}^\top \hat{U}_2) \otimes \hat{U}_3\otimes \cdots \otimes \hat{U}_d\right)\right\|\\
		\leq & \|U_{2\perp}^\top \hat{U}_2\|\cdot M_1  = M_1 \left\|\sin\Theta\left(U_2,\hat{U}_2\right)\right\| \leq \frac{M_1 E_2}{\lambda_2};
	\end{split}
\end{equation*}
Similarly,
\begin{equation*}
	\left\|Z_{1, [I_1, :]} P_{U_2\otimes \cdots \otimes U_{d-1}\otimes U_{d\perp}} \hat{U}_{-1}\hat{V}_1\right\| \leq \frac{M_1 E_d}{\lambda_d}.
\end{equation*}
To sum up, we have
\begin{equation*}
	\begin{split}
		\left\|Z_{1, [I_1, :]}\hat{U}_{-1} \hat{V}_1\right\| \leq & \left\|Z_{1, [I_1, :]}P_{U_{-1}V_1} \hat{U}_{-1}\hat{V}_1\right\| + \left\|Z_{1, [I_1, :]}P_{U_{-1}V_{1\perp}}\hat{U}_{-1}\hat{V}_1\right\| + \left\|Z_{1, [I_1, :]}P_{U_{2\perp} \otimes I \otimes \cdots \otimes I} \hat{U}_{-1}\hat{V}_1\right\|\\
		& + \cdots  + \left\|Z_{1, [I_1, :]}P_{U_2\otimes \cdots \otimes U_{d-1}\otimes U_{d\perp}}\right\|\\
		\leq & N_1 + M_1\left(\frac{E_2}{\lambda_2} + \cdots + \frac{E_d}{\lambda_d} + K_1\right),
	\end{split}
\end{equation*}
which has finished the proof for this lemma.\quad $\square$

\begin{Lemma}\label{lm:rate-analysis}Assume the following hold for some $t$:
	\begin{equation*}
		\begin{split}
			E_k^{(t-1)} &\leq 30\sqrt{s_kr_k} + 30\sqrt{s_k\log p} + \frac{12\sqrt{ds\log p}}{2^{(t-1)}},\quad k=1,\ldots, d,\\
			M_1 & \leq 2\left(\sqrt{s_1}+\sqrt{r_{-1}}+\sum_{l=2}^{d}\sqrt{s_lr_l}+\sqrt{\log p}\right)\\
		\end{split}
	\end{equation*}
	If we set $C_{gap}>20(d+2)(d-1)$, then we have:
	\begin{equation}\label{rate-1}
		\begin{split}
			\frac{\sqrt{r_{1}}M_1E_k^{(t-1)}}{\lambda_j} \leq \frac{1}{d-1}\left(3\sqrt{s_1r_1}+3\sqrt{s_1\log p}+\frac{12\sqrt{ds\log p}}{2^{t+2}}\right),\quad 2\leq k\leq d,\quad 1\leq j\leq d.
		\end{split}
	\end{equation}
	\begin{equation}\label{rate-2}
		\begin{split}
			\frac{M_1\sqrt{2s_1r_1\eta_1}}{\lambda_1} \leq \frac{1}{5}\sqrt{s_1r_1}+\frac{2}{5}\sqrt{s_1\log p}
		\end{split}
	\end{equation}
	\begin{equation}\label{rate-3}
		\begin{split}
			\frac{\sqrt{6}r_1M_1^2}{\lambda_1} \leq \frac{3}{5}\sqrt{s_1r_1}+\frac{3}{5}\sqrt{s_1\log p}
		\end{split}
	\end{equation}
\end{Lemma}

{\noindent Proof of Lemma \ref{lm:rate-analysis}.} Since
\begin{equation}\label{M/lambda}
	\begin{split}
		\frac{M_1}{\lambda_j} \leq \frac{\sqrt{r_1}M_1}{\lambda_j} \leq \frac{2\left(\sqrt{s_1r_1}+r_{-1}+\sum_{l\geq2}\frac{r_1+s_lr_l}{2}+\sqrt{r_1\log p}\right)}{\lambda_j} \leq \frac{2(d+2)\lambda_j}{C_{gap}\lambda_j} \leq \frac{1}{10(d-1)},
	\end{split}
\end{equation}
\eqref{rate-1} essentially follows. Note that $M_1\sqrt{2s_1r_1\eta_1} \leq \sqrt{2}M_1\sqrt{s_1r_1r_{-1}}+2M_1\sqrt{s_1r_1\log p}$, and
\begin{equation*}
	\begin{split}
		\sqrt{2}M_1\sqrt{s_1r_1r_{-1}} & \leq 2\sqrt{2}\left(\sqrt{s_1r_1}\sqrt{s_1r_{-1}}+\sqrt{s_1r_1}r_{-1}+\sum_{l\geq2}\sqrt{s_1r_1}\sqrt{r_{-1}s_lr_l}+\sqrt{s_1\log p}\sqrt{r_{1}r_{-1}}\right)\\
		& \leq 2\sqrt{2}\left(\sqrt{s_1r_1}\sqrt{s}+\sqrt{s_1r_1}r_{-1}+\sum_{l\geq2}\sqrt{s_1r_1}\frac{r_{-1}+s_lr_l}{2}+\sqrt{s_1\log p}\frac{r_1+r_{-1}}{2}\right)\\
		& \leq 2\sqrt{2}\left((d+1)\sqrt{s_1r_1}+\sqrt{s_1\log p}\right)\frac{\lambda_1}{C_{gap}} \leq \frac{\lambda_1}{5}\left(\sqrt{s_1r_1}+\sqrt{s_1\log p}\right).
	\end{split}
\end{equation*}
Thus, we have:
\begin{equation*}
	\begin{split}
		\frac{M_1\sqrt{2s_1r_1\eta_1}}{\lambda_1} & \leq  \frac{1}{5}\left(\sqrt{s_1r_1}+\sqrt{s_1\log p}\right) +  \frac{2M_1\sqrt{r_1}\sqrt{s_1\log p}}{\lambda_1} \\
		& \overset{\eqref{M/lambda}}{\leq} \frac{1}{5}\left(\sqrt{s_1r_1}+\sqrt{s_1\log p}\right) + \frac{1}{5}\sqrt{s_1\log p}\\
		& = \frac{1}{5}\sqrt{s_1r_1} + \frac{2}{5}\sqrt{s_1\log p}.
	\end{split}
\end{equation*}
which proves \eqref{rate-2}. Now we turn to $\sqrt{6}r_1M_1^2$,
\begin{equation*}
	\begin{split}
		\sqrt{6}r_1M_1^2 &\leq 4\sqrt{6}r_1\left(s_1+r_{-1}+\sum_{l\geq2}s_lr_l+\log p\right) \\
		& \leq 4\sqrt{6}\left(\sqrt{s_1r_1}\sqrt{s}+\sqrt{s_1r_1}r_{-1}+\sqrt{s_1r_1}\sum_{l\geq2}s_lr_l+\sqrt{s_1\log p}\sqrt{s\log p}\right) \\
		& \leq 4\sqrt{6}\left((d+1)\sqrt{s_1r_1}+\sqrt{s_1\log p}\right)\frac{\lambda_1}{C_{gap}} \leq \frac{3\lambda_1}{5}\left(\sqrt{s_1r_1}+\sqrt{s_1\log p}\right)
	\end{split}
\end{equation*}
Then \eqref{rate-3} essentially follows.$\quad \square$